\documentclass[3p,times,english,sort&compress]{elsarticle}
\usepackage[T1]{fontenc}
\usepackage[utf8]{inputenc}
\usepackage{babel}
\usepackage{float}
\usepackage{textcomp}
\usepackage{graphicx}
\usepackage{amsmath}
\usepackage{esint}
\usepackage{colortbl}
\usepackage{booktabs}
\usepackage{multirow}
\usepackage{setspace}
\usepackage{amsthm}
\usepackage{amssymb}
\usepackage{latexsym}
\usepackage{color} 
\usepackage{xfrac}
\usepackage{algorithm}
\usepackage{algpseudocode}
\usepackage{verbatim}

\usepackage[final]{changes}

\everymath{
    \fontdimen14\textfont2=1.0ex
    \fontdimen17\textfont2=0.3ex
}

\newcolumntype{L}{>{\centering\arraybackslash}m{2.5cm}}
\definecolor{aliceblue}{rgb}{0.94,0.87,0.8}
\definecolor{anti-flashwhite}{rgb}{0.95, 0.95, 0.96}
\definecolor{antiquewhite}{rgb}{0.7,0.75,0.71}

\usepackage[unicode=true,
 bookmarks=true,bookmarksnumbered=false,bookmarksopen=false,
 breaklinks=false,pdfborder={0 0 1},backref=false,colorlinks=false]
 {hyperref}
\hypersetup{
    unicode=true,
    colorlinks,
    linkcolor={dark-red},
    citecolor={dark-blue},
    urlcolor={medium-blue},
}
\makeatletter

\pdfpageheight\paperheight
\pdfpagewidth\paperwidth

\@ifundefined{showcaptionsetup}{}{%
\PassOptionsToPackage{caption=false}{subfig}}
\usepackage{subfig}
\makeatother

\begin{document}

\title{A Hybrid \textcolor{black}{Finite-Volume} Reconstruction Framework for Efficient High-Order Shock-Capturing on Unstructured Meshes}

\author[CU]{Yiren Tong}
\ead{yiren.tong.502@cranfield.ac.uk}
\author[CU]{Panagiotis Tsoutsanis\corref{cor1}}
\ead{panagiotis.tsoutsanis@cranfield.ac.uk}

\cortext[cor1]{Corresponding author}

\address[CU]{Faculty of Engineering and Applied Sciences, Cranfield University, Cranfield MK43 0AL, United Kingdom}

\begin{abstract}
In this paper, we present a multi-dimensional, arbitrary-order hybrid reconstruction framework for compressible flows on unstructured meshes. The proposed method advances state-of-the-art high-resolution schemes by combining the efficiency of linear reconstruction with the robustness of high-order non-oscillatory formulations, activated only where necessary through a novel a priori detection strategy. This approach minimises the use of costly Compact Weighted Essentially Non-Oscillatory (CWENOZ) or Monotonic Upstream-centered Scheme for Conservation Laws (MUSCL) reconstructions, thereby substantially reducing computational overhead without compromising accuracy or stability.
The framework integrates the strengths of CWENOZ formulations and the Multi-dimensional Optimal Order Detection (MOOD) paradigm, while introducing a redesigned Numerical Admissibility Detector (NAD) that classifies the local flow field in a single step into smooth, weakly non-smooth, and discontinuous regions. Each region is then reconstructed using an optimal method: a high-order linear scheme in smooth areas, CWENOZ in weakly non-smooth zones, and a second-order MUSCL scheme near discontinuities. This targeted, a priori allocation preserves high-order accuracy where possible and guarantees non-oscillatory, stable solutions near shocks and strong gradients.
The proposed hybrid strategy is implemented within the open-source unstructured finite-volume solver UCNS3D and supports arbitrary-order reconstructions on mixed-element meshes. Comprehensive two- and three-dimensional benchmark tests demonstrate that the method maintains the designed order of accuracy in smooth regions while significantly enhancing robustness in shock-dominated flows. Owing to the reduced frequency of expensive nonlinear reconstructions, the framework achieves up to a $2.5\times$ speed-up compared to a CWENOZ scheme of the same order in 3D compressible turbulence simulations. Overall, this hybrid framework brings high-order accuracy closer to in industrial-scale CFD simulations through its combination of reduced computational cost, improved robustness, and reliability.

\end{abstract}
\maketitle

\section{Introduction}
Compressible flow has long been a challenging topic within computational fluid dynamics. In practical aero engineering applications, an accurate prediction of compressible phenomena is crucial for design optimisation, flow control, and safety certification. However, the presence of strong gradients and discontinuities in such flows renders their numerical simulation particularly demanding. Lower-order numerical methods, while robust and capable of effectively suppressing spurious oscillations induced by discontinuities, often suffer from insufficient accuracy. This limitation primarily arises due to inherent numerical dissipation, which prevents these methods from accurately capturing small-scale flow structures in the smooth flow regions. Moreover, modern engineering problems frequently involve complex geometries. The geometric flexibility offered by unstructured meshes is essential for accurately representing complicated geometries efficiently, and reducing the time for grid-generation, in particular when several iterations of the design will be needed in an industrial setting. However, adapting numerical schemes to these unstructured meshes introduces additional complexities. To simultaneously capture fine-scale flow structures and prevent spurious oscillations, numerous high-order non-oscillatory schemes have been proposed. Considerable work has also been done to extend these schemes to unstructured grid frameworks.

\textcolor{black}{{One class of approaches enriches the local solution representation within each element, exemplified by the discontinuous Galerkin (DG) \citep{DUMBSER20088209,ZHU20084330,ZHU2009293321,XU20092194,MALTSEV2023111755,DG_FV_MULTI} and Flux Reconstruction (FR) frameworks \citep{Castonguay201210915,WILLIAMS201353,VINCENT2015248}, where high-order accuracy is obtained by evolving element-local polynomial degrees of freedom. A second class of approaches retains a finite-volume (FV) conservation update but enhances the reconstruction strategy to achieve high order on unstructured grids. Representative examples include the spectral finite-volume (SFV) method \citep{WANG2002210,MENGALDO201556,XU20095787}, which introduces subcell control volumes within each element, and compact finite-volume \citep{PONT201745,WANG20171,XIE2022105436} that target high-order accuracy using compact face-neighbour stencils. Another important direction is the multi-moment FV method \citep{XIE2017637,DENG201785,DENG2019404,XIE2021109841}, which augments the traditional volume-integrated average (VIA) with additional moments such as point values (PV) at selected locations to enable high-order reconstructions.}}

Among these high-order non-oscillatory approaches, one of the earliest and influential development is the Monotone Upstream-centered Schemes for Conservation Laws (MUSCL), The MUSCL scheme was first introduced by Bram van Leer \citep{VANLEER1974361,VANLEER1977263,VANLEER1977276,VANLEER1979101}, which extends Godunov's first-order non-oscillatory method to second-order spatial accuracy by employing a piecewise-linear reconstruction with a slope limiter. Many different slope limiters have been designed after that; for example, the widely used second-order limiter for unstructured grids proposed by Barth and Jespersen \citep{Barth1989366}, and Venkatakrishnan's limiter \citep{VENKATAKRISHNAN1995120}, which improves the Barth-Jespersen limiter by making it differentiable. In addition to these second-order limiters, numerous studies have focused on high-order limiters \citep{MICHALAK20098693,LI20124053,LIU201788,TSOUTSANIS201869}. 

Following previously published work on MUSCL schemes, Ami Harten \citep{HARTEN1983357} introduced the idea of total variation diminishing (TVD), which states that if the total variation of the numerical solution never increases from one time level to the next, then the numerical scheme is TVD. It should be noted that a MUSCL reconstruction can be made fully compatible with Harten's TVD framework whenever an admissible limiter \citep{SWEBY1984215} is chosen. However, although TVD schemes guarantee the elimination of spurious oscillations, the solution degenerates to first-order accuracy at smooth extrema, which strictly restricts the accuracy of TVD schemes. To overcome this limitation, Chi-Wang Shu \citep{Shu198710521} proposed the total variation bounded (TVB) framework, which relaxes the monotonicity requirement by introducing a mesh dependent threshold that preserves high-order accuracy at non-sonic critical points. Around the same time, Harten et al. \citep{HARTEN1987231} addressed the same issue from another angle with the essentially non-oscillatory (ENO) scheme, which avoids discontinuities by adaptively choosing the smoothest stencil instead of relying on slope limiting. 

Subsequently, Liu et al. \citep{LIU1994200} introduced the original Weighted Essentially Non-Oscillatory (WENO) scheme. Building on that work, Jiang and Shu \citep{JIANG1996202} incorporated a novel local smoothness indicator, establishing the classic fifth-order WENO-JS scheme. Although WENO-JS delivers markedly higher accuracy than the original WENO method, several deficiencies persist. The method experiences order degeneration near non-sonic critical points, and it also exhibits excessive numerical dissipation in moderately smooth flow regions. Moreover, the resulting solution is highly sensitive to the parameter used during the weight calculations. 

Over the following decades, several important works have been conducted to address these shortcomings. The first kind of approach involves directly improving the performance of WENO-JS scheme. Henrick et al. \citep{HENRICK2005542} introduced a mapped weights strategy, successfully resolving the issue of order degeneration near smooth extrema, although at the expense of increased computational cost. Borges et al. \citep{BORGES20083191} proposed a global smoothness indicator that reduces dissipation and restores accuracy near critical points without increasing the computational cost compared to the original WENO-JS scheme. Levy et al. \citep{LevyPuppoRusso1999} first introduced the Central WENO scheme, which achieves a high-order accurate reconstruction at the cell center while maintaining compact stencils. Nevertheless, the classical WENO-JS framework still encounters challenges related to robustness and computational efficiency. This has led to the development of an expanding family of WENO variants, which aim to enhance accuracy, robustness, and efficiency in demanding flow problems \citep{GEROLYMOS20098481,TITAREV2004238,TAYLOR2007384,MARTIN2006270,FU2016333,HU20108952,TSOUTSANIS20111585,TSOUTSANIS2021104961,CiCP252}.

Besides enhancements to the WENO scheme itself, another popular strategy is to develop hybrid methods. Because WENO reconstructions are computationally expensive and exhibit poor resolution qualities for shock-turbulence problems, using a \textcolor{black}{{high-order}} linear scheme in smooth regions can significantly improve overall computational efficiency and accuracy. In practice, this approach typically employs a shock detector to distinguish between regions that require a non-linear scheme owing to the presence of discontinuities and those that are smooth, thereby realising a hybrid approach. Adams and Shariff \citep{ADAMS199627} introduced a hybrid compact ENO method, which achieves high spectral resolution in turbulent regions while retaining a sharp, non-oscillatory character in the vicinity of discontinuities. Pirozzoli \citep{PIROZZOLI200281} derived a hybrid compact-WENO scheme that employs the same strategy and  compact stencils  to improve the resolution properties and efficiency of the algorithm. Since then, several hybrid WENO schemes have been developed \citep{REN2003365,COSTA2007970,ZHAO2019422,Fidalgo2018335,Lin20191007,2010JCoPh2291213J}, but the performance of these hybrid methods is often strongly dependent on the shock detector \citep{Lin20191007,LI20108105} employed. Indicators that use low-order information are computationally inexpensive yet highly case-sensitive, whereas those that use high-order information tend to exhibit increased computational cost.

In addition, some modern numerical methods employ adaptive criteria to adjust the solution's order of accuracy when discontinuities or strong gradients arise. The most noteworthy of these is the MOOD (Multi-dimensional Optimal Order Detection) paradigm, introduced by Clain et al. \citep{CLAIN20114028}. Unlike most high-order methods using a priori limitation of reconstructed values, MOOD paradigm operates a posteriori limitation: the high-order update is first accepted as a candidate solution; if this candidate violates a discrete maximum principle (DMP) or positivity check, the algorithm recomputes only the troubled cell with a lower-order, more dissipative scheme. The order is progressively deducted until the modified update passes all admissibility tests. The MOOD paradigm has subsequently undergone a series of improvements \citep{DIOT201243,tsoutsanis2022posteriori,FARMAKIS2020112921,TSOUTSANIS2023127544}, demonstrating excellent robustness and efficiency in a range of difficult flow configurations. Although the MOOD algorithm generally exhibits greater robustness than hybrid WENO formulations, the MOOD paradigm still has some shortcomings. First, its computational overhead increases quickly with the percentage of cells that are flagged as troubled. In some shock-dominated flows, the iterative reduction of spatial order triggers repeated local re-computations and, as a result, erodes efficiency. Second, the discrete maximum principle admissibility test (DMP) used in MOOD remains too restrictive. Even the relaxed variants proposed in later studies still impose more numerical dissipation than that found in well-tuned hybrid WENO schemes.

In this paper, inspired by the high-order hybrid WENO schemes and the MOOD paradigm, we propose a multi-dimensional arbitrary-order hybrid WENO strategy within a finite volume method framework on unstructured meshes. Our goal is to design a method that combines the high efficiency of the hybrid WENO scheme with the robustness provided by the MOOD method. We therefore follow the hierarchical structure of the MOOD paradigm but, for increased efficiency, introduce a redesigned DMP-based indicator that allows us to complete the hierarchical process in a single step. All schemes and paradigms utilised in this work have been implemented within the open-source solver UCNS3D \citep{ANTONIADIS2022108453}. This work evaluates the performance characteristics of the proposed method through a series of challenging 2-D and 3-D numerical tests.

The paper is organised as follows: in Section \ref{sec:framework} we introduce the numerical framework used to describe the high-order finite-volume method on unstructured meshes, detailing the reconstruction procedures for Linear, MUSCL, and CWENOZ, as well as the hybrid approach employed. The numerical results for all test cases are presented in Section \ref{sec:results} and are compared with analytical, reference, or experimental solutions whenever possible. Finally, the conclusions drawn from this study are summarised in the last section.

\section{Numerical Framework}\label{sec:framework}
We consider the 3D compressible Navier-Stokes equations, which can be written in their conservative
form:
\begin{equation}
\frac{\partial\mathbf{U}(\mathbf{x},{t})}{\partial t}+\nabla \cdot ( \vec{\mathbf{F}}_{c}(\mathbf{U})-\vec{\mathbf{F}}_{v}(\mathbf{U,\nabla{U}}))=0, \label{eq:NAVIERCOMPLETE}
\end{equation}
 where $\mathbf{U}$ is the vector of conserved variables, $\vec{\mathbf{F}}_{c}$ and $\vec{\mathbf{F}}_{v}$ are the vectors of inviscid and viscous flux, respectively, as
 \begin{equation}
\textbf{U}= \begin{bmatrix}
\rho \\ \rho u \\ \rho v \\\rho w \\ E 
\end{bmatrix}, 
\vec{\mathbf{F}}_{c}= \begin{bmatrix}
\rho u_n \\ \rho u u_n +n_x p \\ \rho v u_n +n_y p \\ \rho w u_n +n_z p \\ u_n(E+p) 
\end{bmatrix},
\vec{\mathbf{F}}_{v}= \begin{bmatrix}
0 \\ n_x \tau_{xx}+n_y \tau_{xy}+n_z \tau_{xz} \\ n_x \tau_{yx}+n_y \tau_{yy}+n_z \tau_{yz} \\ n_x \tau_{zx}+n_y \tau_{zy}+n_z \tau_{zz} \\ n_x \Theta_{x}+n_y \Theta_{y}+n_z \Theta_{z}
\end{bmatrix},
\label{Terms1}
\end{equation}
where $\rho$ is the density; $u,v,w$ are the velocity components in $x,y$
and $z$ Cartesian coordinates respectively, and $u_n$ is the velocity normal to the bounded surface area, defined by $u_n=n_{x}u+n_{y}v+n_{z}w$. Perfect gas assumption is used, and the total energy per unit mass is calculated by  $E=p/\left(\gamma-1\right)+(1/2)\rho(u^{2}+v^{2}+w^{2})$, where $p$ is the pressure, $\gamma=1.4$ is the ratio of specific heats for air at normal atmospheric conditions; The Sutherland's law is used for the computation of viscosity as a function of temperature:
\begin{equation}
\frac{\mu_{l}}{\mu_{0}}=\left(\frac{T}{T_{0}}\right)^{\frac{3}{2}}\frac{T_{0}+S}{T+S},
\label{eq:sutherlandlaw}
\end{equation}
$S$ is the Sutherland temperature and the subscript ${0}$ denotes a reference state for the corresponding variables. 
The work of viscous stresses and heat conduction, $\Theta$, is given by:
\begin{equation}\begin{split}
&\Theta_{x}=u\tau_{xx}+v\tau_{xy}+w\tau_{xz}+\frac{\mu_{l}}{Pr}\frac{\gamma}{\left(\gamma-1\right)}\frac{\partial T}{\partial x},\\
&\Theta_{y}=u\tau_{yx}+v\tau_{yy}+w\tau_{yz}+\frac{\mu_{l}}{Pr}\frac{\gamma}{\left(\gamma-1\right)}\frac{\partial T}{\partial y},\\
&\Theta_{z}=u\tau_{zz}+v\tau_{zy}+w\tau_{zz}+\frac{\mu_{l}}{Pr}\frac{\gamma}{\left(\gamma-1\right)}\frac{\partial T}{\partial z}.\end{split}\label{eq:HeatConduct}\end{equation}
The viscous stress tensor $\tau_{ij}$  is defined by is 
\begin{equation}
\tau_{ij} = \mu_{l} \left(\frac{\partial{{\mathbf{u}}_{i}}}{\partial{{\mathbf{x}}_{j}}}+\frac{\partial{{\mathbf{u}}_{j}}}{\partial{{\mathbf{x}}_{i}}}-\frac{2}{3} \frac{\partial{{\mathbf{u}}_{k}}}{\partial{{\mathbf{x}}_{k}}}\delta_{ij}\right),
\label{eq:Stresses} \end{equation}
where $\delta_{ij}$ is the Kronecker delta and the subscripts $i,j,k$ refer to the Cartesian coordinate components ${\mathbf{x}} =(x,y,z)$. 

The domain $\omega$ consisting of any combination of tetrahedral, hexahedral, prism or pyramid elements in 3D. And quadrilateral or triangular elements in 2D. Each of these elements is uniquely identified by a single index $i$. By integrating Eq.(\ref{eq:NAVIERCOMPLETE}) over a mesh element and applying a high-order explicit finite-volume formulation, the following expression is obtained:

\begin{equation}\begin{split}
\frac{d{\mathbf{U}}_{i}}{dt}= & -\frac{1}{|V_{i}|}\sum\limits _{l=1}\limits^{N_{f}}\sum\limits _{\alpha=1}\limits^{N_{qp}}\vec{\mathbf{F}}_{c_{l}}\left({\mathbf{U}}^{n}_{l,L}(\mathbf{x}_{l,\alpha},t),{\mathbf{U}}^{n}_{l,R}(\mathbf{x}_{l,\alpha},t)\right)\omega_{\alpha}|S_{l}|\\ & +\frac{1}{|V_{i}|}\sum\limits _{l=1}\limits^{N_{f}}\sum\limits _{\alpha=1}\limits^{N_{qp}}\vec{\mathbf{F}}_{v_{l}}\left({\mathbf{U}}^{n}_{l,L}(\mathbf{x}_{l,\alpha},t),{\mathbf{U}}^{n}_{l,R}(\mathbf{x}_{l,\alpha},t),
\nabla{\mathbf{U}}^{n}_{l,L}(\mathbf{x}_{l,\alpha},t),\nabla{\mathbf{U}}^{n}_{l,R}(\mathbf{x}_{l,\alpha},t)\right)\omega_{\alpha}|S_{l}|,\end{split}
\label{eq:EulerEquationsDetails}
\end{equation}

Where $\mathbf{U}_i$ are the conserved variables averaged over the element volume, $V_i$ is the area/volume of the element, $N_f$ is the total number of faces/sides associated with the element, $N_{qp}$ is the number of quadrature points used to approximate the surface integrals, $\vec{\mathbf{F}}_{c_{l}}$ and $\vec{\mathbf{F}}_{v_{l}}$ is the numerical convective and diffusive flux respectively evaluated normal to the interface between element $i$ and its neighbouring element $j$, $ \mathbf{U}_{ij,L}^n(\mathbf{x}_\alpha, t) $ and $\mathbf{U}_{ij,R}^n(\mathbf{x}_\alpha, t) $ are the high-order approximations of the solution on the left (considered cell) and the right (in the neighboring cell) sides of the interface, and $\nabla{\mathbf{U}}^{n}_{l,L}(\mathbf{x}_{l,\alpha},t)$ and ${\nabla\mathbf{U}}^{n}_{l,R}(\mathbf{x}_{l,\alpha},t)$  the high-order approximations of their gradients respectively. \textcolor{black}{A suitable Gaussian quadrature rule is employed for each element type when evaluating the interface integrals. Gauss-Legendre quadrature is used for tensor-product elements, such as quadrilateral and hexahedral cells, while symmetric Gaussian quadrature rules are adopted for simplex elements, including triangular and tetrahedral cells \citep{WILLIAMS201418}}. $\alpha$ denotes individual Gaussian quadrature points $\mathbf{x}_\alpha$ and their associated weights $\omega_\alpha $ along the interface, and $|S_{l}|$ is the length/area of the shared interface. The inviscid Euler equations are obtained by setting the viscous fluxes $\vec{\mathbf{F}}_{v}=0$.

\subsection{Spatial Reconstruction}
The least squares reconstruction method used in this study follows the approach of Tsoutsanis et al. \citep{TSOUTSANIS2021104961,TSOUTSANIS20111585,TSOUTSANIS2014254} and Titarev et al. \citep{TITAREV2010}, which have shown success in addressing smooth and discontinuous flow problems of different Mach numbers \citep{TSOUTSANIS2018157, TSOUTSANIS2023127544, ASFSILVA2022107401, tsoutsanis2022posteriori, adebayo2022implementation, Silva2022unstructured, TSOUTSANIS20218964, SILVA2021106518, FARMAKIS2020112921, RICCI2020105648, ANTONIADIS201786, TSOUTSANIS2015207, fluids9020033, Tommaso2024dynamic,takis_adda,Silva2024}. Here, we summarise only the essential components, while referring the reader to the aforementioned studies for comprehensive implementation details.

For each considered element $i$, we  build a high-order polynomial $P_{i}(x,y,z)$ of arbitrary order r. The high-order polynomial $P_{i}(x,y,z)$ has the same average as a general quantity $\mathbf{U}_i$ and the order of accuracy is $r+1$, which can be written as:

\begin{equation}
\mathbf{U}_i \;=\; \frac{1}{\lvert V_i \rvert}\int_{V_i} P_i(x,y,z)\ dV.
\label{eq:AVGVariables1}
\end{equation}
The transformation used by Dumbser et al. \citep{DUMBSER2007204} and Tsoutsanis et al. \citep{DUMBSER2007204} is used in the present study to minimise scaling effects that can occur due to different sizes/shapes of unstructured elements in the stencil. The system will transform elements of different sizes from physical space $(x,y,z)$ to a reference space $(\xi,\eta,\zeta)$. Except from the triangular or tetrahedral elements, all other elements will be decomposed into them for the reconstruction process. The detail of element decomposition can be found in Tsoutsanis et al. \citep{TSOUTSANIS2014254}. The general quantity $\mathbf{U}_i$ will not change during the transformation.

\begin{equation}
\mathbf{U}_i \;=\; \frac{1}{\lvert V_i \rvert}\int_{V_i} \mathbf{U}(x,y,z)\ dV \equiv \frac{1}{\lvert V_i' \rvert}\int_{V_i'} \mathbf{U}(\xi,\eta,\zeta)\ d\xi d\eta d\zeta.
\label{eq:AVGVariables2}
\end{equation}
The $V_i'$ represents the volume/area of the considered element in reference space. To perform the reconstruction on the target element $S_i$, we recursively include neighboring elements, forming a stencil of $M+1$ cells that includes the considered cell $i$. The central (or linear) stencil is then defined as:
\begin{equation}
{\cal S}_{i}^{c}=\bigcup\limits _{m=0}^{M_{c}}V_{m},
\end{equation}
where the index $m$ indicates the local numbering of the elements in the stencil, with the element assigned index 0 corresponding to the considered cell $i$. \textcolor{black}{{The index $c$ referring to the stencil number (in case of multiple stencils).}} There are various stencil selection strategies; in this work, we adopt the stencil-based compact (SBC) algorithm proposed by Tsoutsanis \citep{TSOUTSANIS2019100037}. The number of cells in the stencil $\mathcal{S}$ is chosen as $M = 2k$ to enhance robustness of reconstruction \citep{DUMBSER2007204, FU201925, TSOUTSANIS20111585, TSOUTSANIS2019100037, DIOT201243, DUMBSER20091731, NOGUEIRA20102544}, where $k$ is the number of unknown coefficients in the polynomial approximation, which is written as:
\begin{equation}
K(r, d) = \frac{1}{d!} \prod_{l=1}^{d} (r + l),
\end{equation}
where $r$ is the order of reconstruction polynomial, and $d \in [2,3]$ denotes the space dimensions. The stencil is transferred to the reference space, becoming $S_{i}'$, and the $r^{th}$ order polynomial is expressed by locally expanding the basis functions as follows:
\begin{equation}
p(\xi, \eta, \zeta) = \sum_{k=0}^{K} a_k \phi_k (\xi, \eta, \zeta) 
= \mathbf{U}_0 + \sum_{k=1}^{K} a_k \phi_k (\xi, \eta, \zeta),
\end{equation}
where $a_k$ represents the unknown polynomial coefficients, also referred to as degrees of freedom, $\phi_k$ is the basis function, and the vector of conserved variables in the considered cell $i$ is denoted by $\mathbf{U}_0$. For each cell within the stencil, the degrees of freedom $a_k$ must satisfy the condition that the cell averaged value of the reconstruction polynomial $p(\xi, \eta, \zeta)$ equal to the cell averaged solution $\mathbf{U}_m$, which can be written as:

\begin{equation}
\int_{V_m'} p(\xi, \eta, \zeta) \, d\xi d\eta d\zeta 
= |V_m'| \mathbf{U}_0 + \sum_{k=1}^{K} \int_{V_m'} a_k \phi_k \, d\xi d\eta d\zeta 
= |V_m'| \mathbf{U}_m, \quad m = 1, \ldots, M.
\end{equation}
Since elements other than triangles in 2D and tetrahedra in 3D cannot be guaranteed to map to a unit reference element, basis functions $\phi_k$ should also satisfy equation (\ref{eq:AVGVariables1}):

\begin{equation}
\phi_k(\xi, \eta, \zeta) \equiv \psi_k(\xi, \eta, \zeta) - \frac{1}{|V_0'|} \int_{V_0'} \psi_k \, d\xi d\eta d\zeta, \quad k = 1, 2, \ldots, K.
\end{equation}
In the present study, the Legendre polynomials are employed as basis functions for $\psi_k$. We define the integral of the basis function over cell $m$ as $A_{mk}$, and introduce the right hand side vector, denoted by $b$, as follows:
\begin{equation}
A_{mk} \equiv \int_{V_m'} \phi_k \, d\xi d\eta d\zeta, 
\quad b_m \equiv |V_m'| (\mathbf{U}_m - \mathbf{U}_0).
\end{equation}
Then, the following matrix can be built, as:
\begin{equation}
\sum_{k=1}^{K} A_{mk} a_k = b_m, \quad m = 1, 2, \ldots, M.
\label{eq:dofMatrix}
\end{equation}
Since the matrix $A_{km}^TA_{mk}$ is invertible, the degrees of freedom $a_k$ can be written using the pseudo-inverse $A_{km}^\dagger$:
\begin{equation}
a_k  = A_{km}^\dagger b_m = \left( A_{km}^\mathrm{T} A_{mk} \right)^{-1} A_{km}^\mathrm{T} b_m.
\end{equation}
For solving the overdetermined linear system, QR decomposition based on Householder transformation \citep{Stewart1998} is used in present study. Since the geometry of the stencil elements will not change during the simulation, the pseudo-inverse matrix $A_{km}^\dagger$ is be only computed once at the beginning of the simulation.

\subsection{MUSCL}
In the current study, the degrees of freedom $a_k$ are determined through the least-squares reconstruction process. Based on the \textcolor{black}{{high-order}} linear scheme introduced above, the MUSCL scheme can be written as:
\begin{equation}
\mathbf{U}_{l,\alpha} = \mathbf{U}_i + \theta_i \cdot \sum_{k=1}^{K} a_k \phi_k (\xi_a, \eta_a, \zeta_a),
\end{equation}
Where $\mathbf{U}_{l,\alpha}$ represents the extrapolated reconstructed solution on the face $l$, the quadrature point $\alpha$, $\mathbf{U}i$ is the value of the conserved variable in element $i$,  and $\theta_i$ is the slope limiter. By selecting different types of limiters, we can impose prior constraints on the solution at each quadrature point to suppress oscillations. Finally, \textcolor{black}{{$\xi_\alpha, \eta_\alpha, \zeta_\alpha$ }}specifies the coordinates of the quadrature point located on the face $l$. In the current study, the second-order Barth \& Jespersen limiter \citep{Barth1989366} (minmod equivalent for unstructured meshes)  and the Venkatakrishnan limiter \citep{VENKATAKRISHNAN1995120} are evaluated to determine how various slope limiters affect our hybrid method. Furthermore, we also examine the high-order slope limiter MOGE \citep{TSOUTSANIS201869}, which is designed for MUSCL schemes with accuracy up to the fourth order.

\subsection{CWENOZ}
The CWENOZ scheme developed by Tsoutsanis and Dumbser \citep{TSOUTSANIS2021104961} improves both the computational efficiency and the robustness of the original WENO schemes and has enabled a more compact framework for hybridisation of DG-FV schemes \citep{DG_FV_MULTI,MALTSEV2023111755}. The key difference between the classic WENO scheme \citep{JIANG1996202} and the CWENOZ scheme can be summarised as follows: first, the CWENOZ scheme combines an optimal (high-order) polynomial $p_{opt}$ with lower-order polynomials. The high-order polynomial is used by the central stencil, while the lower-order polynomials are applied to the directional stencils. Employing these lower-order directional stencils, which are contained within the high-order central stencil, results in notable time savings compared to the original WENO scheme. Moreover, when discontinuities occur in the fluid domain, some of the lower-order directional stencils can contain smooth data. With smooth data present, the high-order optimal polynomial is restored to attain the desired order of accuracy. All polynomials involved are required to satisfy the previously stated condition of matching the cell averages of the solution. The definition of an optimal polynomial is given by:
\begin{equation}
p_{\text{opt}}(\xi, \eta, \zeta) = \sum_{s=1}^{s_t} \lambda_s p_s(\xi, \eta, \zeta),
\end{equation}
where $s$ denotes the stencil index, with $s = 1$ corresponding to the central stencil, and $s = 2, 3, \ldots, s_t$ representing the directional stencils. Here, $s_t$ is the total number of stencils, and $\lambda_s$ is the linear coefficient associated with each stencil, whose sum is equal to 1. The $p_1$ polynomial is not evaluated directly. It is obtained by subtracting the lower-order polynomials from the optimal polynomial $p_{opt}$, as follows:
\begin{equation}
p_1(\xi, \eta, \zeta) = \frac{1}{\lambda_1} \left( p_{\text{opt}}(\xi, \eta, \zeta) 
- \sum_{s=2}^{s_t} \lambda_s p_s(\xi, \eta, \zeta) \right).
\end{equation}
The CWENOZ reconstruction polynomial can be written as
\begin{equation}
p(\xi, \eta, \zeta)^{\text{cweno}} = \sum_{s=1}^{s_t} \omega_s p_s(\xi, \eta, \zeta).
\end{equation}
The reconstruction polynomial is a non-linear combination of all polynomials from $s = 1$ to $s = s_t$, where $\omega_s$ is the non-linear weight assigned to each polynomial, which also constitutes another main difference from the original WENO scheme. The original WENO non-linear weight, denoted by $\omega_s^{weno}$, is defined as:

\textcolor{black}{{
\begin{equation}
\omega_s^{\mathrm{weno}} =
\frac{\tilde{\omega}_s}{\sum\limits_{s=1}^{s_t} \tilde{\omega}_s}
\quad \text{where} \quad
\tilde{\omega}_s = \frac{\lambda_s}{\left(\epsilon + \mathcal{SI}_s\right)^b}.
\end{equation}
}}
In the CWENOZ scheme, the non-linear weights are modified by multiplying the linear weight $\lambda_s$ with a correction factor, and are defined as:
\textcolor{black}{{
\begin{equation}
\omega_s^{cwenoz} = 
\frac{\tilde{\omega}_s}{\sum\limits_{s=1}^{s_t} \tilde{\omega}_s}
\quad \text{where} \quad
\tilde{\omega}_s = \lambda_s \left( 1 + \frac{\tau}{\epsilon + \mathcal{SI}_s} \right),
\end{equation}
}}
where $\tau$ is the universal oscillation indicator, defined as the absolute difference between the smoothness indicators, given by:
\begin{equation}
\tau = \left( 
\frac{\displaystyle \sum_{s=2}^{s_t} \left|\mathcal{SI}_s - \mathcal{SI} _1 \right|}
{s_t - 1}
\right)^{b},
\end{equation}
and the smoothness indicator $\mathcal{SI_m}$ can be written as:
 \textcolor{black}{{
\begin{equation}
\mathcal{SI}_s = \sum_{1 \leq |\beta| \leq r} 
\int_{V_0'} |h_c|^{|\beta|-1} 
\left( \mathcal{D}^\beta p_s(\xi, \eta, \zeta) \right)^2 
\, (d\xi\, d\eta\, d\zeta),
\label{eq:SI1}
\end{equation}
}}
where $\beta$ is a multi-index, $r$ is the order of the polynomials, $h_c$ is the characteristic cell size of the considered cell, \textcolor{black}{{and $\mathcal{D}$ is the differential operator.}} In the present work, we set $\epsilon = 10^{-3}$ and $b = 4$. We also employ $r = 1$ for the directional polynomials, which yield second-order accuracy, and allow for any arbitrary order of accuracy in the polynomial associated with the central stencil. The linear weights are determined by first assigning an arbitrary non-normalised linear weight $\lambda_1'$ to the central stencil, followed by normalisation as defined by:
\begin{equation}
\lambda_1 = 1 - \frac{1}{\lambda_1'},
\end{equation}
with the linear weights associated with lower-order polynomials being assigned the same linear weights as follows:
\begin{equation}
\lambda_s = \frac{1 - \lambda_1}{s_t - 1}.
\end{equation}
Since the smoothness indicator is a quadratic function of the degrees of freedom $a_k^s$, Eq. (\ref{eq:SI1}) can also be written as:
\begin{equation}
\mathcal{SI}_s = \sum_{k=1}^{K} a_k^s 
\left( \sum_{q=1}^{K} \mathcal{O I}_{kq} a_q^s \right),
\end{equation}
where the oscillation indication matrix $\mathcal{O I}_{kq}$ is defined as:
\textcolor{black}{{
\begin{equation}
\mathcal{O} \mathcal{I}_{kq} = \sum_{1 \leq |\beta| \leq r} 
\int_{V_0'} |h_c|^{|\beta|-1} 
\left( \mathcal{D}^\beta \phi_k(\xi, \eta, \zeta) \right) 
\left( \mathcal{D}^\beta \phi_q(\xi, \eta, \zeta) \right) 
\, (d\xi\, d\eta\, d\zeta),
\end{equation}
}}
and can be precomputed and stored at the beginning of the simulation. This reconstruction can be performed with respect to the characteristic, primitive, or conserved variables and the reader is referred to \citep{TSOUTSANIS20218964,TSOUTSANIS20111585,TSOUTSANIS2014254,fluids9020033, TSOUTSANIS2023127544} and references therein for further details of the implementation.

\subsection{Fluxes approximation \& Temporal discretisation}
\textcolor{black}{It needs to be stressed that all the reconstructed solutions are subject to the Harten \citep{HARTEN1983357} bounds (referred to as the TVD check hereforth), which bounds the reconstructed solutions.
Let $\bar{q}_i$ denote the cell-average of a reconstructed variable $q$ in cell $i$, and
$q_{i,f}^{\text{rec}}$ the corresponding reconstructed value at face $f$ of cell $i$.
We apply an order-reduction (jump) fix such that
\begin{equation}
q_{i,f} \;=\;
\begin{cases}
q_{i,f}^{\text{rec}}, & \text{if } \left|q_{i,f}^{\text{rec}}-\bar{q}_i\right|
\le \theta\,\max\!\left(|\bar{q}_i|,\,\varepsilon\right), \\[6pt]
\bar{q}_i, & \text{otherwise},
\end{cases}
\label{eq:harten_jump_fix}
\end{equation}
where $\theta=0.9$ and $\varepsilon$ is a small threshold (e.g.\ $10^{-14}$) preventing
division-by-zero or overly strict limiting when $\bar{q}_i\approx 0$.
When the inequality is violated, the state for the entire cell is set to the cell average (first-order). The two variables that are checked are density and pressure, if any of them violates this inequality the reconstructed values for this cell for all the variables is set to the cell average.}

For inviscid fluxes we employ the (Harten-Lax-van Leer-Contact) HLLC approximate Riemann solver of Toro et al. \cite{HLLC}, or the HLL Riemann solver \citep{HLL}. 
For viscous fluxes, the gradients are computed using a constrained least-square reconstruction for the enforcement of the boundary condition, as detailed in \citep{TSOUTSANIS2014254,ANTONIADIS2022108453}. 
The gradients of the discontinuous states for the approximation of the viscous fluxes are averaged by including the penalty terms similar to previous approaches \cite{Nakahashi199997} in the following manner:
\begin{equation}
  \nabla{\mathbf{U}}=\frac{1}{2}\left(\nabla{\mathbf{U}_L}+\nabla{\mathbf{U}_R}\right)+\frac{\alpha}{L_{int}} \left({\mathbf{U}_R}-{\mathbf{U}_L}\right) \vec{n},
  \label{eq:viscous_jump}
 \end{equation}
where ${L_{int}}$ is the distance between the cell centres of adjacent cells, and we set $\alpha=4/3$ similarly to previous approaches \cite{gooch4}. 
For this study, we employ the explicit 4th-order Strong-Stability-Preserving (SSP) Runge-Kutta of Spiteri and Ruuth \cite{SSP4}, which is stable for $(CFL\approx 1.5)$, where the readers are referred for further details \cite{SSP4}.

\subsection{Hybrid Method}
In this section, we propose a hybrid reconstruction strategy that fuses the high-order hybrid WENO framework with some of the key ingredients of the a posteriori MOOD method \citep{DIOT201243,CLAIN20114028,tsoutsanis2022posteriori,FARMAKIS2020112921,TSOUTSANIS2023127544}. Rather than progressively lowering the polynomial degree of troubled cells, as the original MOOD algorithm does in a posteriori fashion, we deploy an a priori troubled-cell indicator that partitions the computational domain into three areas: smooth, weakly nonsmooth, and discontinuous. Each class is then assigned a dedicated reconstruction: a high-order linear polynomial in smooth regions, the high-order CWENOZ scheme in weakly non-smooth regions, and a second-order MUSCL scheme across discontinuities. This targeted allocation maintains optimal accuracy where the solution is smooth that avoids the expensive CWENOZ reconstruction, while enhancing robustness and non-oscillatory behaviour in the presence of sharp gradients or shocks.

\begin{figure}[h!]
\begin{centering}
\captionsetup[subfigure]{width=0.9\textwidth}
  \subfloat[Default setting of the Hybrid Method: $\alpha_m=5^{-3}$, $\beta_m=0.5$, $\alpha_w=0$, and $\beta_w=0$]
 {\includegraphics[angle=0,width=0.99\textwidth]{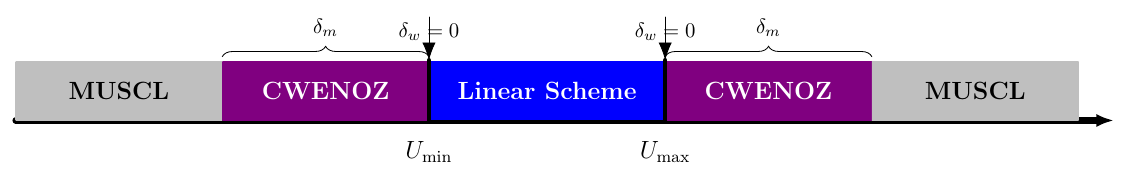}} \\
 \subfloat[Relaxed setting of the Hybrid Method: $\alpha_m=10^{-4}$, $\beta_m=10^{-1}$, $\alpha_w=10^{-4}$, and $\beta_w=10^{-3}$]
 {\includegraphics[angle=0,width=0.99\textwidth]{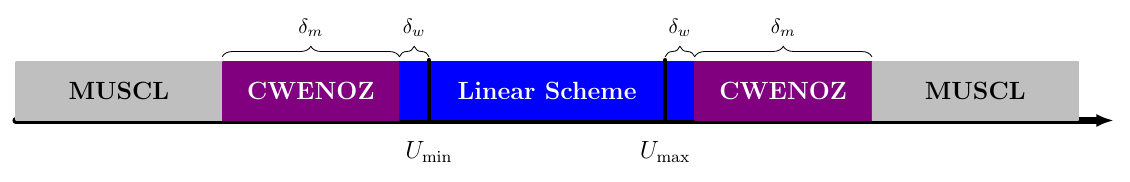}} \\
 \subfloat[Linear-CWENOZ setting of the Hybrid Method:  $\alpha_m=10^{6}$, $\beta_m=10^{9}$, $\alpha_w=0$, and $\beta_w=0$]
 {\includegraphics[angle=0,width=0.99\textwidth]{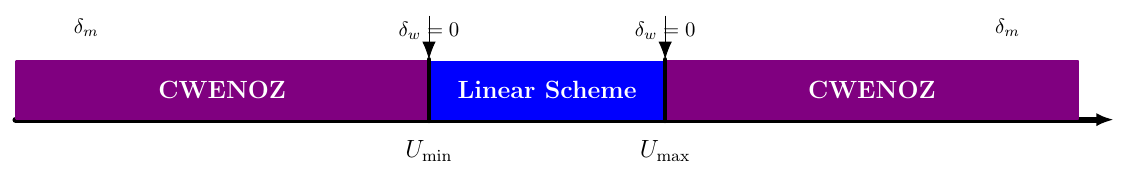}} \\
  \subfloat[Linear-MUSCL setting of the Hybrid Method:  $\alpha_m=0$, $\beta_m=0$, $\alpha_w=0$, and $\beta_w=0$]
 {\includegraphics[angle=0,width=0.99\textwidth]{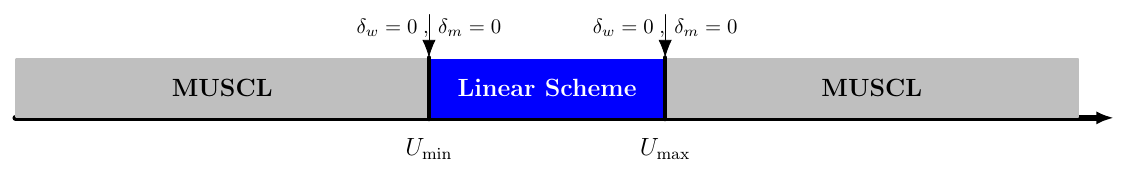}} \\
  \subfloat[MUSCL-CWENOZ setting of the Hybrid Method:  $\alpha_m=0$, $\beta_m=0$, $\alpha_w=0$, and $\beta_w=-0.5$]
 {\includegraphics[angle=0,width=0.99\textwidth]{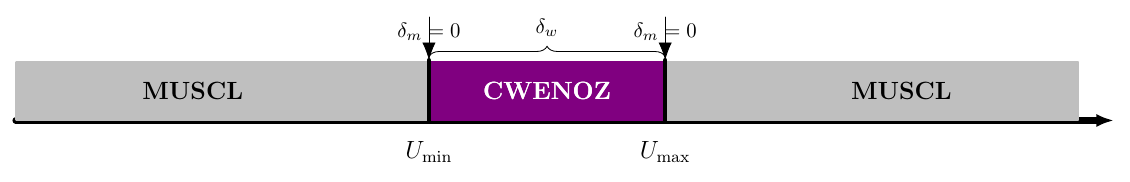}} \\
 
\par\end{centering}\caption{Illustration of some typical combinations of different NAD margins. By adjusting the margin parameters, we can easily control how we mix different schemes, and thereby achieve the designed order of accuracy and lower computational cost for different problems. The plots presented above may not be strictly to scale.}\label{fig:NADMargin}\end{figure}

Troubled cell indicators critically influence both the accuracy and robustness of hybrid reconstruction schemes. In this study, we employ a MOOD type Numerical Admissibility Detector (NAD) \citep{tsoutsanis2022posteriori, TSOUTSANIS2023127544}, chosen for its seamless integration into our hybrid reconstruction workflow and its low computational cost.  NAD follows the discrete maximum principle (DMP) \citep{CLAIN20114028,FARMAKIS2020112921}, DMP will ensure that the reconstructed values of the cells lie within the physically admissible limits inferred from its neighbors,  preventing the creation of new extrema. Since we divide the flow field into three regions, we modify this indicator accordingly. In its original form, the DMP criterion is expressed as
\begin{equation}\label{eq:DMP}
    \displaystyle{\min_{y\in \mathcal{V}_i}}(U^n(y)) \leq U(x) \leq \displaystyle{\max_{y\in \mathcal{V}_i}}(U^n(y)).
\end{equation}
\\ \textcolor{black}{{The set of cells $\mathcal{V}_i$ denotes the von Neumann neighbours of the cell $i$ selected by the algorithm.}} When our solution satisfies the DMP criterion, the solution can be considered monotonic. Accordingly, we can naturally use linear reconstruction in this fluid region. We also need two other criteria for the weakly non-smooth, and discontinuous regions. The NAD as a relaxed version of DMP, using an extra $\delta$ to relax the limitation bound. In our study, we introduce the auxiliary parameters $\delta_w$ and $\delta_m$ to control the region classification. When designing these two parameters, we aim for high flexibility that can be adjusted to different situations. $\delta_m$ determines the area of discontinuity, which can be defined as:
\begin{equation}\label{eq:NADMUSCL}
\begin{cases}
\displaystyle \min_{y \in \mathcal{V}i} U^n(y) - \delta_m \geq U(x), \\
\displaystyle \max_{y \in \mathcal{V}_i} U^n(y) + \delta_m \leq U(x),
\end{cases}
\quad \text{(Apply MUSCL scheme)}
\end{equation}
And together with $\delta_w$, these two parameters determine the area of the weakly non-smooth region:
\begin{equation}\label{eq:NADCWENOZ}
\begin{cases}
\displaystyle \min_{y \in \mathcal{V}i} U^n(y) - \delta_m \leq U(x) \leq \min_{y \in \mathcal{V}i} U^n(y) + \delta_w, \\
\displaystyle \max_{y \in \mathcal{V}i} U^n(y) + \delta_m \geq U(x) \geq \max_{y \in \mathcal{V}_i} U^n(y) - \delta_w,
\end{cases}
\quad \text{(Apply CWENOZ scheme)}
\end{equation}
The MUSCL margin $\delta_m$ and WENO margin $\delta_w$ can be expressed as:
\begin{equation}\label{eq:marginm}
\delta_m = \max\left(\alpha_m,\beta_m \cdot \left[\max_{y \in \mathcal{V}_i} (U^n(y)) - \min_{y \in \mathcal{V}_i} (U^n(y))\right]\right).
\end{equation}
\begin{equation}\label{eq:marginw}
\delta_w = \operatorname{sign}(\beta_w)\cdot\max\left(\alpha_w,\left|\beta_w \cdot \left[\max_{y \in \mathcal{V}_i} (U^n(y)) - \min_{y \in \mathcal{V}_i} (U^n(y))\right]\right|\right).
\end{equation}
The MUSCL margin design $\delta_m$ considers only positive values to define the MUSCL region, and for $\delta_w$ we allow both positive and negative values to control the CWENOZ reconstruction region. Adjusting the margin parameters $\alpha_w$, $\beta_w$, $\alpha_m$, and $\beta_m$, we can flexibly modify the regions for different reconstruction methods, and the algorithm can better handle different cases. In Fig. \ref{fig:NADMargin}, we present visualisations of different combinations of margin parameters. \textcolor{black}{{First, panel (a) presents the default parameter setting adopted in most of our tests. In this default configuration, we set $\alpha_m=5^{-3}$ and $\beta_m=0.5$ to delineate the MUSCL and CWENOZ subregions, while $\alpha_w$ and $\beta_w$ were set to 0. Relative to previous MOOD studies, these values correspond to more permissive (relaxed) margin parameters, which is enabled by our proposed three-region design. In the classical MOOD paradigm, a DMP-based indicator determines whether the scheme switches from a high-order WENO discretization to a more robust MUSCL discretization. In contrast, our framework introduces an additional WENO region to treat weakly non-smooth yet non-discontinuous structures. This intermediate layer mitigates spurious oscillations, thereby allowing larger margin parameters without compromising robustness.}}

 \textcolor{black}{{Panel (b) corresponds to a more aggressive bias toward the high-order linear scheme, with $\alpha_m=10^{-4}$, $\beta_m=10^{-1}$, $\alpha_w=10^{-4}$, and $\beta_w=10^{-3}$. By tuning $\alpha_w$ and $\beta_w$, the method can further expand the linear-scheme region and thus reduce numerical dissipation when higher spectral fidelity is required, for example in turbulence-dominated zones.}}

 \textcolor{black}{{Panels (c), (d), and (e) further illustrate the flexibility of the proposed hybrid strategy. By choosing very large margin parameters, $\alpha_m=10^{6}$ and $\beta_m=10^{9}$, the algorithm effectively restricts the hybridization to the high-order linear scheme and the CWENOZ scheme. Conversely, setting all margin parameters to zero, $\alpha_m=0$, $\beta_m=0$, $\alpha_w=0$, and $\beta_w=0$, yields a two-scheme formulation that employs only the linear and MUSCL schemes. Finally, selecting $\alpha_w=0$ and $\beta_w=-0.5$ suppresses the linear-scheme region, such that the hybrid method transitions only between the CWENOZ and MUSCL schemes.}}

 \textcolor{black}{{In this work, all test cases, except those explicitly identified for comparison, use the same hybrid parameters, which are designated as the default setting, which can be found in Fig. \ref{fig:NADMargin} (a).}}

In addition, to prevent the troubled cell indicator from erroneously activating in completely smooth regions, such as the fluid far field, and because this type of troubled cell indicator is sometimes too sensitive \citep{FARMAKIS2020112921,TSOUTSANIS2023127544}, we use a deactivation criterion similar to the Venkatakrishnan limiter \citep{VENKATAKRISHNAN1995120}. The additional safeguard avoids unnecessary limiting and preserves high-order accuracy where the flow remains smooth. The deactivation criterion can be written as: 

\begin{equation}\label{eq:NADDeactivatedCriteria}
\delta\mathbf{U} \;\equiv\; \delta\mathbf{U}_i^{\max} \;-\; \delta\mathbf{U}_i^{\min}
\;<\;
\bigl(\kappa\,\Delta x\bigr)^{n}.
\end{equation}

\textcolor{black}{{
Where $\kappa$ is a tunable parameter, and $\Delta x$ is also represents characteristic cell size of the considered cell. The maximum and minimum deviations of the solution between the current control volume and its immediate neighbours can be written as: }}

\textcolor{black}{{
\begin{equation}\label{eq:NADDeactivatedCriteria_2}
\delta \mathbf{U}^{\min}(x)=\min_{\substack{y\in V_x\\y\neq x}}\left\|\mathbf{U}^{n}(y)-\mathbf{U}^{n}(x)\right\|,
\qquad
\delta \mathbf{U}^{\max}(x)=\max_{\substack{y\in V_x\\y\neq x}}\left\|\mathbf{U}^{n}(y)-\mathbf{U}^{n}(x)\right\|.
\end{equation}
}}
Following the classification of different regions, the next step involves assigning a reconstruction method for each one of them. We will use a high-order linear scheme for the smooth-regions as the baseline. There are several reasons for this. Firstly, in smooth regions the high-order WENO family of methods offers no additional accuracy over the high-order linear method, yet the WENO scheme incurs much higher computational cost. Second, using a \textcolor{black}{{high-order}} linear scheme avoids the limitations inherent in the a priori criteria used by WENO type of schemes. For example, in steep but smooth regions, the traditional WENO-JS scheme will degenerate its order of accuracy. Although many improved variants of WENO schemes, such as the CWENOZ scheme we implemented, mitigate this issue, linear reconstructions still retain a clear advantage because of its low dissipation. We will not need to worry about spurious oscillations introduced by the \textcolor{black}{{high-order}} linear scheme. Our trouble-cell indicator detects these locations and reconstructs them with more dissipative methods, keeping the whole simulation robust.

Another important point is that when the solution computed from \textcolor{black}{{high-order}} linear scheme violates the DMP, we prefer the CWENOZ scheme instead of immediately switching to the MUSCL scheme. The primary consideration of our hybrid method is to apply the \textcolor{black}{{high-order}} linear scheme whenever the flow field is completely smooth. When the flow becomes turbulent and new extrema appear, the DMP type trouble cell indicator may fail to recognise these new extrema. However, because of the high dissipation inherent in the MUSCL scheme, applying it prematurely would erase these small structure details. Hence, the NAD detector used in the previous study \citep{TSOUTSANIS2023127544} introduces an extra margin to preserve these small structures whenever possible. Therefore, in our default hybrid method setting, we allocate a large operating region to the CWENOZ scheme. As mentioned earlier, the improved CWENOZ scheme possesses relatively low dissipation, enabling the accurate capture of small-scale turbulent structures. Moreover, due to the excellent performance of the CWENOZ scheme, we can be generous by expanding the region where the CWENOZ scheme is used, as this will not compromise the robustness of the overall method. We will only use the MUSCL scheme when there are drastic changes in the local solution, such as when a discontinuity occurs, to avoid the CWENOZ scheme producing non-physical solutions, thereby further enhancing the robustness of the method.

\begin{figure}[h!]
    \begin{centering}

    {\includegraphics[angle=0,width=0.90\textwidth]{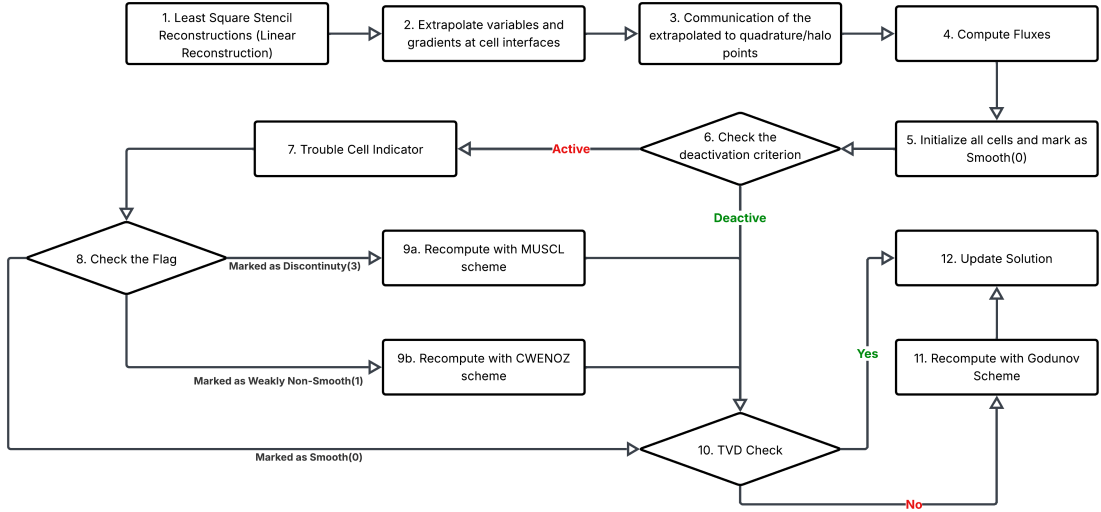}}
    \caption{The flowchart of the implementation of the Hybrid method}
    \label{fig:HybridWENOFlowchart}
    \end{centering}
\end{figure}

The flowchart of the hybrid method can be seen in Fig. \ref{fig:HybridWENOFlowchart}, and the overall steps can be described as follows: First, the high-order linear reconstruction is used to compute the candidate solution. Next, the candidate solution passes through the trouble cell indicator, which tags the flagged cells and their Von Neumann neighborhood as CWENO or MUSCL cells. Subsequently, the CWENOZ weights are calculated for those cells marked as CWENOZ cells, and we recompute the solution with CWENOZ reconstruction. After that, the cells tagged as MUSCL undergo MUSCL reconstruction and are computed again. Finally, all cells are checked again against the condition introduced by Harten et al. \citep{HARTEN1987231}, if any cell fails at this step, they are recalculated with the first-order Godunov scheme. Through the above process, we expect that our hybrid method will exhibit enhanced robustness compared with previous methods and will flexibly address different problems by adjusting the hybrid parameters $\delta_w$ and $\delta_m$, thus achieving significantly lower computational cost while maintaining similar or better accuracy than the CWENOZ scheme.

\section{Applications}\label{sec:results}
\subsection{2D Vortex Evolution}

\begin{figure}[h!]
    \begin{centering}
    \captionsetup[subfigure]{width=0.4\textwidth}

    \subfloat[Density contour of 2D Vortex Evolution on uniform triangular meshes ($ \mathcal{P}3 / 32 $)]
    {\includegraphics[angle=0,width=0.47\textwidth,trim={6cm 6cm 6cm 6cm},clip]{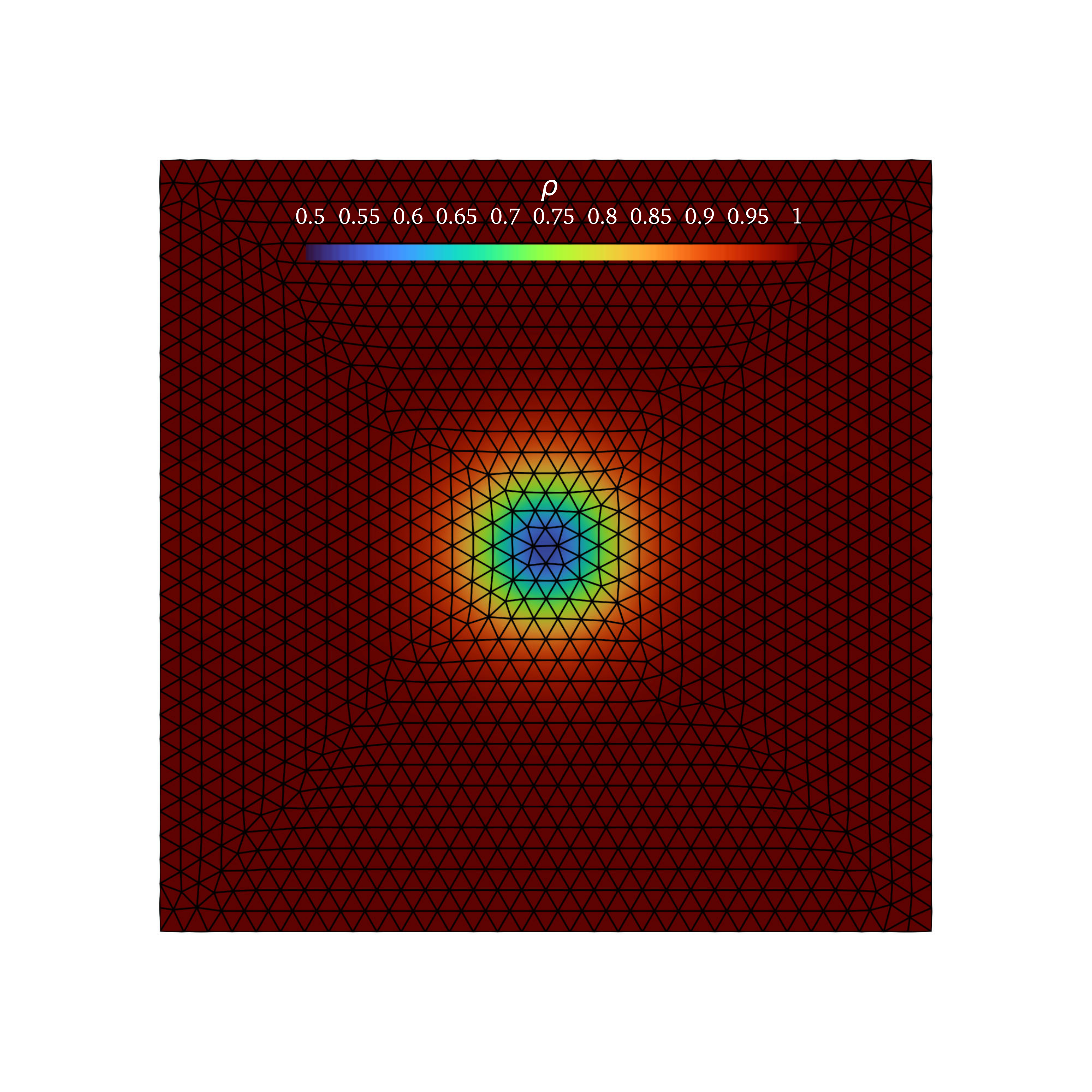}}
    \subfloat[Compasion of CPU Time for different methods]
    {\includegraphics[angle=0,width=0.47\textwidth,trim={0cm 0cm 0cm 0cm},clip]{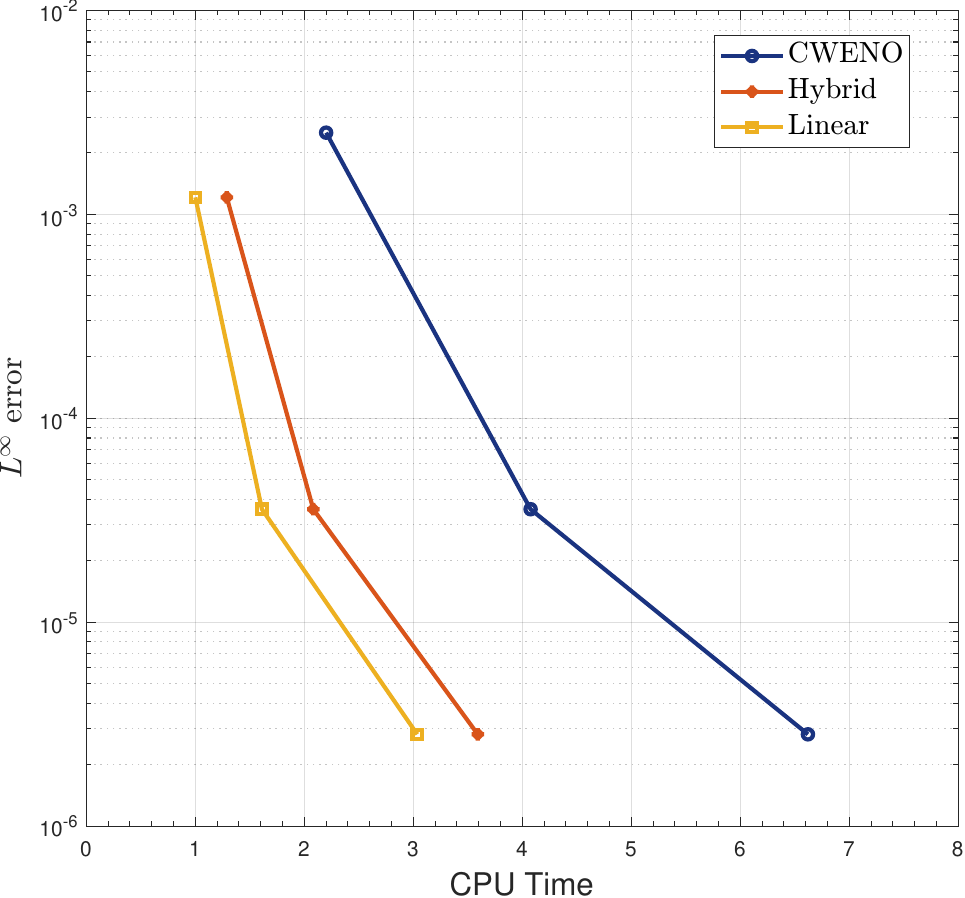}}

    \caption{(a) Density contour at $t=10$ using the Hybrid method on a $128 \times 128$ mesh with seventh-order accuracy.
    (b) Efficiency comparison of three methods (Linear, CWENOZ, and Hybrid) on the same $128 \times 128$ mesh. The data points indicate results at $\mathcal{P}3$, $\mathcal{P}5$, and $\mathcal{P}7$ (from top to bottom). The CPU time is normalised by the time required for the linear, third-order scheme on a coarse mesh which has 16 edges per side.}
    \label{fig:2D_Vortex_DensityCPUTime}
    \end{centering}
\end{figure}

The two-dimensional vortex evolution problem \citep{Balsara2000405} is used to assess both the accuracy and the convergence rate of our numerical method, as well as to evaluate the performance of the shock detector. The computational domain is a $10 \times 10$ square with periodic boundary conditions, and the total simulation time is set to 10. During this interval, a supersonic vortex propagates diagonally from the bottom-left to the top-right corner at a $45^\circ$ angle, ultimately returning to its original position by the end of the simulation due to the domain's periodicity. The flow is initialised uniformly with $(\rho, u, v, p) = (1, 1, 1, 1)$. The vortex is imposed as a fluctuation on this mean flow, characterised by:

\begin{equation}\label{eq:2DVortex}
\qquad \left(\delta u,\delta v\right)=\frac{\epsilon}{2\pi}e^{0.5\left(1-r^{2}\right)}\left(-\left(y-5\right),\left(x-5\right)\right), \quad \delta T=-\frac{\left(\gamma-1\right)\epsilon^{2}}{8\gamma\pi^{2}} e^{\left(1-r^{2}\right)}
\end{equation}

where $r^2 = (x-5)^2 + (y-5)^2$, the vortex strength is $\epsilon = 5$, and the specific heat ratio is $\gamma = 1.4$. The temperature and entropy are defined respectively by $T = p / \rho$ and $S = p / \rho^\gamma$. This problem is tested in four mesh resolutions, corresponding to the edges of $16$, $32$, $64$, and $128$ per side. All meshes are the triangular, which is presented in Fig. \ref{fig:2D_Vortex_DensityCPUTime}. For each of these resolutions, we perform spatial discretisations using $\mathcal{P}3$, $\mathcal{P}5$, and $\mathcal{P}7$. The error norms are then computed as follows:

\begin{equation}\label{eq:NormsofErrors}
e_{L^{\infty}}=\mathcal{M}ax \left|(\mathbf{U}_{e}\left(x,t_{f}\right)-\mathbf{U}_{c}\left(x,t_{f}\right)\right|
\quad\text{and}\quad
{e_{L^{2}}=\sqrt{\frac{\sum_i \int_{\Omega_{i}} \left({U}_{e}\left(x,t_{f}\right)-{U}_{c}\left(x,t_{f}\right)\right)^{2}dV}{\sum_i \left|\Omega_{i}\right|}}}
\end{equation}
where $\mathbf{U}_{e}$ denotes the exact solution derived from the initial condition, and $\mathbf{U}_{c}$ is the computed solution. The linear weight for the CWENOZ central stencil is set to $\lambda_1 = 10^{4}$, and a constant CFL number of 0.6 is used in all simulations. From Table \ref{tab:2D_Vortex}, it is observed that for lower-order polynomials, the CWENOZ scheme does not fully achieve its theoretical order of convergence. Moreover, on coarser meshes, its accuracy lags behind that of the \textcolor{black}{{high-order}} linear scheme. This disparity diminishes as the mesh is refined; at sufficiently high resolutions and polynomial orders, both methods yield nearly identical results. We did not observe the erroneous detection of troubled cells reported by Maltsev et al.\citep{MALTSEV2023111755}. The hybrid method effectively operates as a \textcolor{black}{{high-order}} linear scheme in smooth regions, producing results identical to those of the linear method. This confirms that in regions free of shocks or discontinuities, the hybrid scheme correctly defaults to a linear approach, as intended. Fig. \ref{fig:2D_Vortex_DensityCPUTime} further illustrates the computational efficiency of each method. For this smoothly varying problem, the hybrid scheme requires only about $18\%-29\%$ more CPU time than the fully \textcolor{black}{{high-order}} linear scheme. In contrast, the CWENOZ scheme requires approximately $71\%-96\%$ more computational resources than the hybrid scheme. Consequently, in the smooth regions of 2D flows, the hybrid approach can reduce computational costs by nearly half compared to CWENOZ, while maintaining comparable accuracy, and even surpassing it at lower polynomial orders.

\begin{table*}[h!]
\caption{Numerical errors and convergence rates for the multi-resolution 2D vortex evolution problem, comparing the Linear, CWENOZ, and Hybrid schemes. In this particular test, the flow is sufficiently smooth that no troubled cells are detected. As a result, the Hybrid scheme effectively operates as a 100\% \textcolor{black}{{high-order}} Linear scheme.}
\begin{center}
\resizebox{\textwidth}{!}{
\begin{tabular}{c@{\qquad}|cccc|cccc|cccc}
\toprule
\multicolumn{1}{c}{Order/Number of Edges} & \multicolumn{4}{c}{Linear} & \multicolumn{4}{c}{CWENOZ} & \multicolumn{4}{c}{Hybrid} \\
\cmidrule{1-13}
     & $e_{L^{\infty}}$ & $\mathcal{O}_{L^{\infty}}$ & $e_{L^{2}}$ & $\mathcal{O}_{L^{2}}$ & $e_{L^{\infty}}$ & $\mathcal{O}_{L^{\infty}}$ & $e_{L^{2}}$ & $\mathcal{O}_{L^{2}}$ & $e_{L^{\infty}}$ & $\mathcal{O}_{L^{\infty}}$ & $e_{L^{2}}$ & $\mathcal{O}_{L^{2}}$ \\
\midrule
$ \mathcal{P}3 / 16 $ & 2.417E-01 & - & 2.936E-02 & - & 2.599E-01 & - & 3.330E-02 & - & 2.417E-01 & - & 2.936E-02 & - \\
$ \mathcal{P}3 / 32 $ & 4.923E-02 & 2.30 & 7.183E-03 & 2.03 & 5.895E-02 & 2.14 & 7.855E-03 & 2.08 & 4.923E-02 & 2.30 & 7.183E-03 & 2.03 \\
$ \mathcal{P}3 / 64 $ & 9.216E-03 & 2.42 & 1.390E-03 & 2.37 & 1.332E-02 & 2.15 & 1.681E-03 & 2.22 & 9.216E-03 & 2.42 & 1.390E-03 & 2.37 \\
$ \mathcal{P}3 / 128 $ & 1.211E-03 & 2.93 & 1.945E-04 & 2.84 & 2.515E-03 & 2.41 & 2.345E-04 & 2.84 & 1.211E-03 & 2.93 & 1.945E-04 & 2.84 \\
\midrule
$ \mathcal{P}5 / 16 $ & 9.349E-02 & - & 1.185E-02 & - & 1.015E-01 & - & 1.313E-02 & - & 9.349E-02 & - & 1.185E-02 & - \\
$ \mathcal{P}5 / 32 $ & 2.894E-02 & 1.69 & 2.357E-03 & 2.33 & 2.889E-02 & 1.81 & 2.313E-03 & 2.51 & 2.894E-02 & 1.69 & 2.357E-03 & 2.33 \\
$ \mathcal{P}5 / 64 $ & 8.260E-04 & 5.13 & 1.426E-04 & 4.05 & 8.260E-04 & 5.13 & 1.426E-04 & 4.02 & 8.260E-04 & 5.13 & 1.426E-04 & 4.05 \\
$ \mathcal{P}5 / 128 $ & 3.579E-05 & 4.53 & 5.642E-06 & 4.66 & 3.579E-05 & 4.53 & 5.642E-06 & 4.66 & 3.579E-05 & 4.53 & 5.642E-06 & 4.66 \\
\midrule
$ \mathcal{P}7 / 16 $ & 4.584E-02 & - & 5.404E-03 & - & 6.119E-02 & - & 1.054E-02 & - & 4.584E-02 & - & 5.404E-03 & - \\
$ \mathcal{P}7 / 32 $ & 1.583E-02 & 1.53 & 1.149E-03 & 2.23 & 1.578E-02 & 1.95 & 1.113E-03 & 3.24 & 1.583E-02 & 1.53 & 1.149E-03 & 2.23 \\
$ \mathcal{P}7 / 64 $ & 3.254E-04 & 5.60 & 3.443E-05 & 5.06 & 3.254E-04 & 5.60 & 3.443E-05 & 5.01 & 3.254E-04 & 5.60 & 3.443E-05 & 5.06 \\
$ \mathcal{P}7 / 128 $ & 2.818E-06 & 6.85 & 4.045E-07 & 6.41 & 2.817E-06 & 6.85 & 4.046E-07 & 6.41 & 2.818E-06 & 6.85 & 4.045E-07 & 6.41 \\
\cmidrule{1-13}
\end{tabular}}
\end{center}
\label{tab:2D_Vortex}
\end{table*}

\pagebreak
\subsection{Shu-Osher Problem}
The Shu-Osher problem \cite{SHU198932}, commonly referred to as the shock-entropy wave interaction case, is a well-known benchmark for evaluating the performance of numerical schemes in capturing shock waves. In this test case, the initial condition features a shock wave traveling from left to right, which interferes with an entropy wave. The initial conditions can be described as:
\begin{equation}\label{eq:2DShuOsherProblem}
(\rho, u, v, p) =
\begin{cases}
(3.857143, 2.629369, 0, 10.3333) & \text{for } x < -4, \\
(1 + 0.2\sin(5x), 0, 0, 1) & \text{for } x \geq -4.
\end{cases}
\end{equation}

The computational domain is defined as \(x \in [-4.5, 4.5]\) and \(y \in [0, 1]\). The upper and lower limits are defined with periodic boundary conditions, the left boundary is set as a supersonic inlet, and the right boundary is specified as an outflow. The total simulation time for the test is set to \(t = 1.8\). The reference solution is obtained by solving the Euler equations using a 5th-order WENO scheme on a mesh with 10,000 grid points. Various deactivation criteria, including \(\Delta x^{\frac{2}{3}}\), \(\Delta x\), \(\Delta x^2\), and \((5\Delta x)^2\), are tested on three levels of triangular meshes, with 225, 450, and 900 grid points in the \(x\)-direction. All simulations are using 5th-order spatial reconstruction. 

\begin{figure}[h!]
\begin{centering}
\captionsetup[subfigure]{width=0.3\textwidth}
  \subfloat[Coarse Mesh]
 {\includegraphics[angle=0,width=0.33\textwidth,trim={0cm 0cm 0cm 0cm},clip]{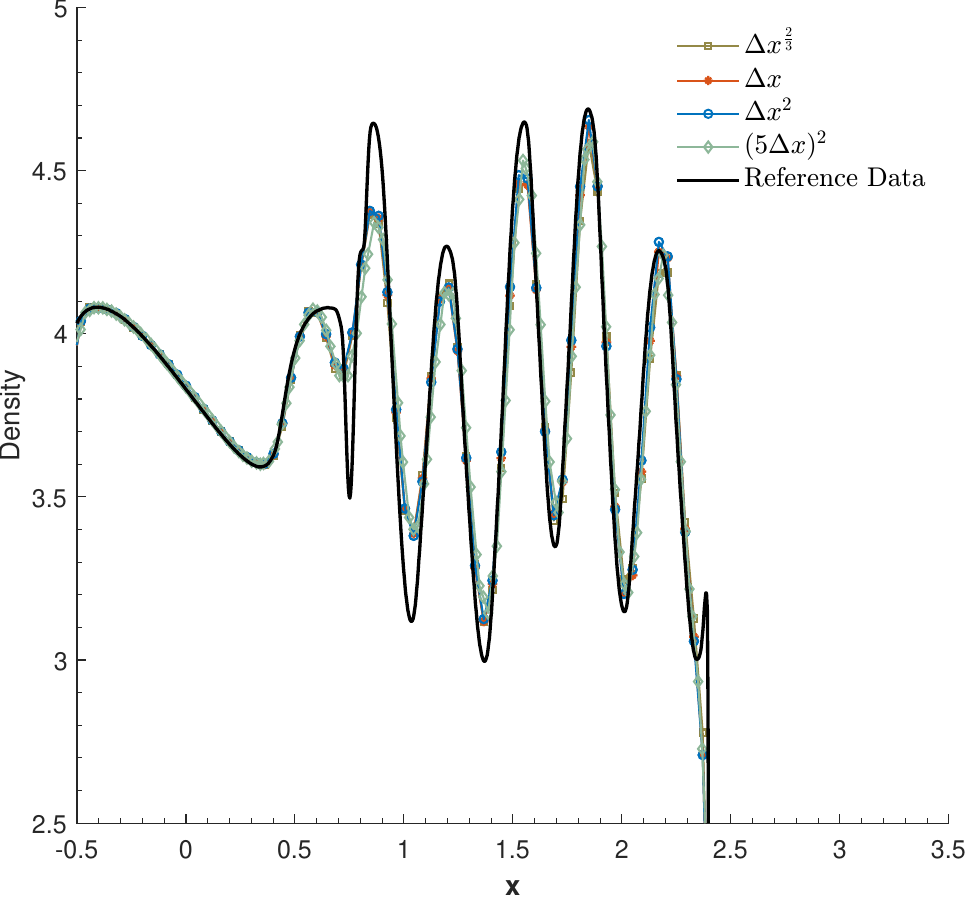}}
 \subfloat[Medium Mesh]
 {\includegraphics[angle=0,width=0.33\textwidth,trim={0cm 0cm 0cm 0cm},clip]{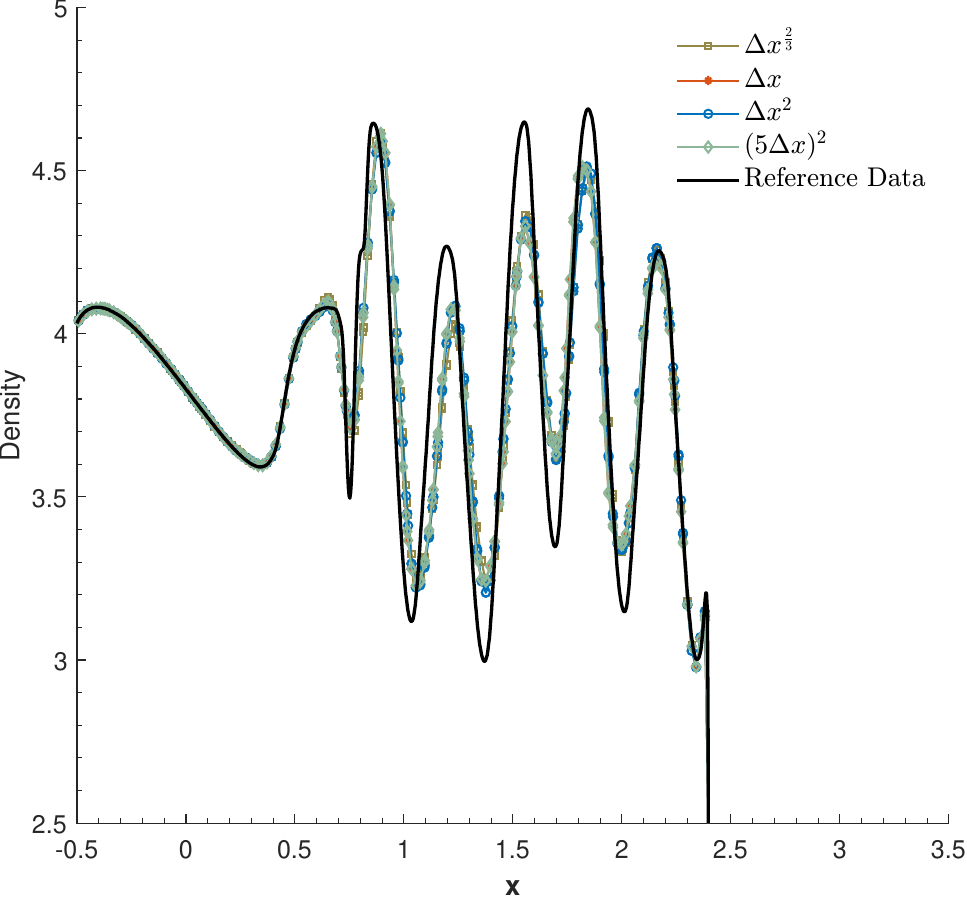}}
 \subfloat[Fine Mesh]
 {\includegraphics[angle=0,width=0.33\textwidth,trim={0cm 0cm 0cm 0cm},clip]{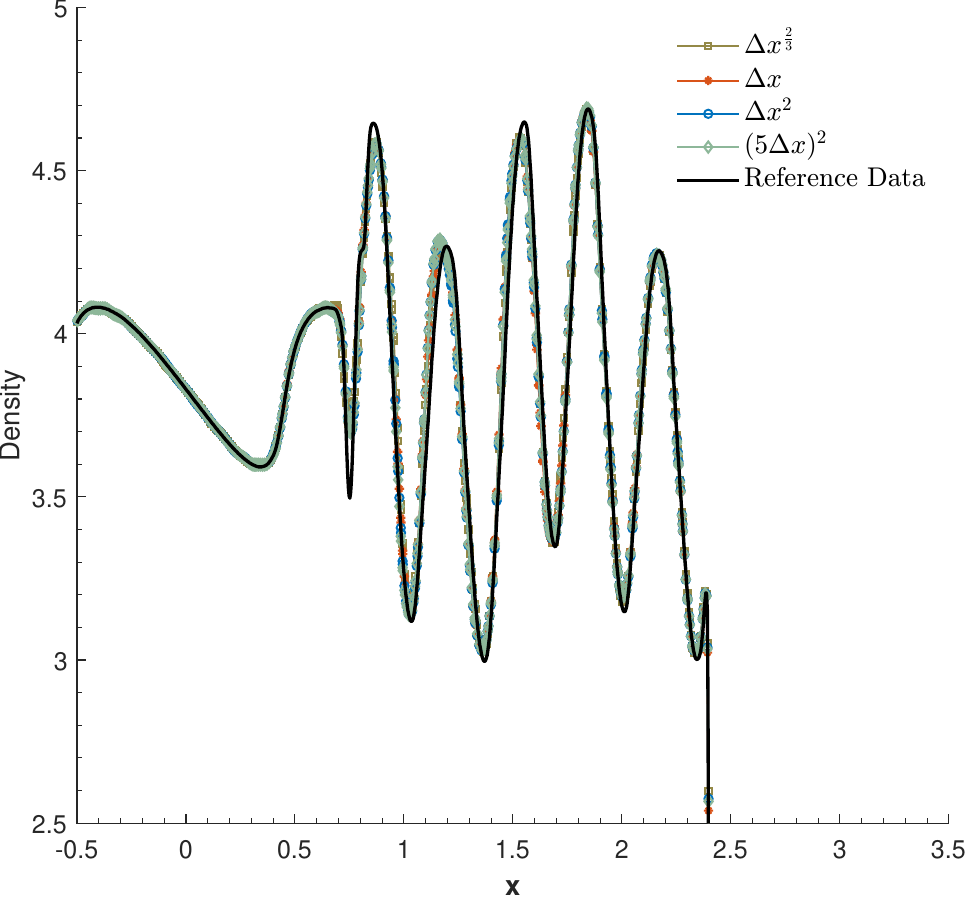}}
 
\par\end{centering}\caption{Density distribution at $t=1.8$ for the Shu–Osher problem, comparing different deactivation criterias across three mesh levels with 5th-order spatial reconstruction}\label{fig:ShuOsher1}\end{figure}

In general, the computed results are consistent with the reference data. As the mesh density increases, the accuracy of the solution improves. Fig. \ref{fig:ShuOsher1} illustrates that the switch value is only used to deactivate the NAD detector in smooth regions, preventing the smooth region from being calculated using the CWENOZ scheme. The results are not sensitive to different deactivation criteria across varying mesh resolutions. On the medium mesh, the hybrid scheme demonstrates performance comparable to the CWENOZ scheme, even exhibiting superior accuracy at certain locations.

\begin{figure}[h!]
\begin{centering}
\captionsetup[subfigure]{width=0.4\textwidth}
  \subfloat[Comparison between different Schemes]
 {\includegraphics[angle=0,width=0.45\textwidth,trim={0cm 0cm 0cm 0cm},clip]{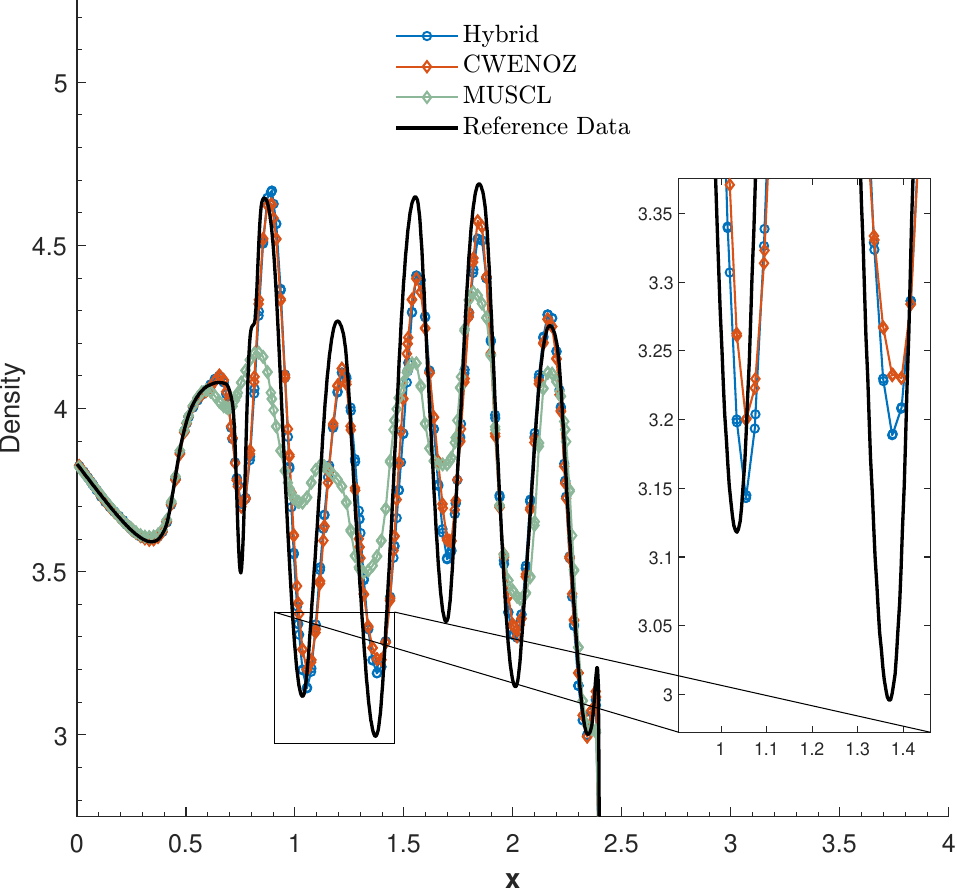}}
 \subfloat[Hybrid scheme with different limiters]
 {\includegraphics[angle=0,width=0.45\textwidth,trim={0cm 0cm 0cm 0cm},clip]{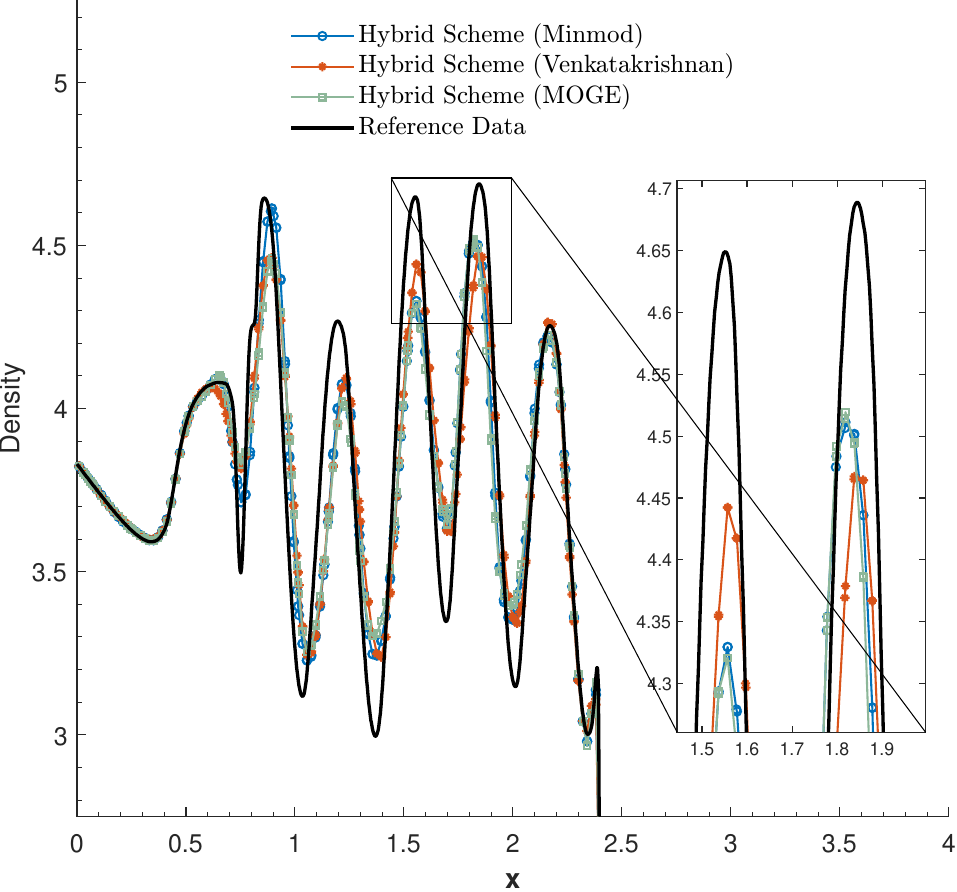}}

\par\end{centering}\caption{(a) Density distribution at $t = 1.8$, obtained using the Minmod limiter within both the Hybrid and MUSCL schemes. (b) Density distribution at $t = 1.8$, obtained using Minmod, MOGE and the Venkatakrishnan limiter with a tunable parameter of $\kappa=5$. All simulations are conducted on the medium mesh.
}\label{fig:ShuOsher2}\end{figure}

A comparison of the classical Barth \& Jespersen (Minmod) limiter \cite{Barth1989366}, the high-order MOGE limiter \cite{TSOUTSANIS201869}, and the Venkatakrishnan limiter \cite{VENKATAKRISHNAN1995120} clarifies how the choice of slope limiter interacts with the NAD shock detector in our hybrid scheme. In principle, the MOGE limiter should deliver higher accuracy than the second-order Barth \& Jespersen and Venkatakrishnan limiters. However, the results in Fig.~\ref{fig:ShuOsher2} contradict this expectation. Barth \& Jespersen enforces strict local monotonicity by constraining each reconstructed state to the range defined by its immediate neighbours, thereby guaranteeing total-variation-diminishing (TVD) behaviour. Although this aggressive clipping introduces additional numerical dissipation, it keeps the shock detector perfectly aligned with the limiter so that only genuinely discontinuous cells are flagged for MUSCL reconstruction. By contrast, MOGE relaxes the TVD constraint by deriving its admissible bounds from a larger stencil. A reconstructed state may therefore overshoot its nearest neighbour, violating the detector's monotonicity criterion; the detector then misclassifies extra cells as discontinuous, unnecessarily expanding the MUSCL and first-order upwind region. Despite MOGE's low intrinsic dissipation, this over-activation of MUSCL degrades its overall accuracy compared with Barth \& Jespersen or Venkatakrishnan limiter, and it still underperforms relative to CWENOZ or \textcolor{black}{{high-order}} linear scheme. The Venkatakrishnan limiter employs a similar stencil that matches the NAD detector's, its behaviour resembles that of Barth \& Jespersen. Nevertheless, Barth \& Jespersen remains the most attractive option due to its lower computational cost.

\newpage
\subsection{Double Mach Reflection}
The double Mach reflection test, introduced by Woodward and Colella \citep{Woodward1984115}, serves as a benchmark for evaluating the accuracy and robustness of numerical algorithms. This test generates strong shock waves and discontinuities, leading to the formation of an intricate flow structure characterised by multiple shock interactions, despite originating from a simple geometric configuration. In this test case, the computational domain is defined as $x \in [0,4]$, $y \in [0,1]$. The x-axis of the computational domain is aligned with the upper edge of a wedge, which forms a $30^\circ$ angle with the original x-axis. The left boundary is assigned supersonic inlet boundary condition, while the right boundary is set as an outflow boundary. The upper boundary is prescribed with the exact solution for the propagating oblique shock wave. The initial position of the shock wave is set at $x = 1/6$, where the reflecting boundary condition begins. The shock propagates at a Mach number of 10, with the pre-shock flow state given by $[\rho, p, \gamma] = [1.4, 0, 1.4]$.

\begin{figure}[h!]
 \begin{centering}
 
 \captionsetup[subfigure]{width=0.30\textwidth}
  \subfloat[$\mathcal{P}3$/ CWENOZ]
 {\includegraphics[angle=0,width=0.33\textwidth,trim={3cm 25cm 0cm 25cm},clip]{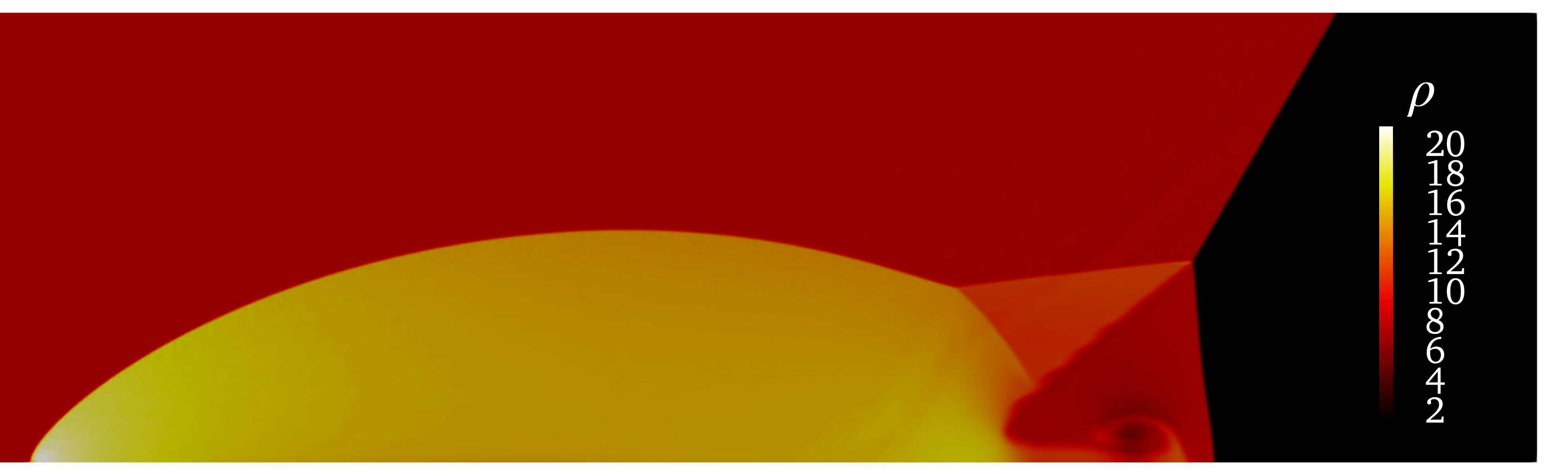}}
 \subfloat[$\mathcal{P}3$/ Hybrid]
 {\includegraphics[angle=0,width=0.33\textwidth,trim={3cm 25cm 0cm 25cm},clip]{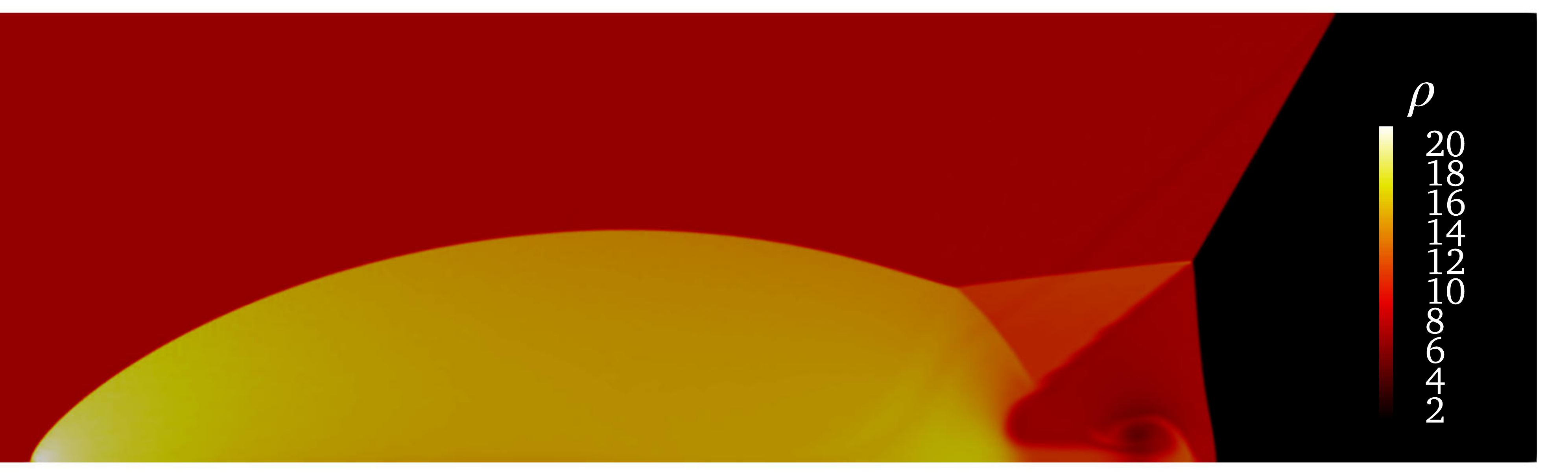}}
 \subfloat[$\mathcal{P}3$/ Hybrid]
 {\includegraphics[angle=0,width=0.33\textwidth,trim={3cm 25cm 0cm 25cm},clip]{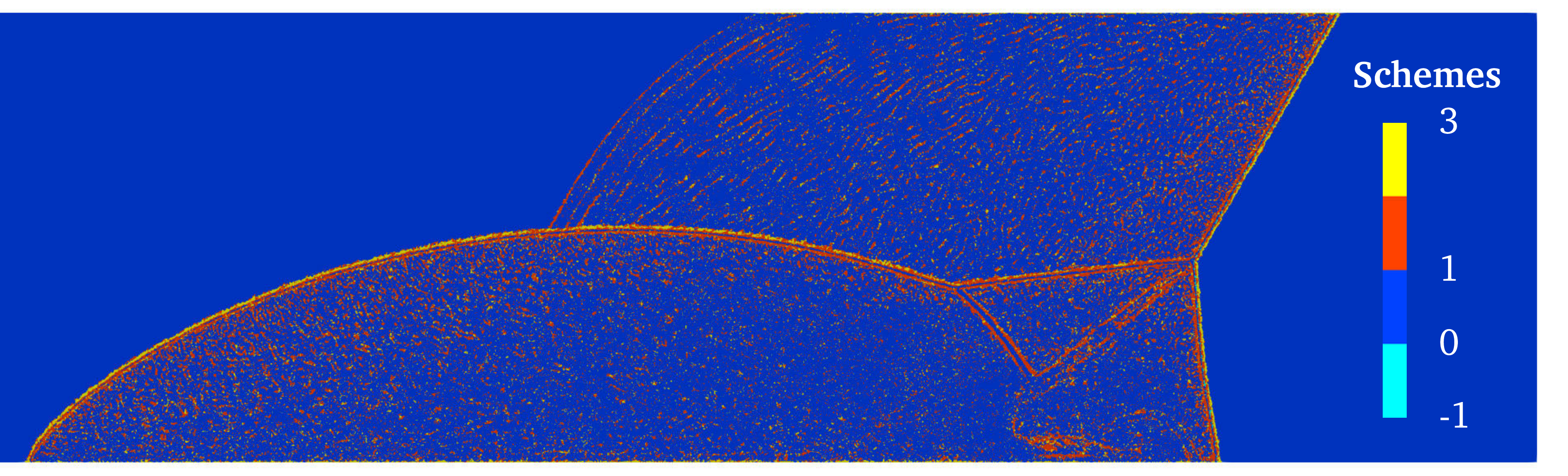}}

  \subfloat[$\mathcal{P}5$/ CWENOZ]
 {\includegraphics[angle=0,width=0.33\textwidth,trim={3cm 25cm 0cm 25cm},clip]{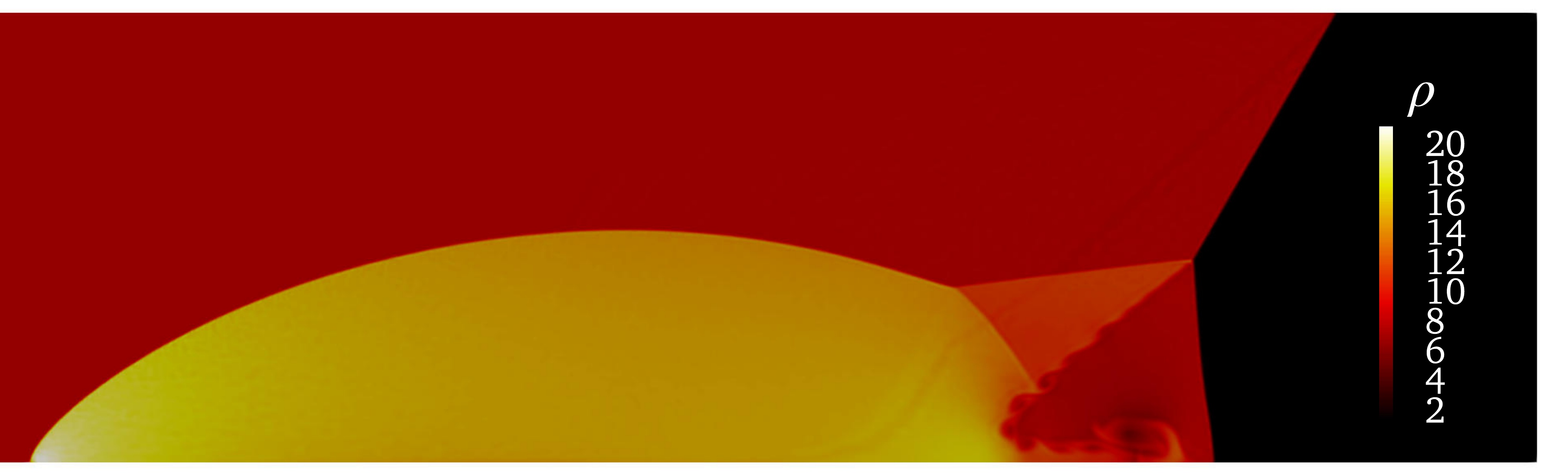}}
 \subfloat[$\mathcal{P}5$/ Hybrid]
 {\includegraphics[angle=0,width=0.33\textwidth,trim={3cm 25cm 0cm 25cm},clip]{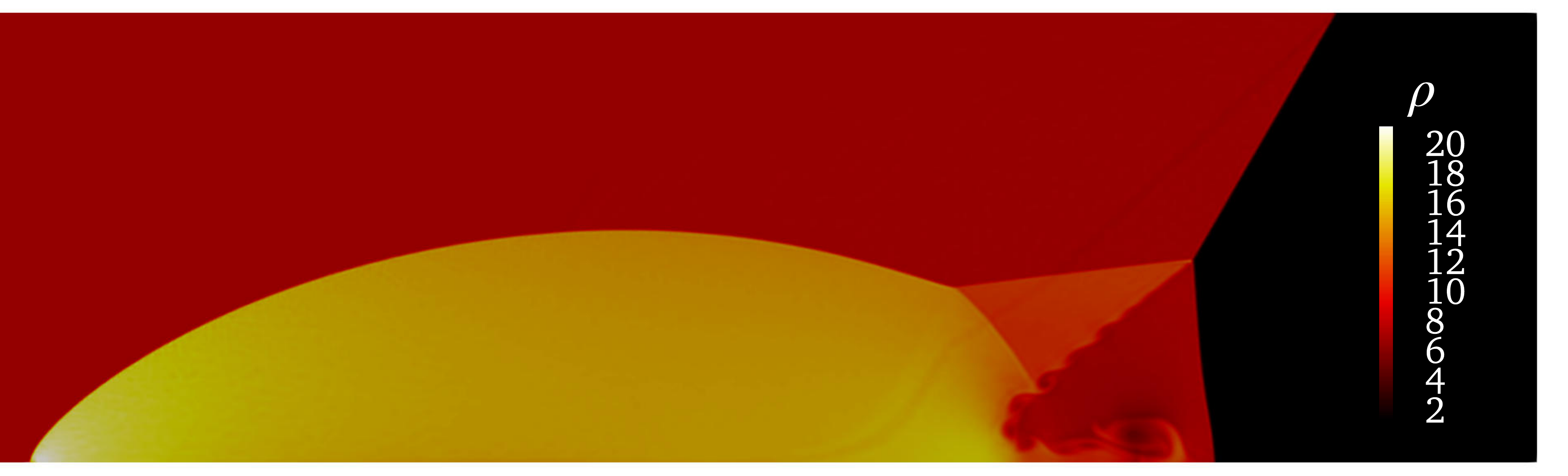}}
 \subfloat[$\mathcal{P}5$/ Hybrid]
 {\includegraphics[angle=0,width=0.33\textwidth,trim={3cm 25cm 0cm 25cm},clip]{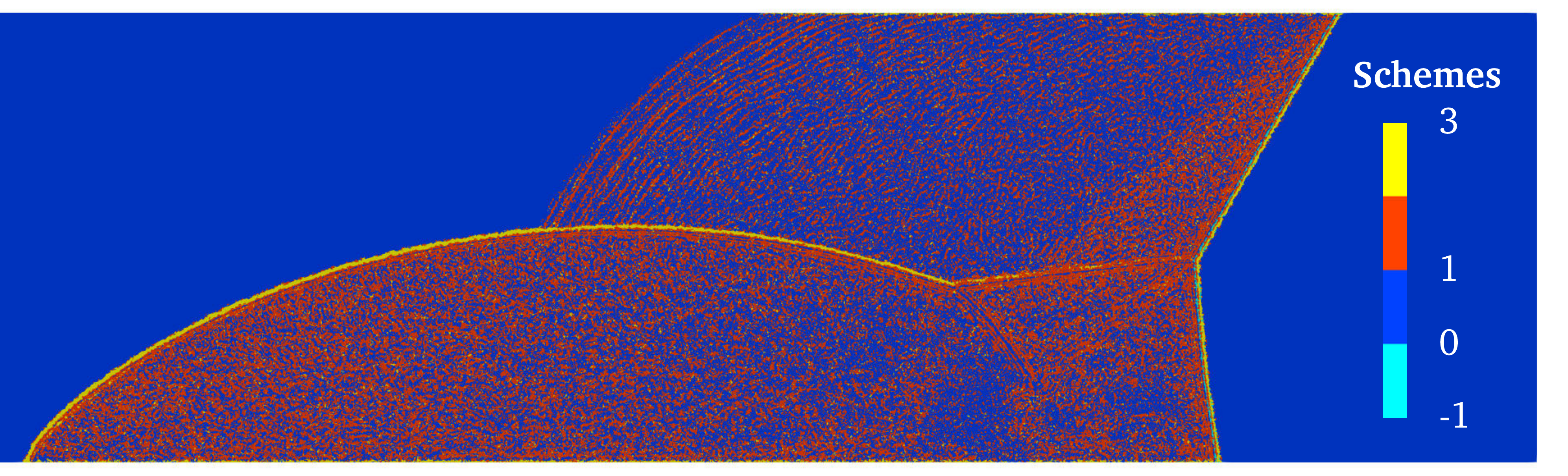}}

  \subfloat[$\mathcal{P}7$/ CWENOZ]
 {\includegraphics[angle=0,width=0.33\textwidth,trim={3cm 25cm 0cm 25cm},clip]{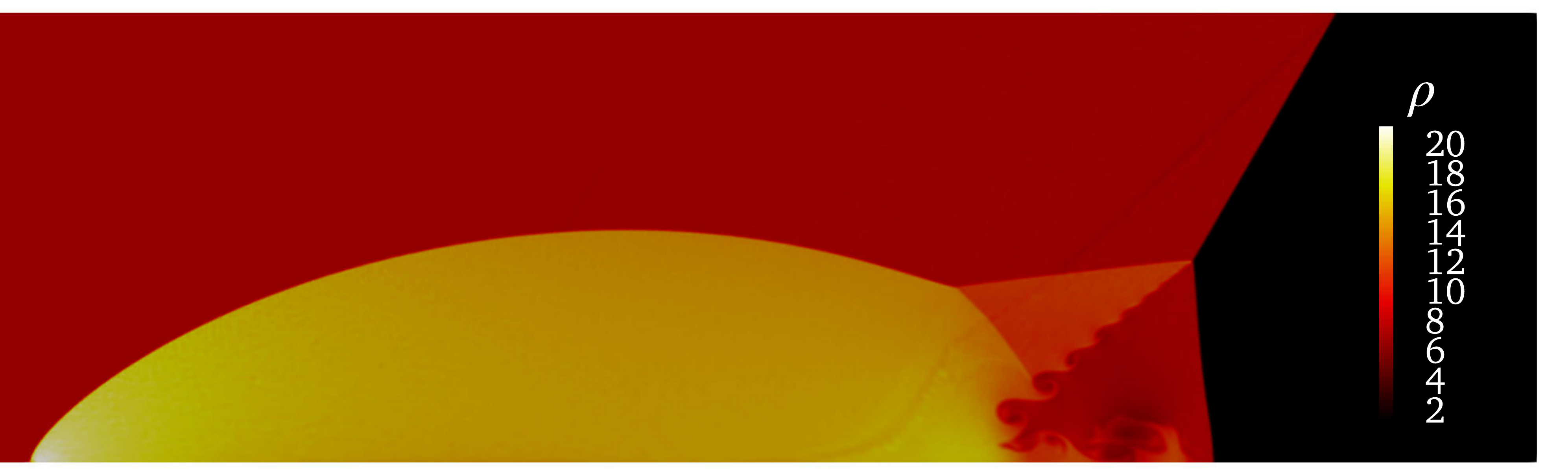}}
 \subfloat[$\mathcal{P}7$/ Hybrid]
 {\includegraphics[angle=0,width=0.33\textwidth,trim={3cm 25cm 0cm 25cm},clip]{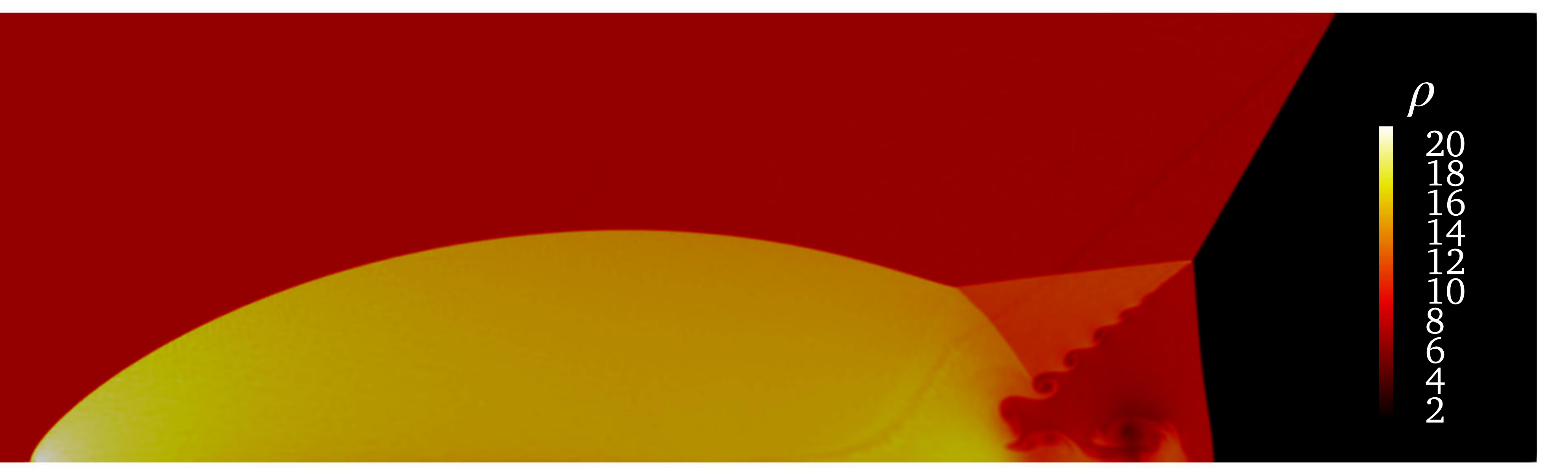}}
 \subfloat[$\mathcal{P}7$/ Hybrid]
 {\includegraphics[angle=0,width=0.33\textwidth,trim={3cm 25cm 0cm 25cm},clip]{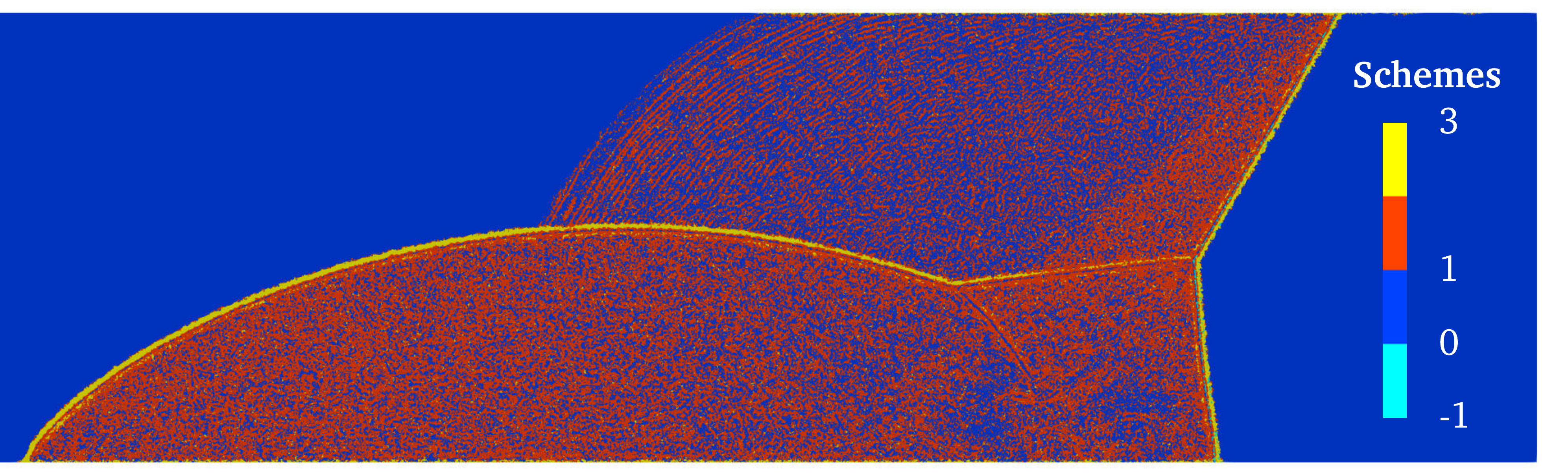}}

 \par\end{centering}\caption{Density contours computed using the CWENOZ scheme (\textit{left}); Density contours obtained via the Hybrid scheme (\textit{middle}); numerical schemes contours in the Hybrid method (\textit{right}). The scheme legend illustrates: yellow regions (Scheme 3) denote cells employing the $2nd-$ order MUSCL scheme, red zones (Scheme 1) indicate high-order CWENOZ scheme, deep blue areas (Scheme 0) represent \textcolor{black}{{high-order}} linear scheme, and light blue sections (Scheme -1) correspond to $1st-order$ upwind treatment. All result were conducted until final time $t=2$.}.
 \label{fig:DMR_DifferentOrders}\end{figure}

\begin{figure}[h!]
    \begin{centering}
    \captionsetup[subfigure]{width=0.30\textwidth}

    \subfloat[$\mathcal{P}3$/ CWENOZ)]
    {\includegraphics[angle=0,width=0.32\textwidth,trim={9cm 9cm 10cm 10cm},clip]{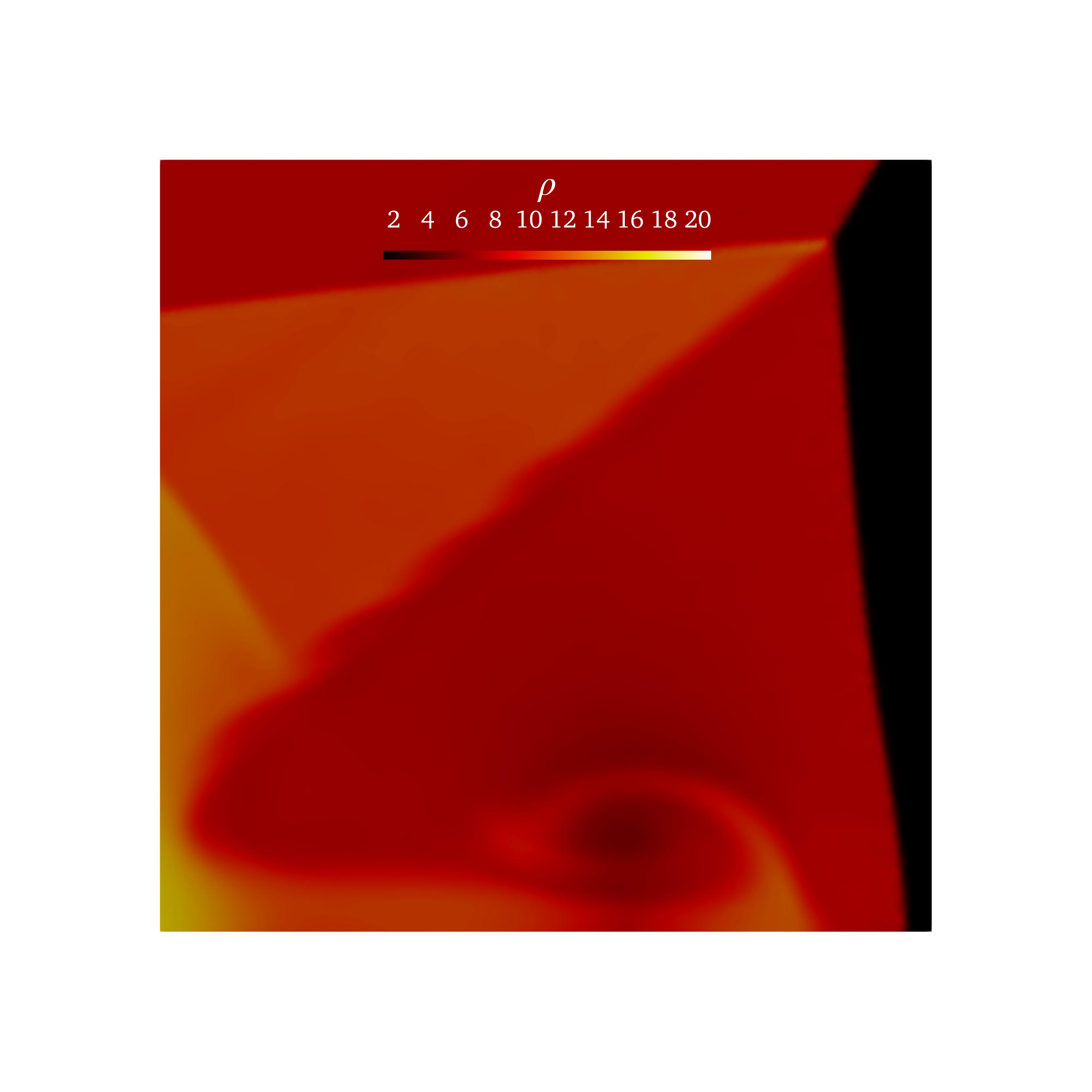}}
    \subfloat[$\mathcal{P}5$/ CWENOZ]
    {\includegraphics[angle=0,width=0.32\textwidth,trim={9cm 9cm 10cm 10cm},clip]{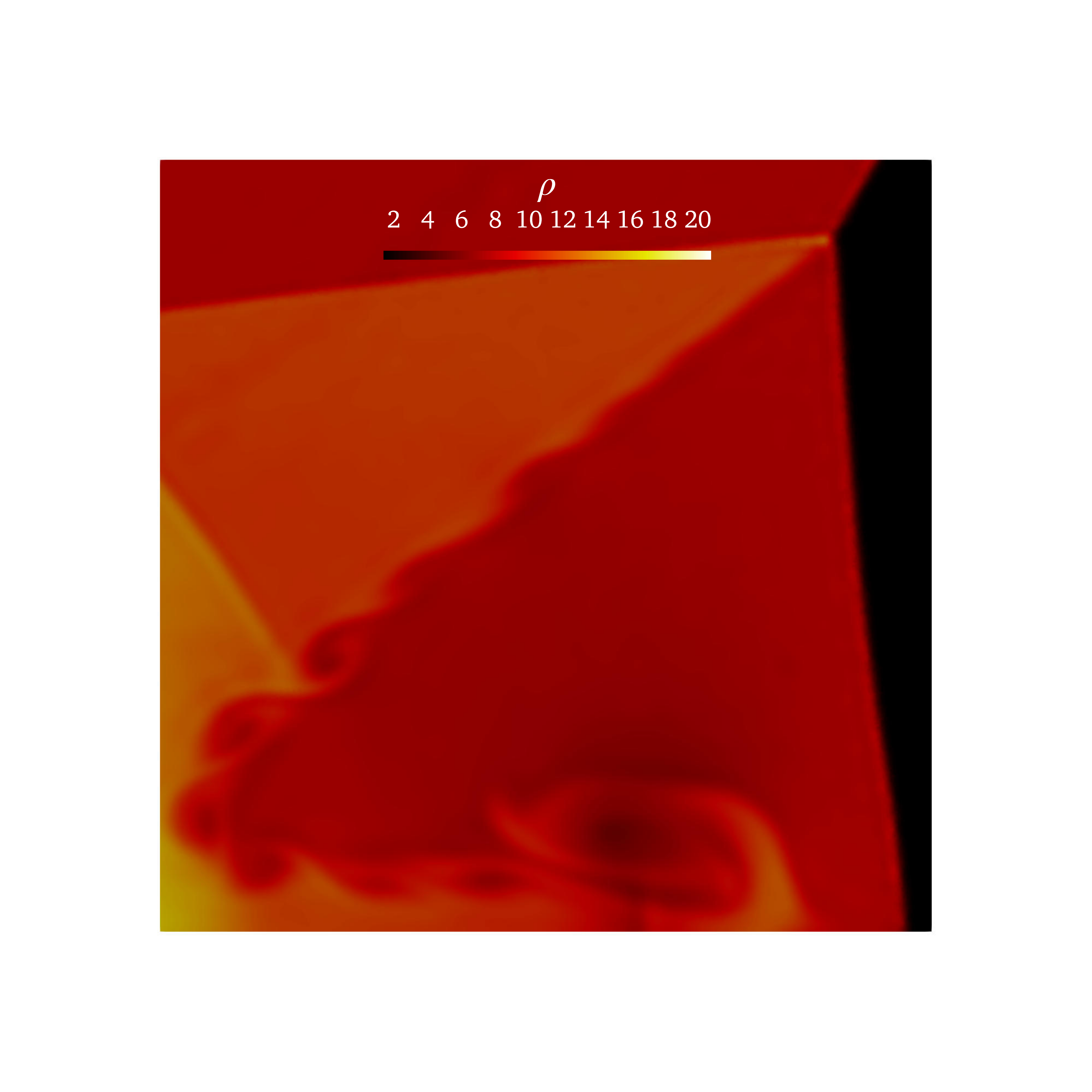}}
    \subfloat[$\mathcal{P}7$/ CWENOZ]
    {\includegraphics[angle=0,width=0.32\textwidth,trim={9cm 9cm 10cm 10cm},clip]{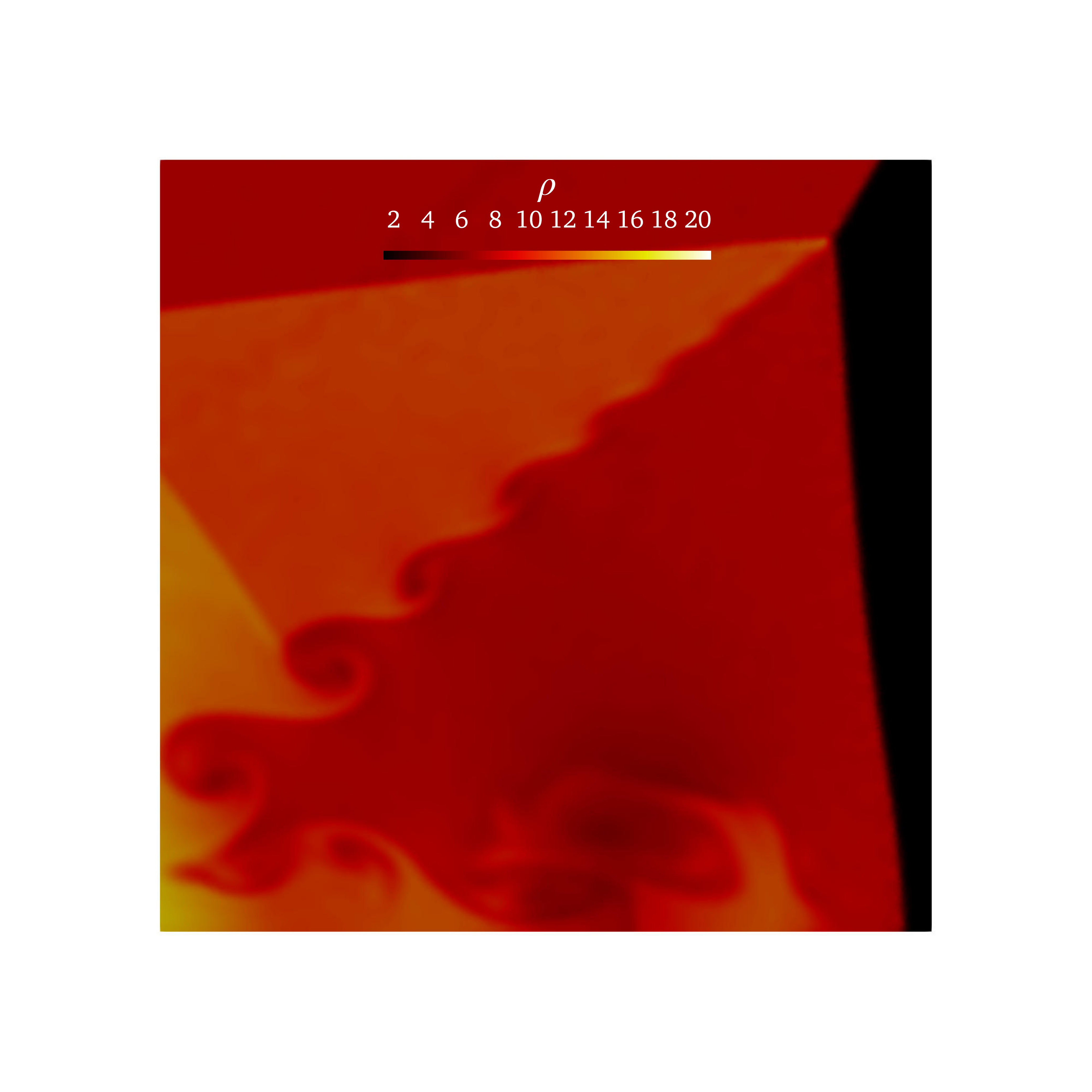}}\\
   
    \subfloat[$\mathcal{P}3$/ Hybrid]
    {\includegraphics[angle=0,width=0.32\textwidth,trim={9cm 9cm 10cm 10cm},clip]{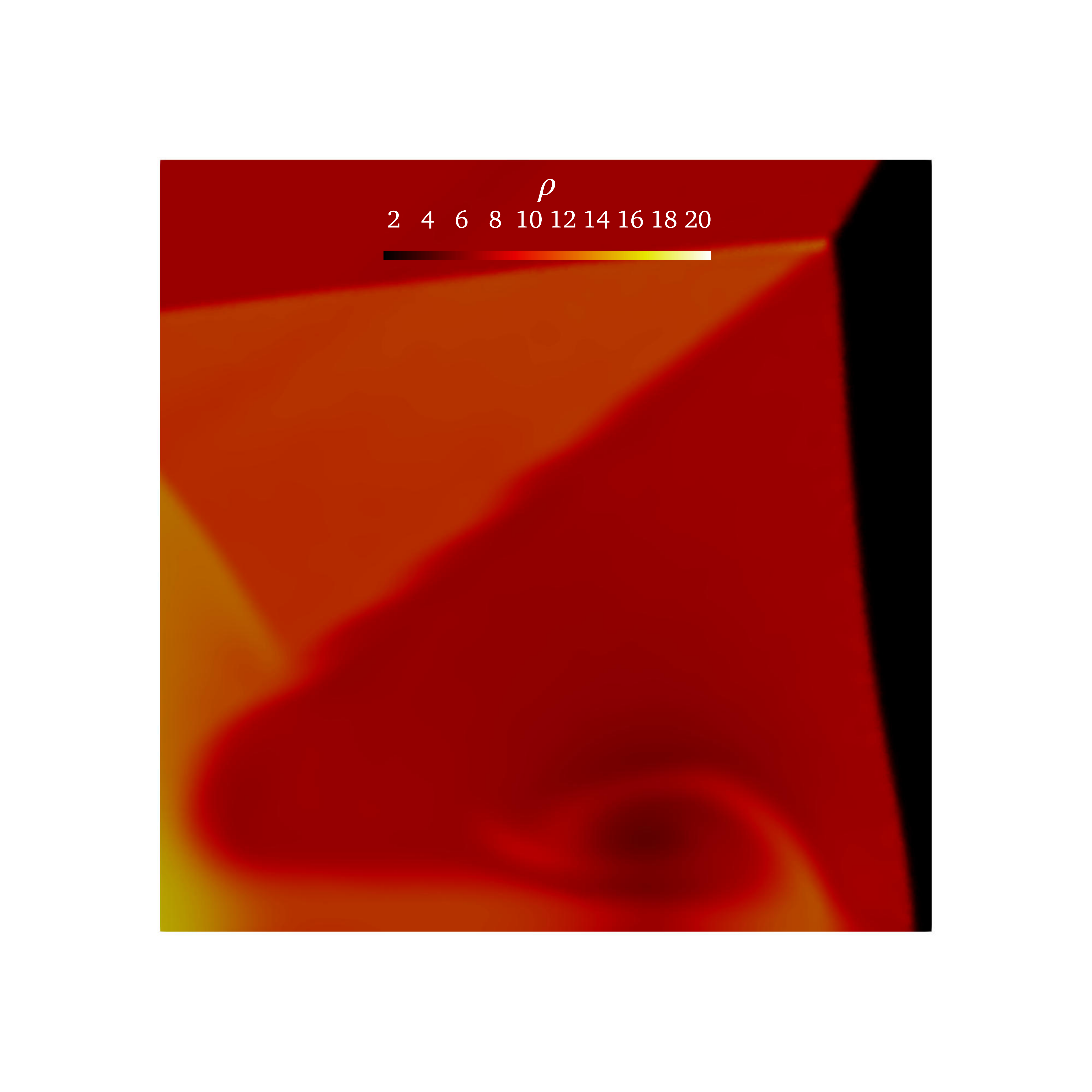}}
    \subfloat[$\mathcal{P}5$/ Hybrid]
    {\includegraphics[angle=0,width=0.32\textwidth,trim={9cm 9cm 10cm 10cm},clip]{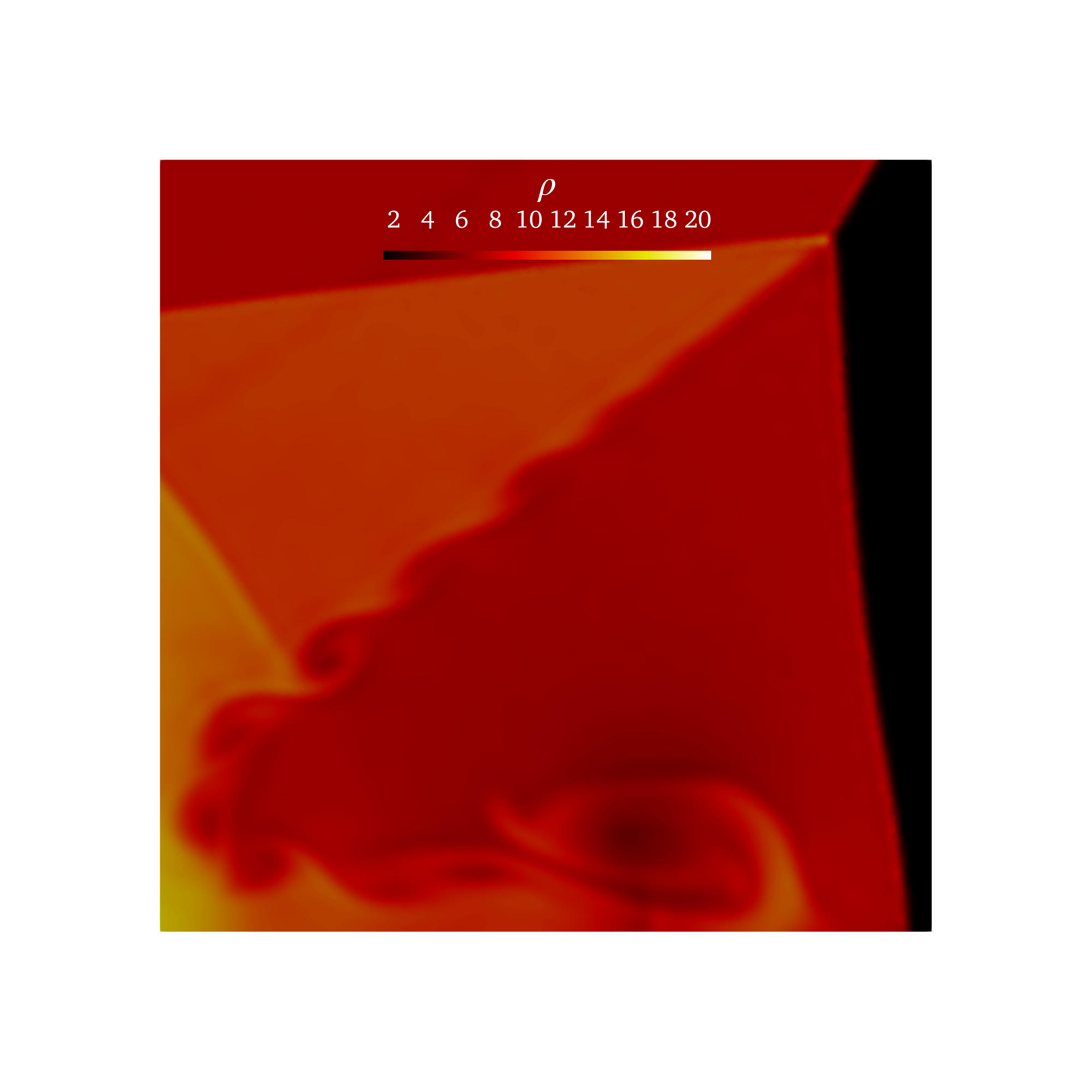}}
    \subfloat[$\mathcal{P}7$/ Hybrid]
    {\includegraphics[angle=0,width=0.32\textwidth,trim={9cm 9cm 10cm 10cm},clip]{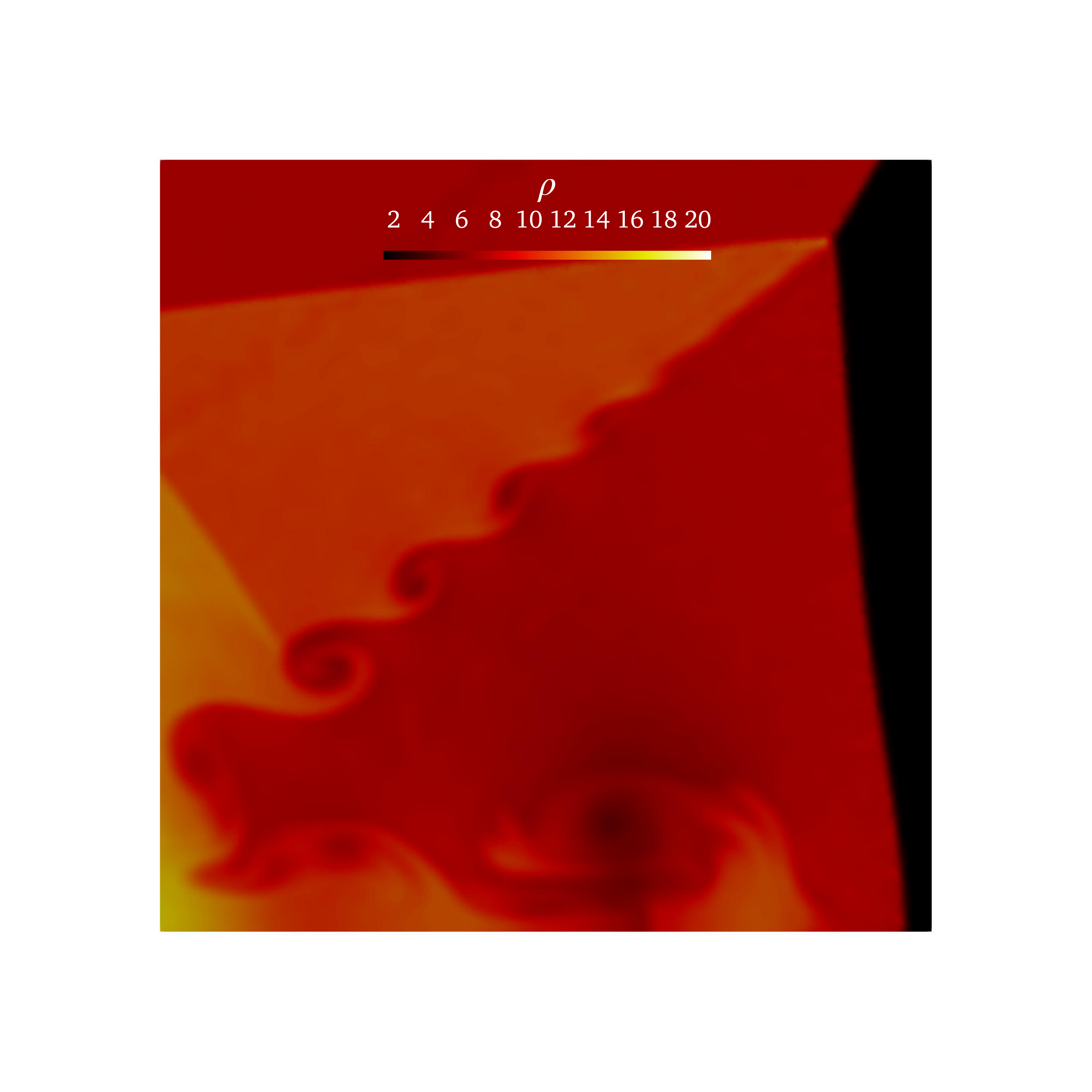}}

    \caption{The zoomed-in density contours reveal the details of the wave interaction zone. It is observed that as the order of accuracy increases, the results obtained from the Hybrid schemes lose some small-scale features. In particular, when comparing the $\mathcal{P}7$ CWENOZ and $\mathcal{P}7$ Hybrid schemes, the vortices along the slipstream are noticeably less pronounced in the Hybrid scheme. This issue can be alleviated by relaxing the Hybrid criteria.}
    \label{fig:DMR_WaveInteractionArea}
    \end{centering}
\end{figure}

The numerical simulations were conducted until time $t = 2$, employing a $CFL=0.6$ and the $\lambda_1 = 10^{3}$ for all simulations. The test cases were performed on a hybrid unstructured mesh (combination of triangles and quadrilaterals) with spatial resolution $h = 1/320$, and with reconstruction polynomial orders $\mathcal{P}3$, $\mathcal{P}5$, and $\mathcal{P}7$. Fig.  \ref{fig:DMR_DifferentOrders} presents the resulting density contours and corresponding scheme selection patterns, while Fig.  \ref{fig:DMR_WaveInteractionArea} provides detailed visualisation of the wave interaction zone in this test case.

\begin{figure}[h!]
    \begin{centering}
    \captionsetup[subfigure]{width=0.4\textwidth}

    \subfloat[Hybrid Setting: $\beta_m=0.5$]
    {\includegraphics[angle=0,width=0.4\textwidth,trim={9cm 9cm 10cm 10cm},clip]{DMR_7th_Detail.pdf}}
    \subfloat[Hybrid setting: $\beta_m=2$]
    {\includegraphics[angle=0,width=0.4\textwidth,trim={9cm 9cm 10cm 10cm},clip]{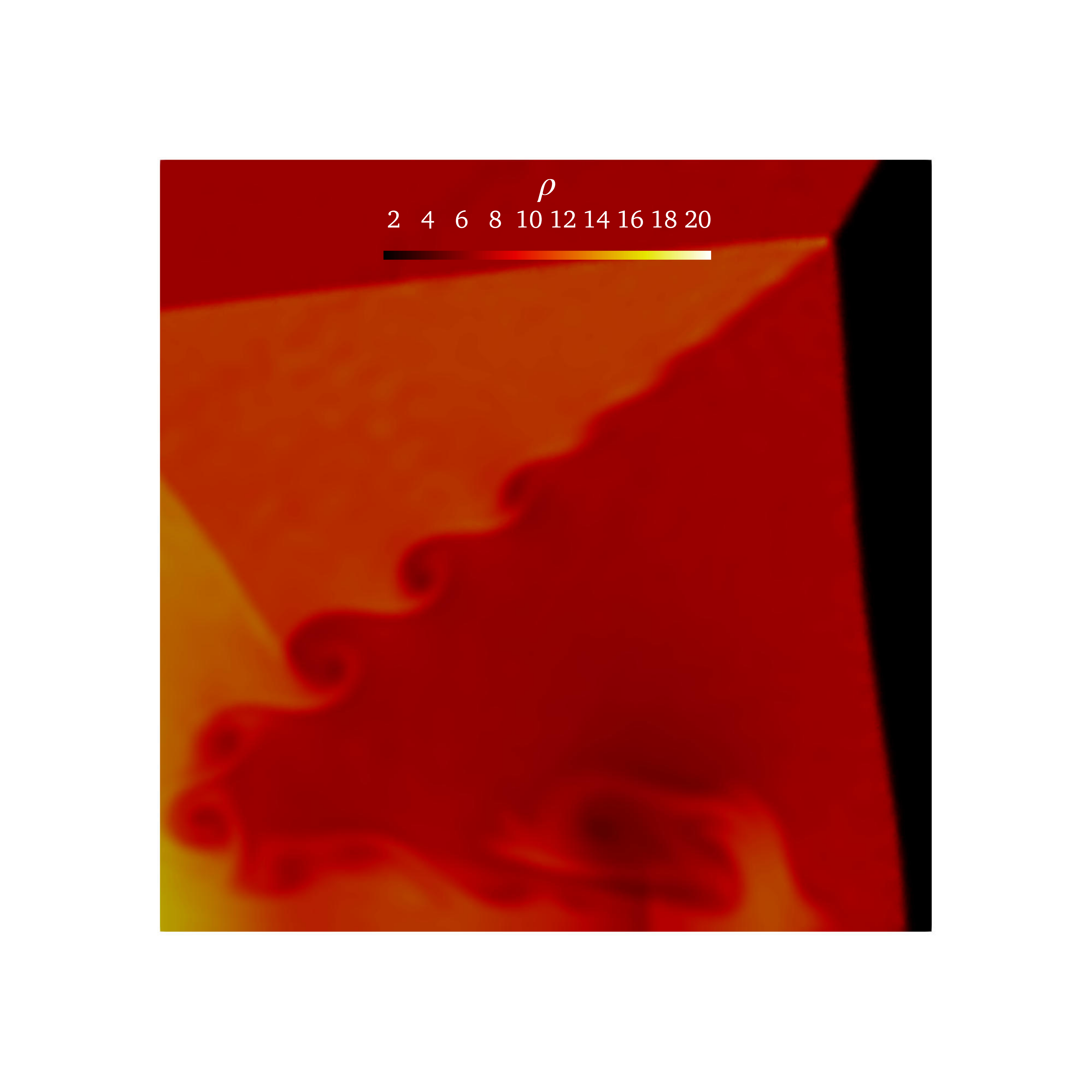}}

    \caption{The density distributions for different Hybrid settings with $\mathcal{P}7$ show that relaxing the criteria reduces the numerical diffusion introduced by the MUSCL scheme. This improvement is achieved by incorporating more contributions from the CWENOZ scheme into the Hybrid method.}
    \label{fig:DMR_DifferentHybridSetting}
    \end{centering}
\end{figure}

The numerical results demonstrate that all flow features, including Mach stems, shock reflections, and contact discontinuities, are consistently resolved across different accuracy orders. Although the hybrid scheme achieves improved accuracy with increasing polynomial order, subtle discrepancies emerge at higher orders ($\mathcal{P}7$) compared to the pure CWENOZ implementation. As shown in Fig. \ref{fig:DMR_DifferentOrders}, higher-order simulations exhibit increased activation of MUSCL and CWENOZ schemes. This occurs because higher-order simulations will be able to resolve finer flow structures, the detailed flow structure like the vortices, pressure waves will activate the NAD shock detector more frequently compared to lower order result. Detailed examination of the wave interaction zone reveals that the hybrid scheme's $\mathcal{P}7$ results exhibit slightly greater dissipation than CWENOZ method. The $\mathcal{P}7$ CWENOZ scheme better preserves small-scale vortical structures, as evidenced by the details of vortices on slipstream. With the excessive use of the MUSCL scheme, the increased numerical dissipation leads to marginally fewer vortices being resolved by the hybrid scheme relative to the CWENOZ method in high-order ($\mathcal{P}7$) simulations. To address this problem, modifications to the hybrid criteria that prioritise CWENOZ and high-order linear scheme utilisation and lower the excess use of the MUSCL scheme at higher polynomial orders, can improve the performance of the Hybrid scheme as seen from the results found on Fig.\ref{fig:DMR_DifferentHybridSetting}.

\newpage

\clearpage
\subsection{Forward Facing Step}
\begin{figure}[h!]
 \begin{centering}
 
 \captionsetup[subfigure]{width=0.30\textwidth}
  \subfloat[CWENOZ (Coarse Mesh)]
 {\includegraphics[angle=0,width=0.32\textwidth,trim={1cm 20cm 1cm 20cm},clip]{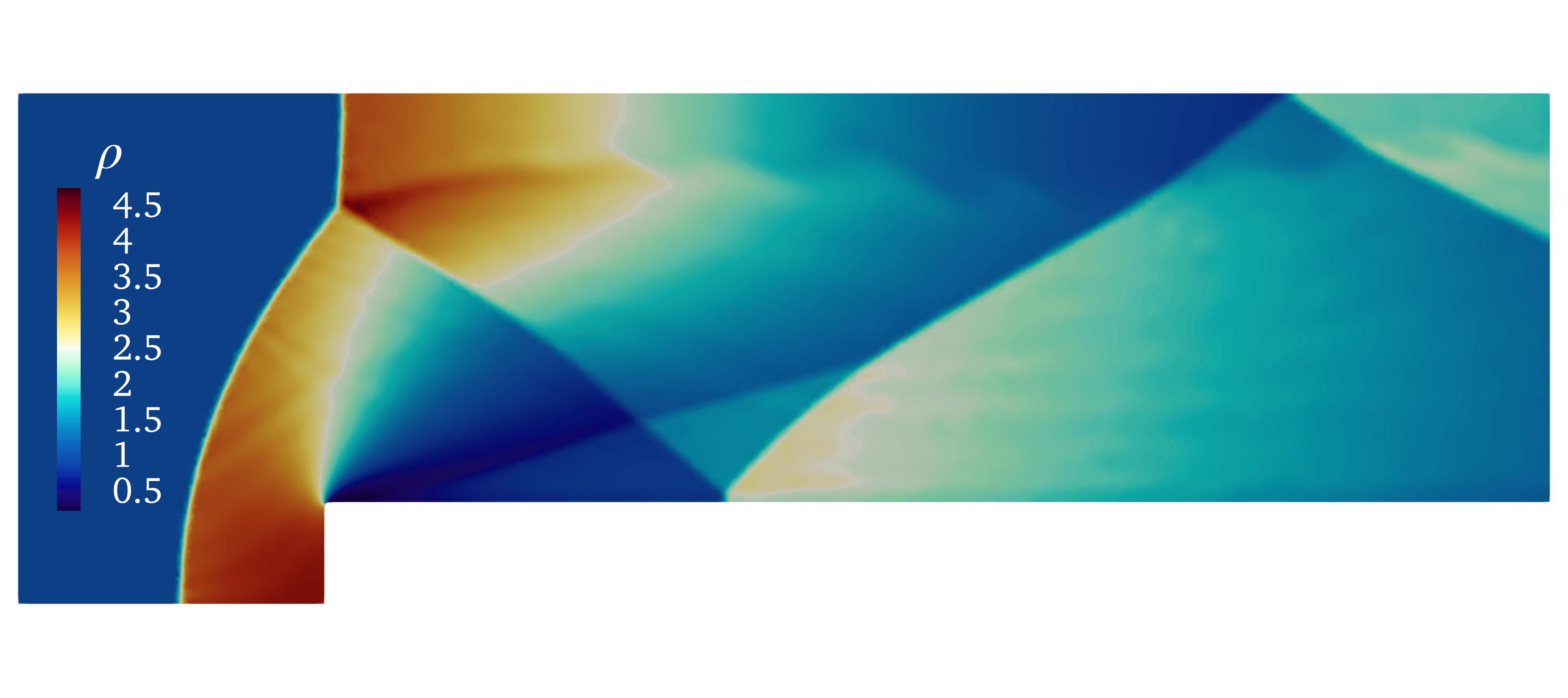}}
 \subfloat[Hybrid (Coarse Mesh)]
 {\includegraphics[angle=0,width=0.32\textwidth,trim={1.5cm 20cm 1.5cm 20cm},clip]{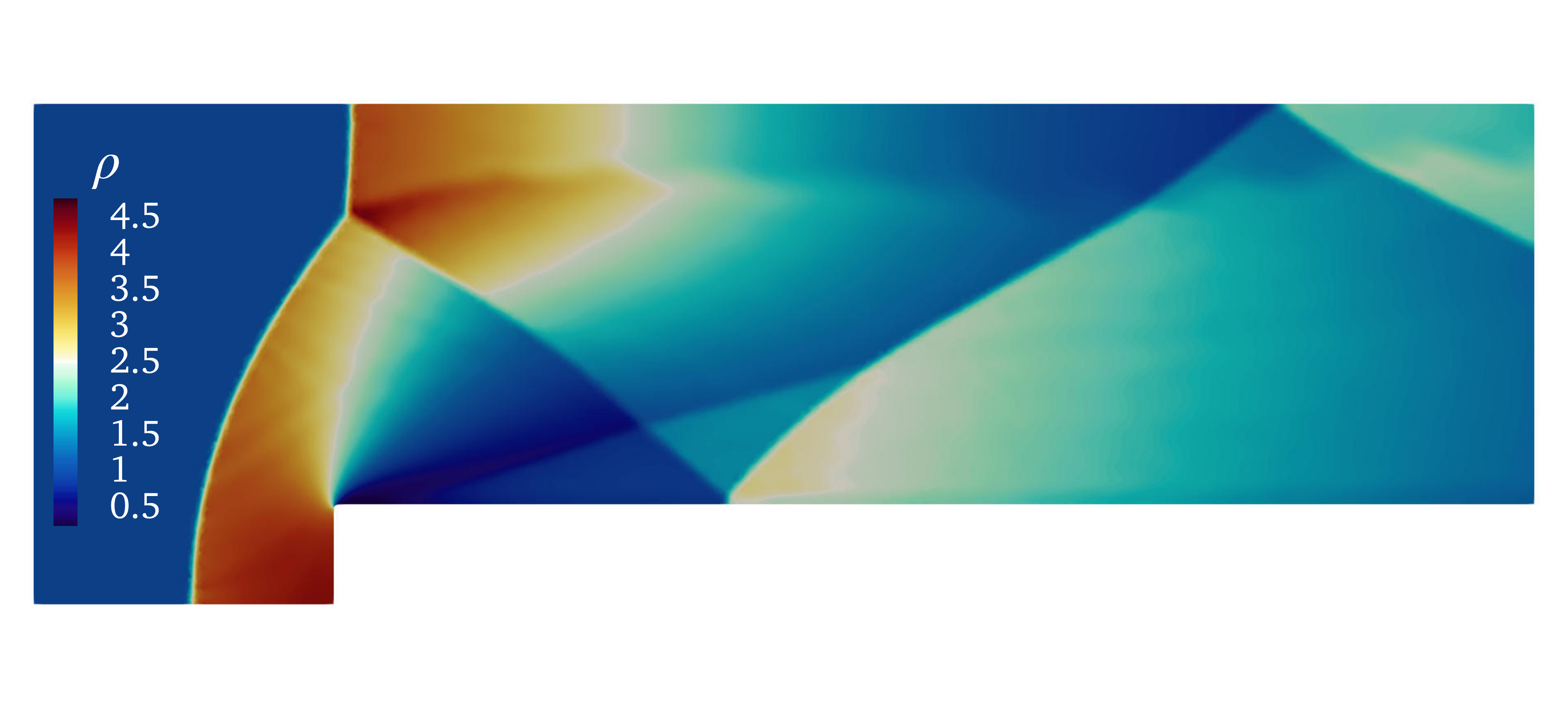}}
 \subfloat[Hybrid (Coarse Mesh)]
 {\includegraphics[angle=0,width=0.32\textwidth,trim={1cm 20cm 1cm 20cm},clip]{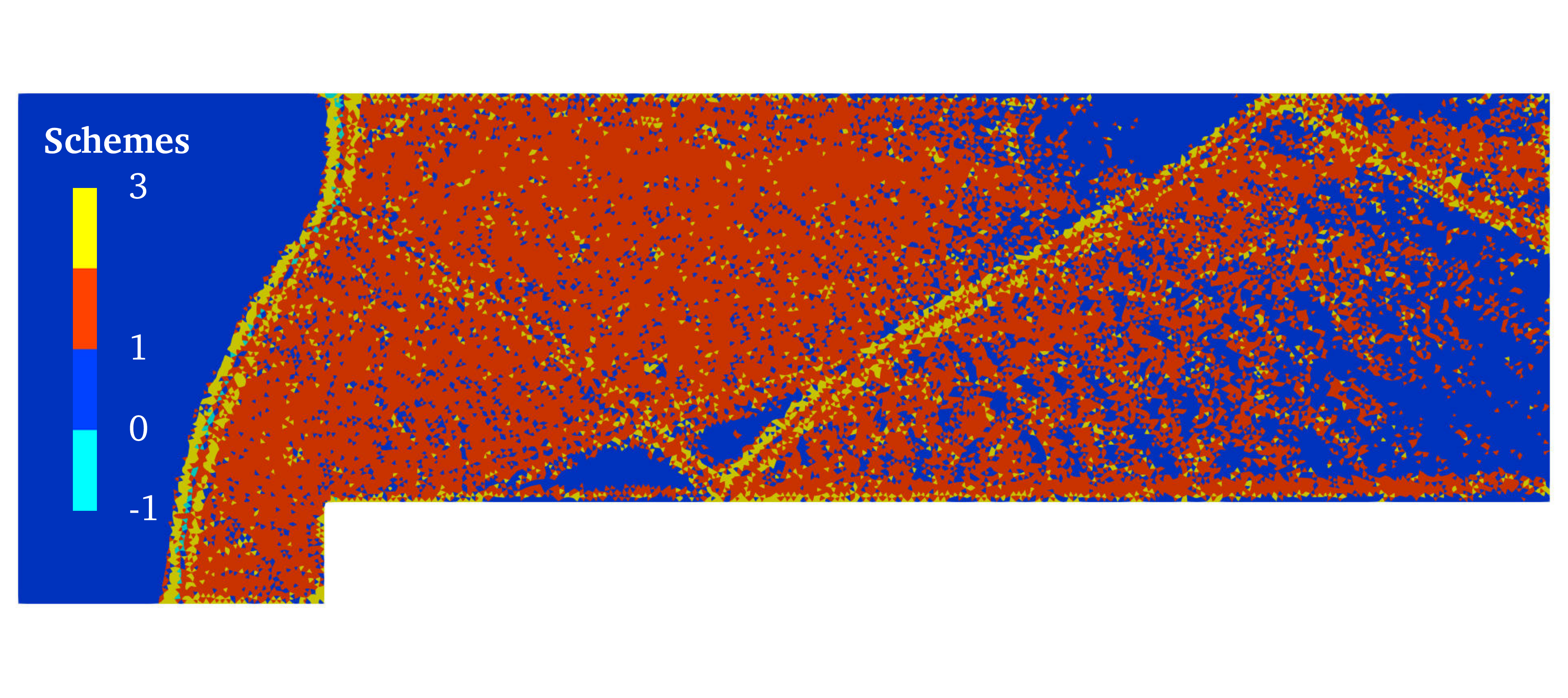}}

   \subfloat[CWENOZ (Medium Mesh)]
  {\includegraphics[angle=0,width=0.32\textwidth,trim={1cm 20cm 1cm 20cm},clip]{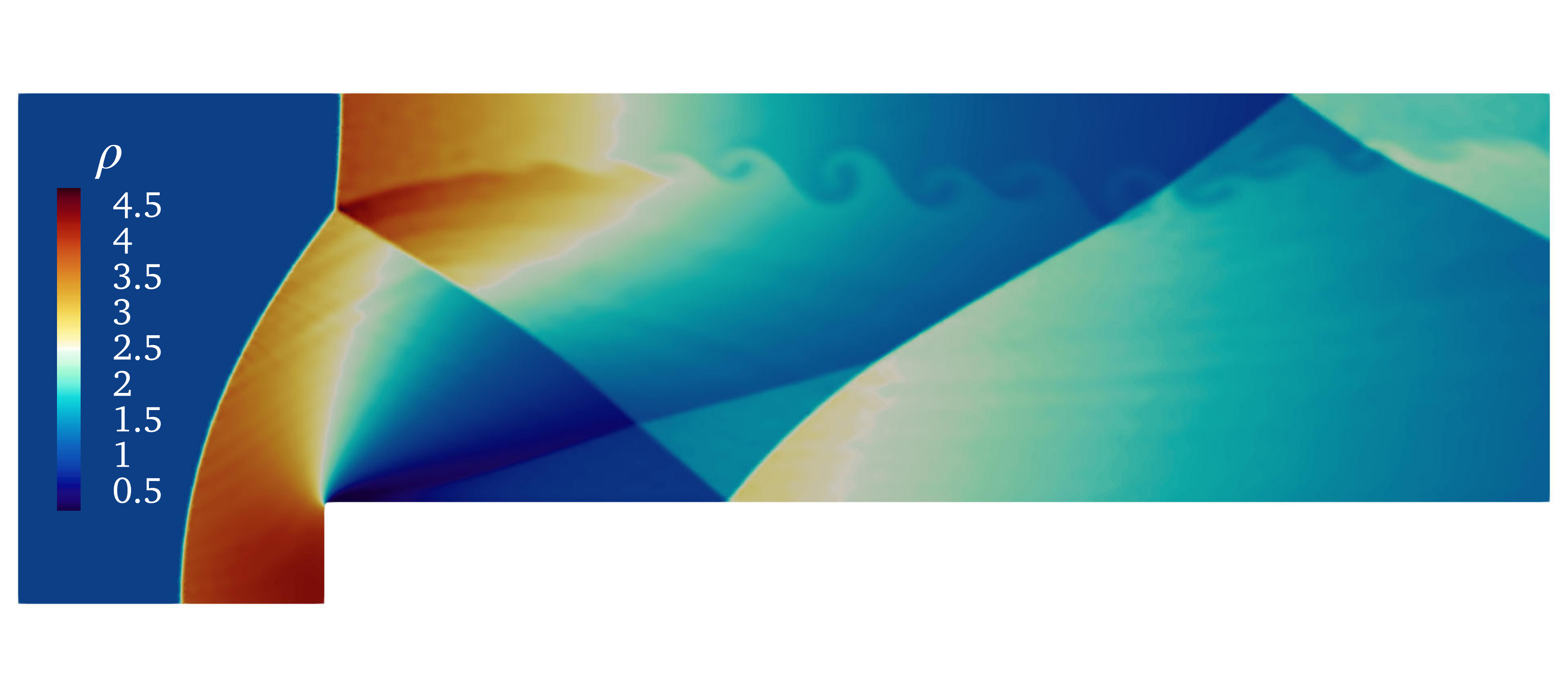}}
 \subfloat[Hybrid (Medium Mesh)]
 {\includegraphics[angle=0,width=0.32\textwidth,trim={1.5cm 20cm 1.5cm 20cm},clip]{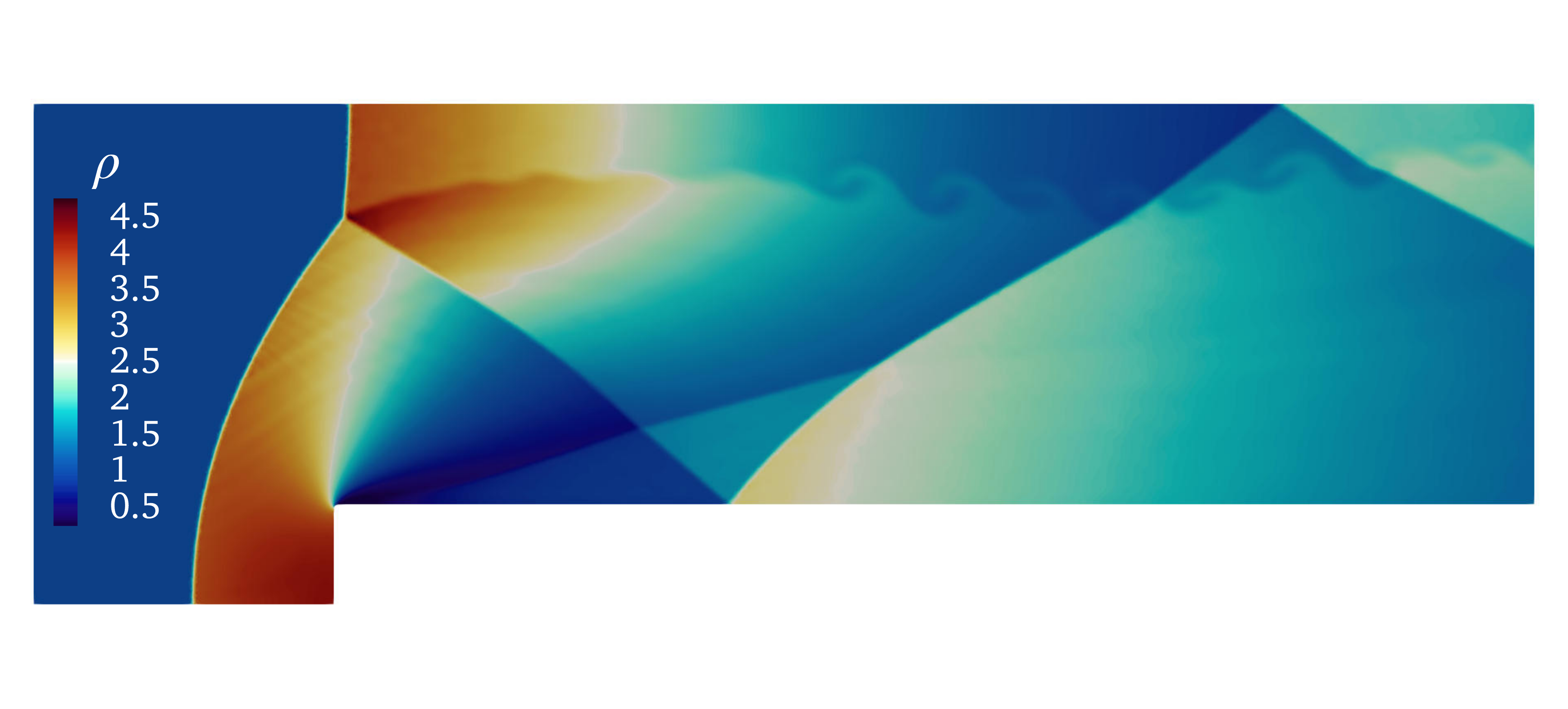}}
 \subfloat[Hybrid (Medium Mesh)]
 {\includegraphics[angle=0,width=0.32\textwidth,trim={1cm 20cm 1cm 20cm},clip]{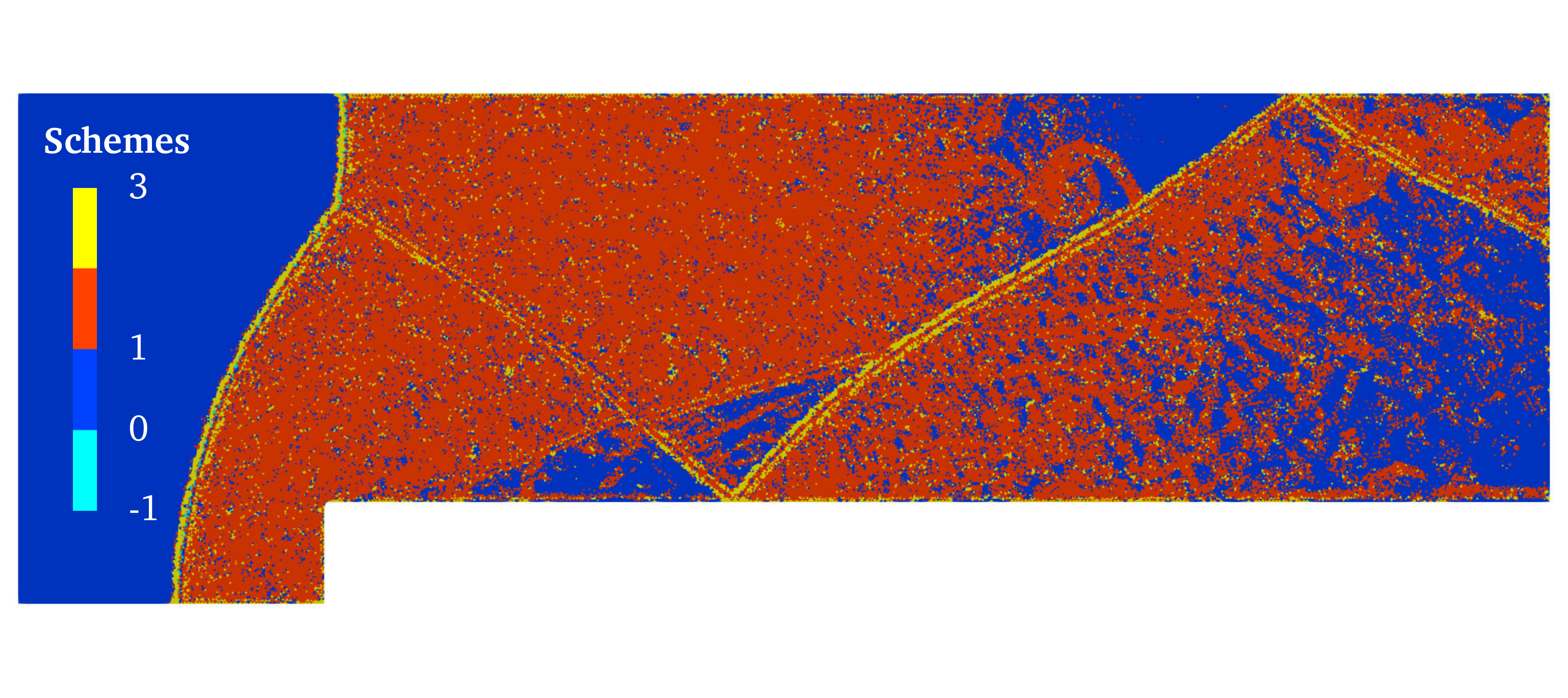}}

   \subfloat[CWENOZ (Fine Mesh)]
  {\includegraphics[angle=0,width=0.32\textwidth,trim={1cm 20cm 1cm 20cm},clip]{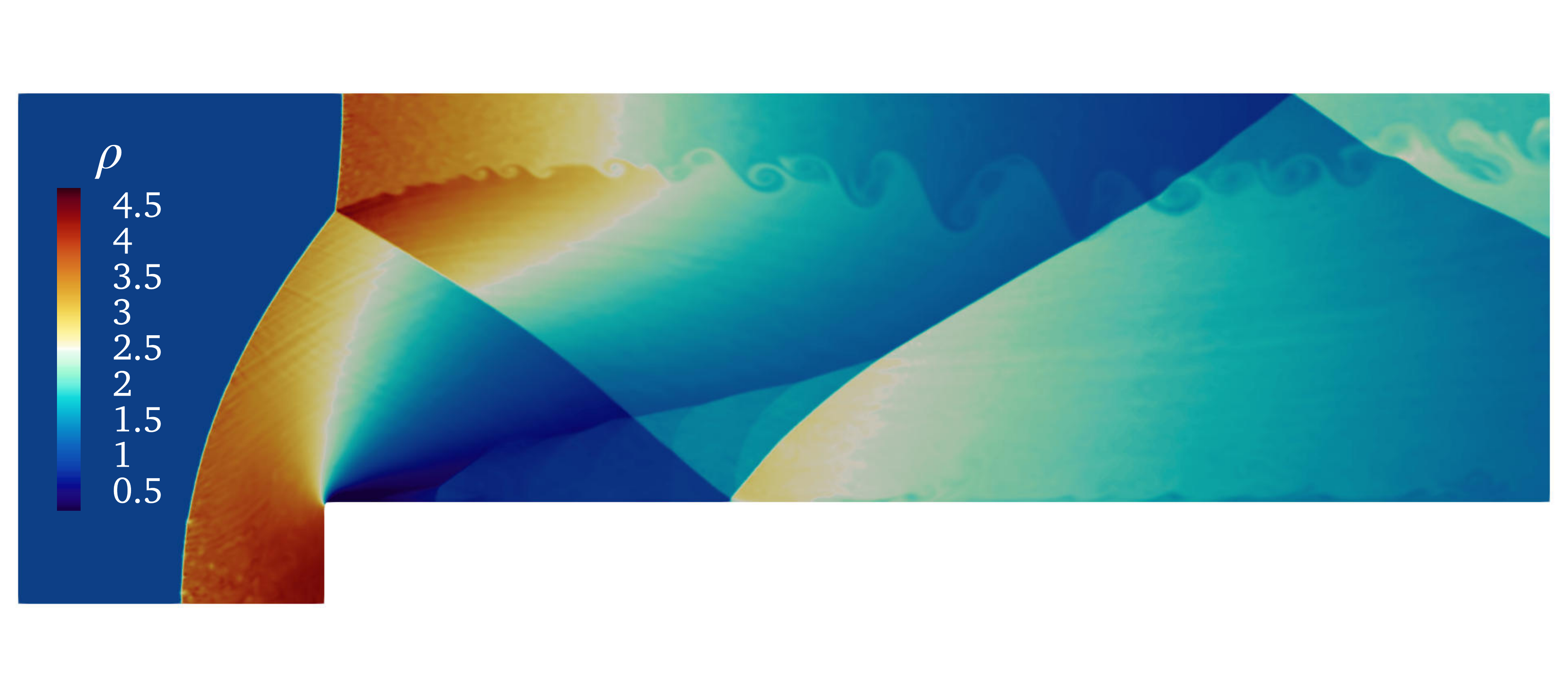}}
 \subfloat[Hybrid (Fine Mesh)]
 {\includegraphics[angle=0,width=0.32\textwidth,trim={1.5cm 20cm 1.5cm 20cm},clip]{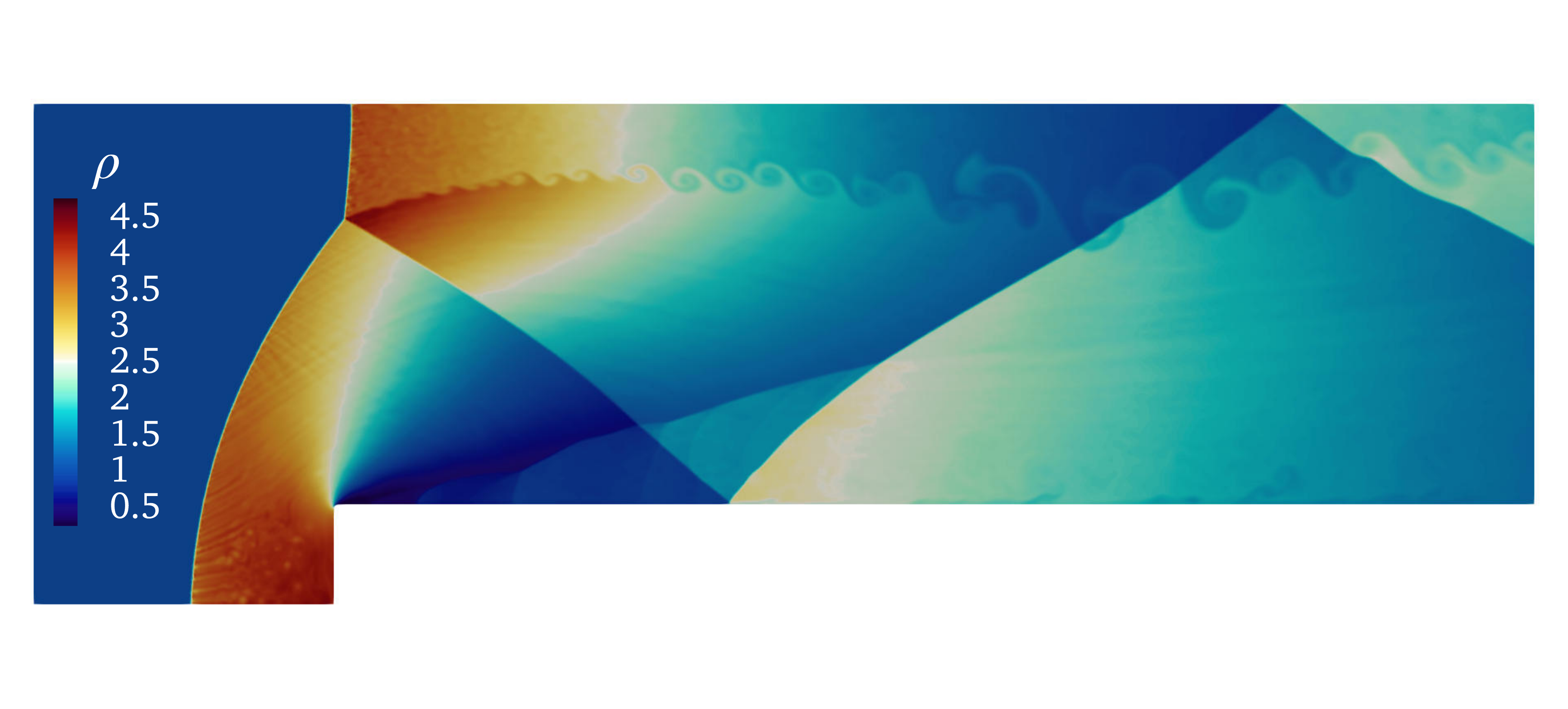}}
 \subfloat[Hybrid (Fine Mesh)]
 {\includegraphics[angle=0,width=0.32\textwidth,trim={1cm 20cm 1cm 20cm},clip]{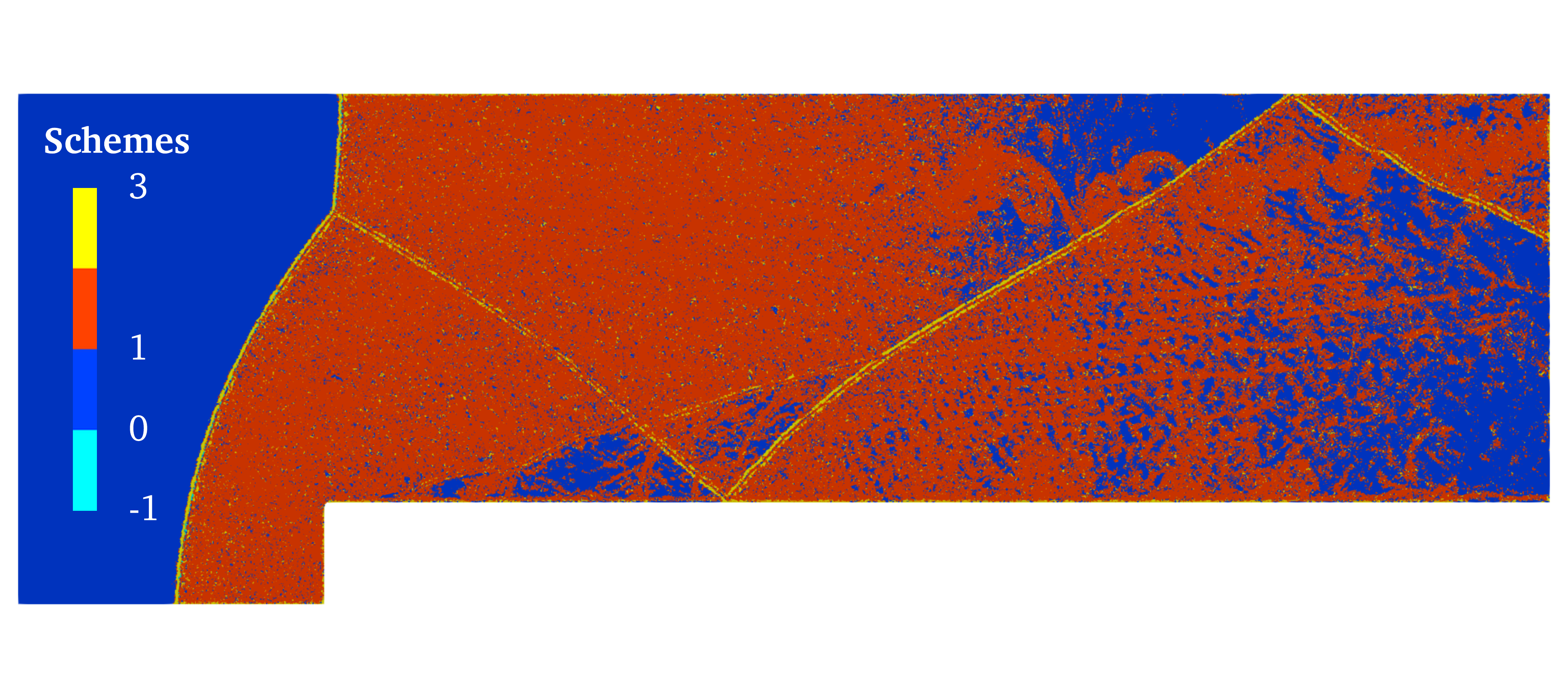}}

 \par\end{centering}\caption{Density contours computed using the CWENOZ scheme (\textit{left}); Density contours obtained via the Hybrid scheme (\textit{middle}); numerical schemes contours in the Hybrid method (\textit{right}). The scheme legend illustrates: yellow regions (Scheme 3) denote cells employing the $2nd-$ order MUSCL scheme, red zones (Scheme 1) indicate high-order CWENOZ scheme, deep blue areas (Scheme 0) represent \textcolor{black}{{high-order}} linear scheme, and light blue sections (Scheme -1) correspond to $1st-order$ upwind treatment.}\label{fig:FFS_DifferentMeshes}\end{figure}

The 2D Forward-Facing Step (FFS) problem, originally introduced by Woodward and Colella \citep{Woodward1984115}, is employed to evaluate the capability of numerical methods in handling complex shock interactions. The computational domain is defined as $x \in [0, 3]$, $y \in [0, 1]$, with a step located at $x = 0.6$ and a height of 0.2 units. Following the approach of Nazarov and Larcher \citep{NAZAROV2017128}, the corner of the step is smoothed with a curvature of radius 0.003 to eliminate non-physical effects arising from the singularity point. The left boundary is specified as a supersonic inflow, while the right boundary is set as a supersonic outflow. Reflective boundary conditions are applied to both the top and bottom boundaries. The initial condition consists of a uniform ideal gas with parameters $[\rho, p, u, v] = [1.4, 1, 3, 0]$. The total simulation time is $t = 4$, and the CFL number for all simulations is fixed at 0.6 and an HLL Riemann solver is used due to its robustness. This test is conducted on three levels of resolution of triangular meshes $h = 1/70$, $1/140$, and $1/280$ with the $5th$ order of accuracy. The corresponding results are shown in Fig. \ref{fig:FFS_DifferentMeshes}.

The results across different level of meshes successfully capture all key flow features, including the incident shock, Mach stem, slipstream, and reflected shocks. A comparison between the CWENOZ and Hybrid schemes reveals that, at coarse and medium resolutions, the Hybrid method exhibits slightly higher numerical dissipation compared to the CWENOZ scheme. However, in the fine mesh, both methods demonstrate excellent capability in resolving instability vortices along the slipstream. Notably, the Hybrid scheme shows a slight advantage in capturing the Kelvin Helmholtz instability at this resolution. The distribution of schemes employed within the Hybrid method remains largely consistent across mesh refinements. As the mesh becomes finer and more flow structures are resolved, the utilisation of the CWENOZ scheme increases. Meanwhile, the MUSCL scheme is predominantly applied in regions with strong discontinuities, such as shock waves precisely as intended in the design of the Hybrid approach.

\newpage
\subsection{2D Riemann Problem}
The two-dimensional Riemann problem introduced by Schulz et al.~\citep{RINNE199311} is used to evaluate the performance of the proposed hybrid scheme. The computational domain is defined as $x \in [-0.5, 0.5]$, $y \in [-0.5, 0.5]$. The domain is divided into four quadrants, each initialised with a distinct fluid state. The initial condition can thus be expressed as:

\begin{equation}\label{eq:2DRiemannProblem}
(\rho, u, v, p) =
\begin{cases}
(0.5323,1.206,0,0.3)      & \text{for } x \leq 0, \; y \geq 0 \\
(1.5,0,0,1.5)             & \text{for } x \geq 0, \; y \geq 0 \\
(0.138,1.206,1.206,0.029) & \text{for } x \leq 0, \; y \leq 0 \\
(0.5323,0,1.206,0.3)      & \text{for } x \geq 0, \; y \leq 0
\end{cases}
\end{equation}

The mesh configurations used in this study follow those described by Tsoutsanis \citep{TSOUTSANIS201869}. Three different levels of mesh refinement are considered. In the refined region corresponding to the lower-left quadrant, the grid resolutions are $h = 1/200$, $1/400$, and $1/800$, respectively. For the remainder of the domain, a coarser resolution of $h = 1/30$ is used. The CFL number is set to 0.6 for all simulations and the HLL Riemann solver is used. The total simulation $t=1$. Different hybrid mesh resolutions were tested on this benchmark case. We also explored several hybrid schemes, such as MUSCL-CWENOZ, MUSCL-Linear, and CWENOZ-Linear, by adjusting the hybrid settings.

\begin{figure}[htbp]
 \begin{centering}
 
 \captionsetup[subfigure]{width=0.30\textwidth}
  \subfloat[CWENOZ (Coarse Mesh)]
 {\includegraphics[angle=0,width=0.33\textwidth,trim={9.8cm 9.8cm 9.8cm 9.8cm},clip]{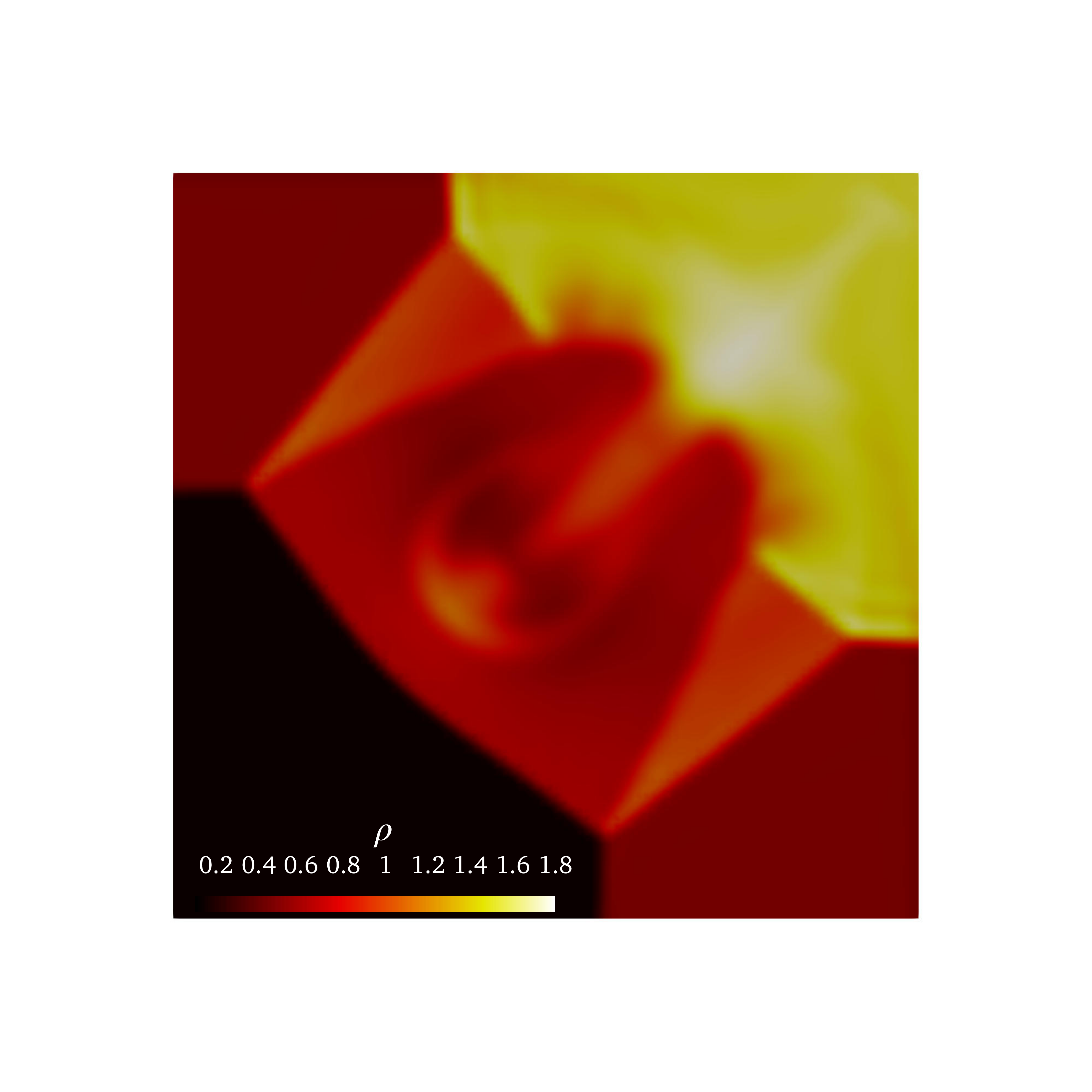}}
 \subfloat[Hybrid (Coarse Mesh)]
 {\includegraphics[angle=0,width=0.33\textwidth,trim={9.8cm 9.8cm 9.8cm 9.8cm},clip]{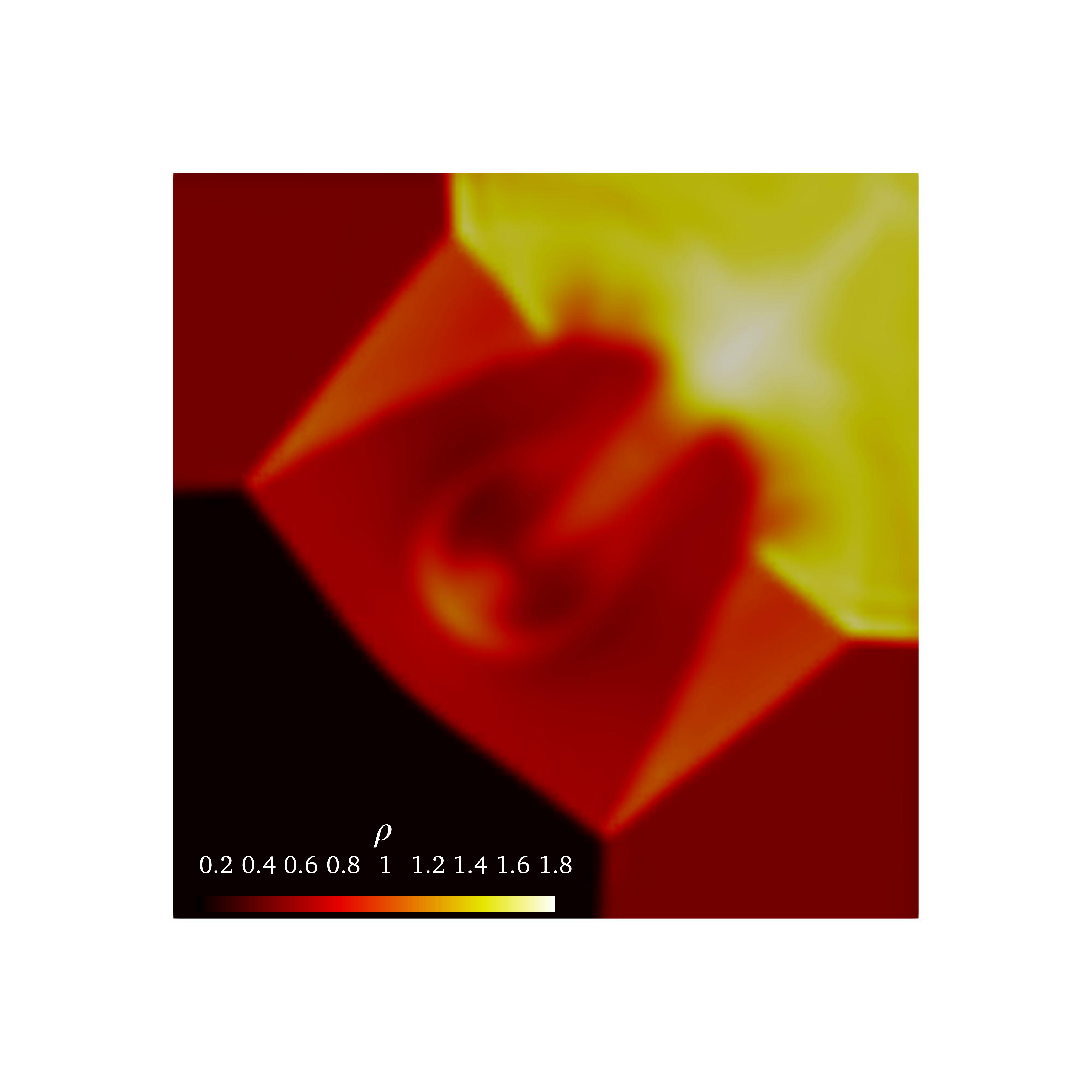}}
 \subfloat[Hybrid (Coarse Mesh)]
 {\includegraphics[angle=0,width=0.33\textwidth,trim={9cm 9cm 9cm 9cm},clip]{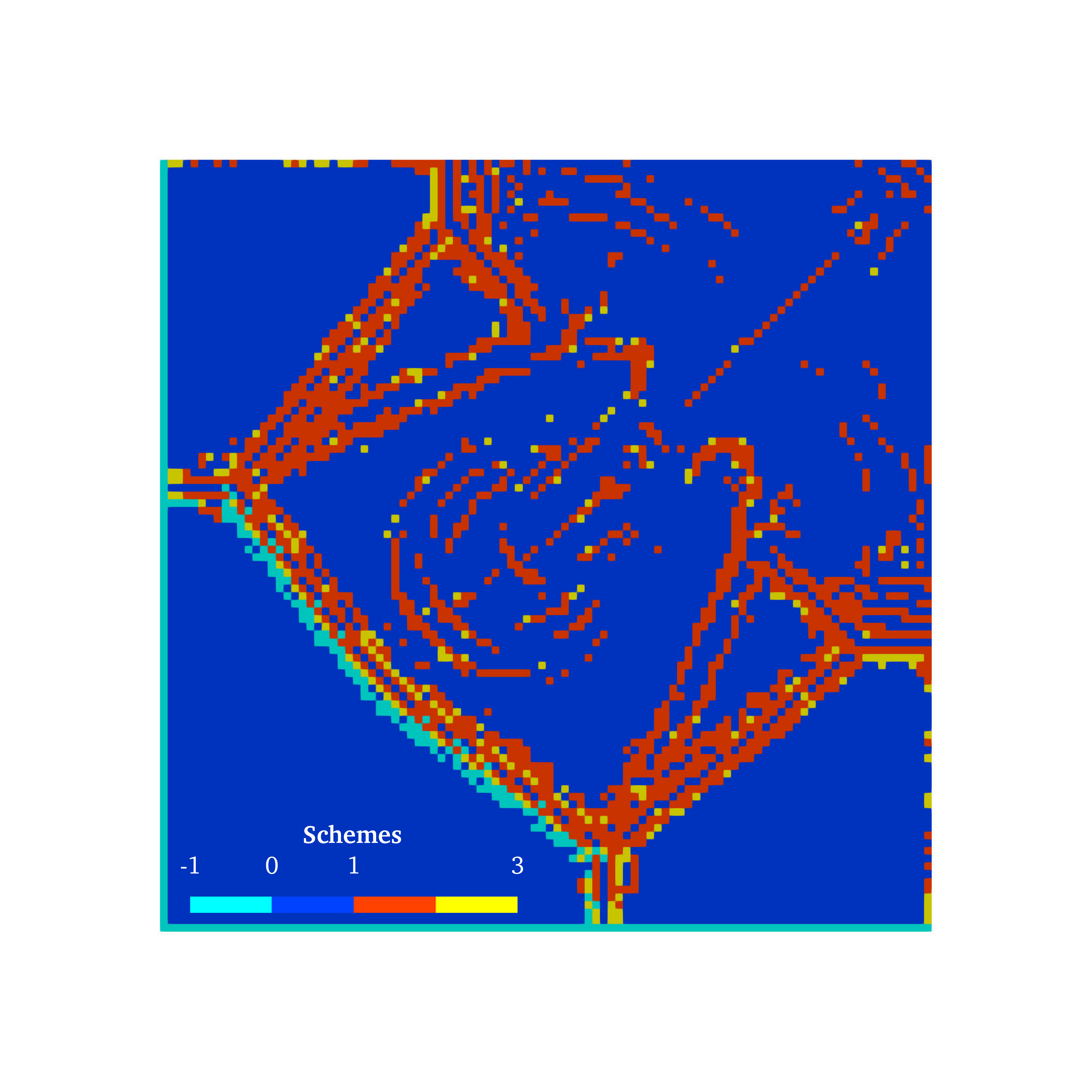}}

   \subfloat[CWENOZ (Medium Mesh)]
 {\includegraphics[angle=0,width=0.33\textwidth,trim={9.8cm 9.8cm 9.8cm 9.8cm},clip]{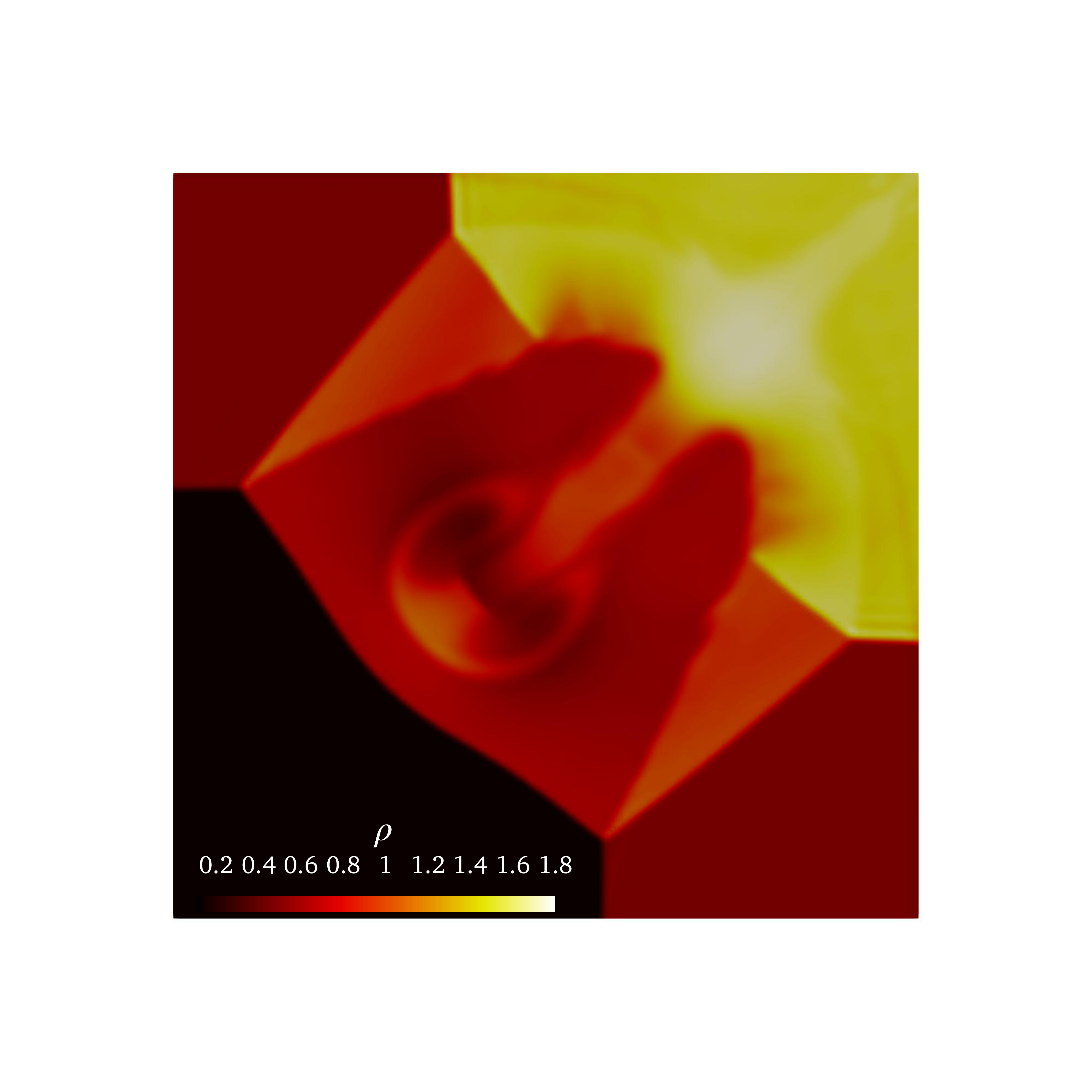}}
 \subfloat[Hybrid (Medium Mesh)]
 {\includegraphics[angle=0,width=0.33\textwidth,trim={9.8cm 9.8cm 9.8cm 9.8cm},clip]{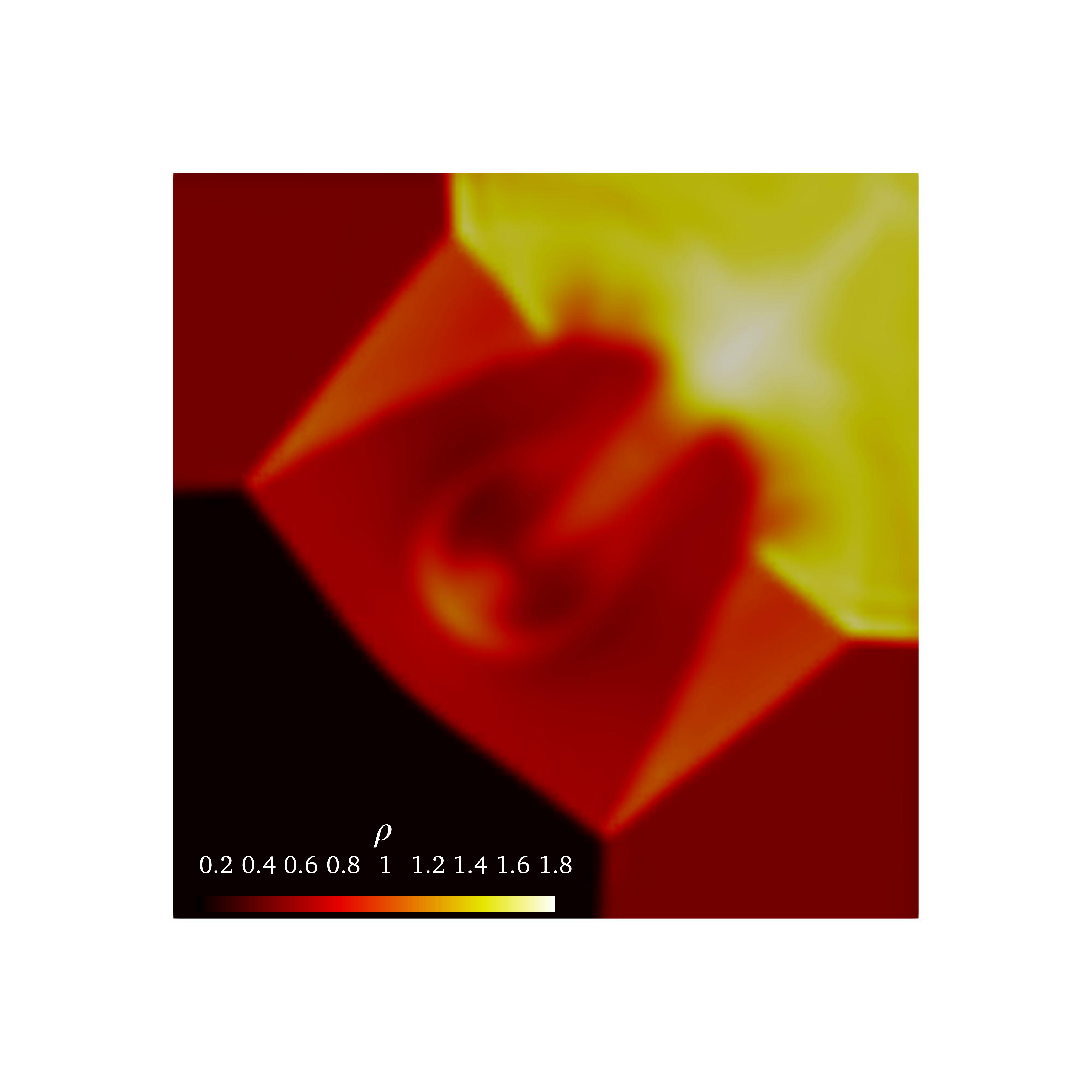}}
 \subfloat[Hybrid (Medium Mesh)]
 {\includegraphics[angle=0,width=0.33\textwidth,trim={9cm 9cm 9cm 9cm},clip]{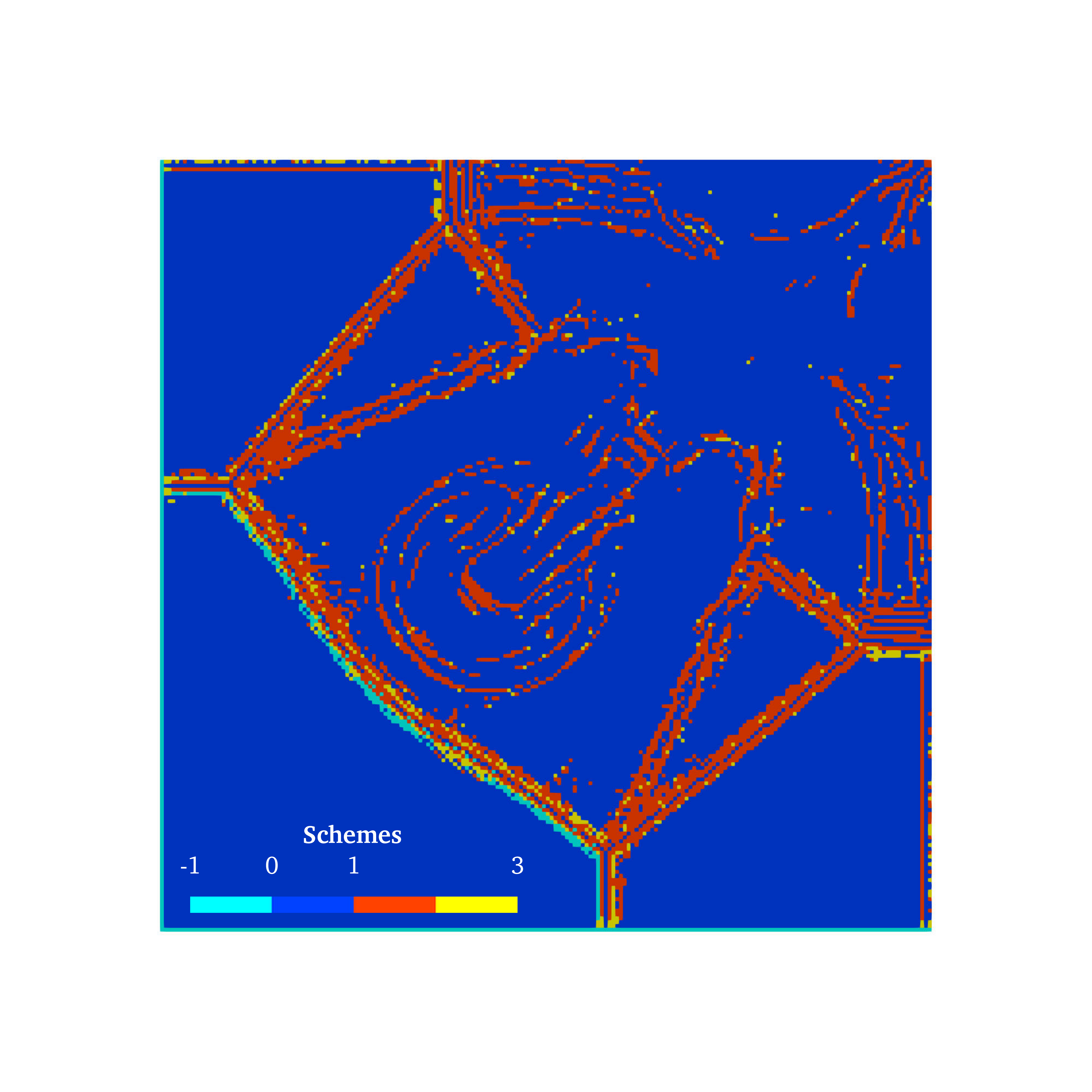}}

    \subfloat[CWENOZ (Fine Mesh)]
 {\includegraphics[angle=0,width=0.33\textwidth,trim={9.8cm 9.8cm 9.8cm 9.8cm},clip]{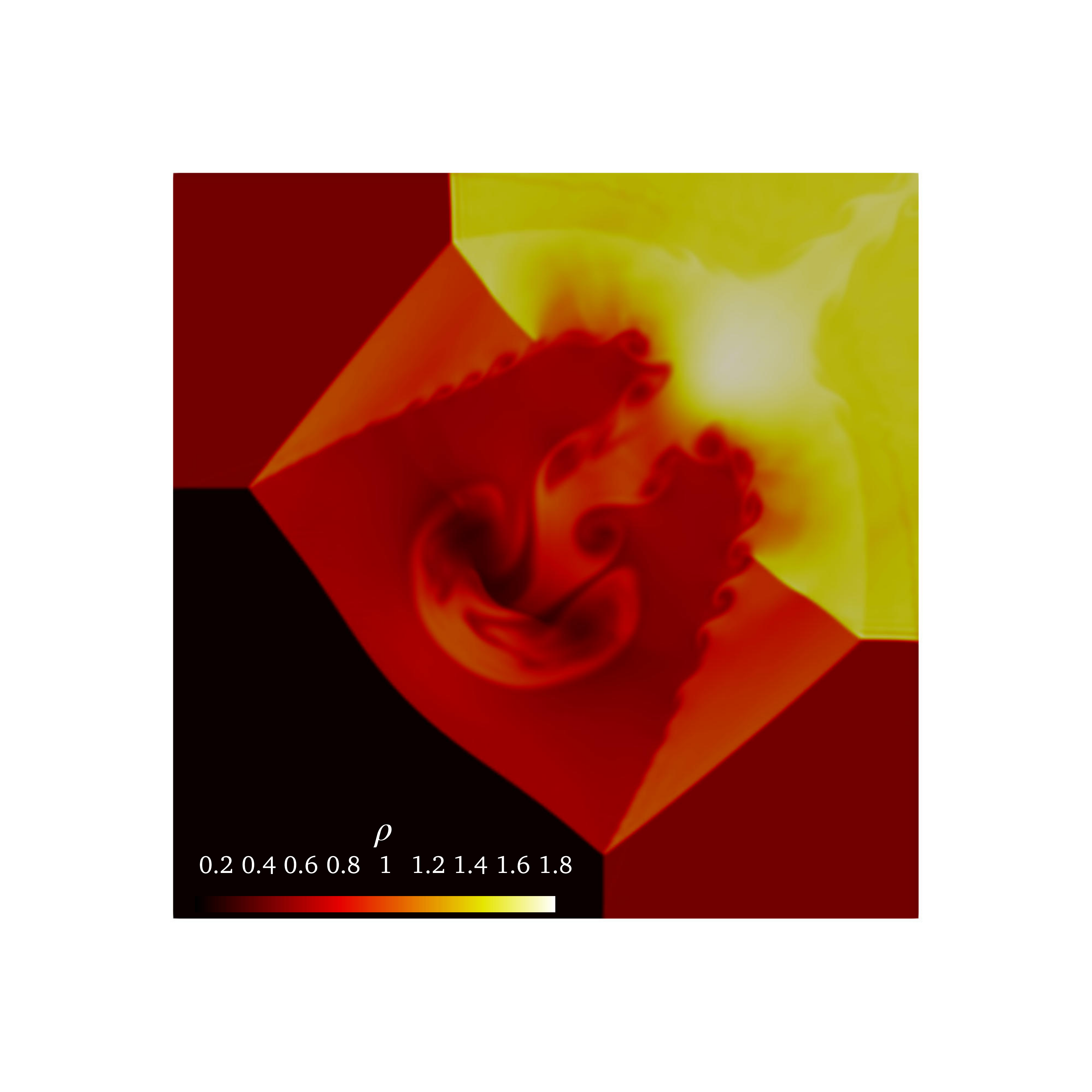}}
 \subfloat[Hybrid (Fine Mesh)]
 {\includegraphics[angle=0,width=0.33\textwidth,trim={9.8cm 9.8cm 9.8cm 9.8cm},clip]{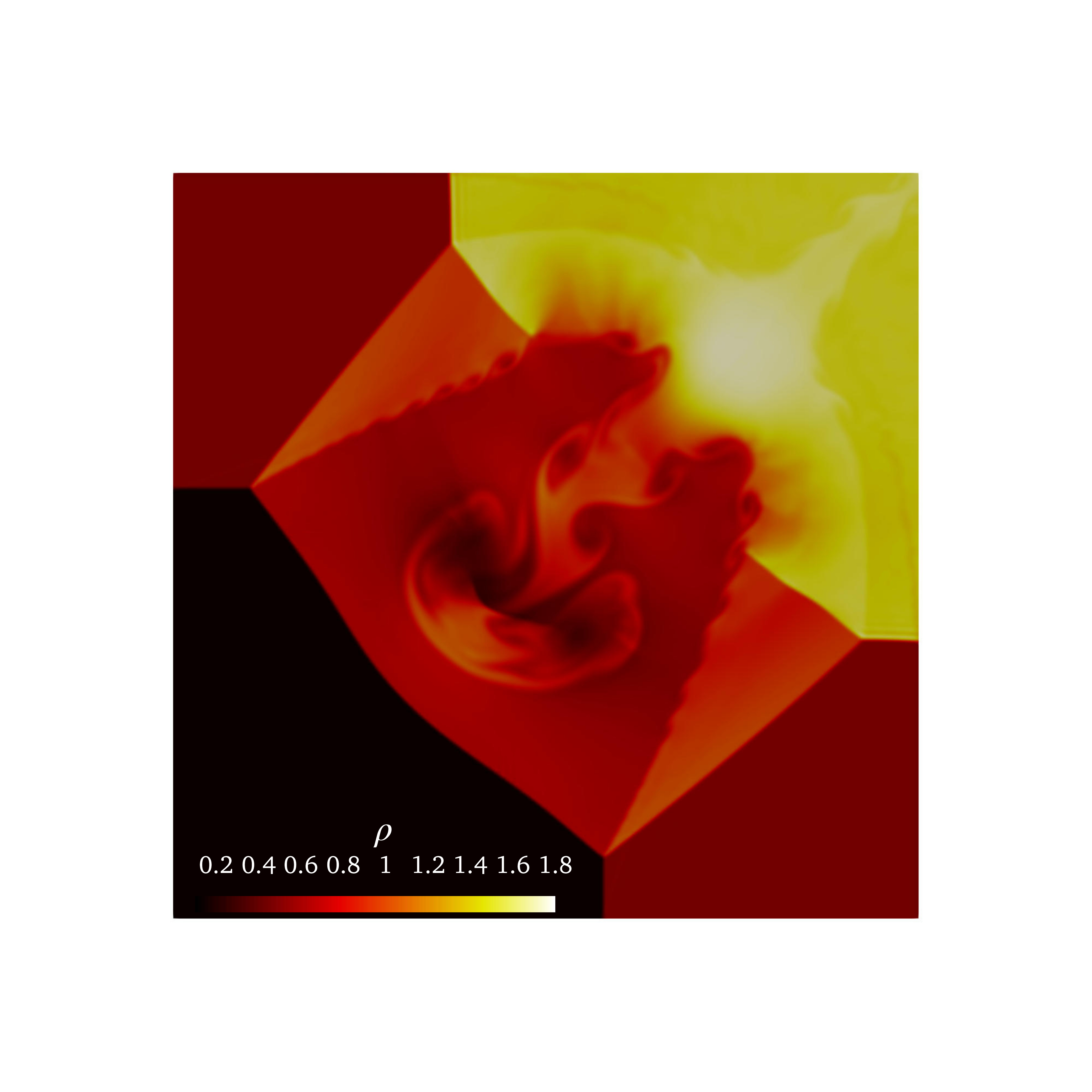}}
 \subfloat[Hybrid (Fine Mesh)]
 {\includegraphics[angle=0,width=0.33\textwidth,trim={9cm 9cm 9cm 9cm},clip]{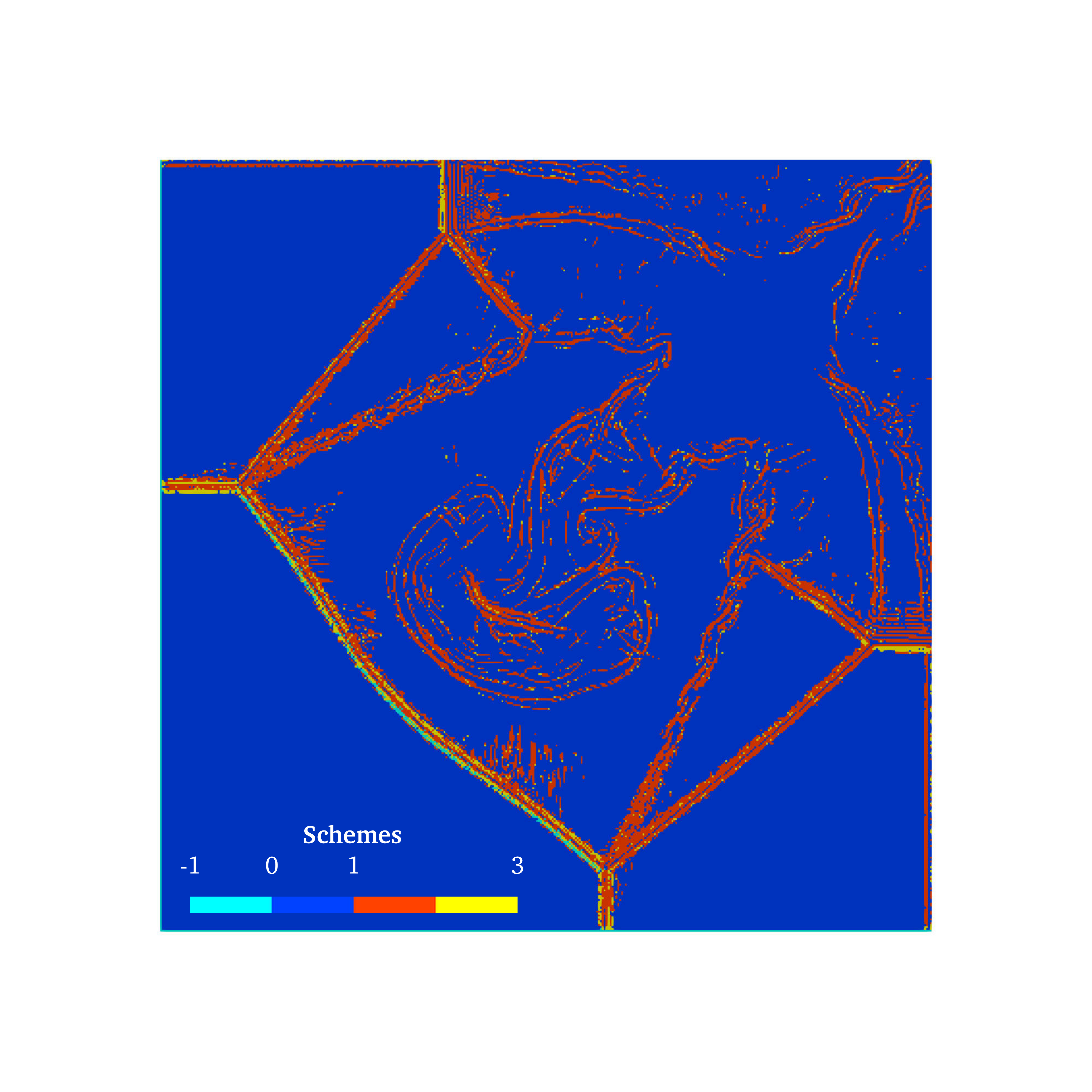}}
 \par\end{centering}\caption{Density contours for CWENOZ scheme \textit({left}) , Density contours for Hybrid scheme \textit({middle}), and Cells using different schemes(\textit{right}),}\label{fig:RiemannDifferentMeshes}\end{figure}

From the results obtained shown in Fig. \ref{fig:RiemannDifferentMeshes}, it can be found that there is no loss of order accuracy in the coarse mesh. All levels of meshes successfully capture all the solution patterns. The scheme distribution contours show that the DMP indicator partitions the flow field into three intended categories: smooth regions, weakly non-smooth regions, and discontinuities. This partitioning guides the selective application of reconstruction methods. High-dissipation schemes are applied only in cells that intersect true discontinuities, which ensures robust shock capturing without unnecessary diffusion of smooth features. On the fine mesh, the default hybrid setting nevertheless fails to resolve the roll ups of the Kelvin Helmholtz instability as well as the CWENOZ baseline. The loss of small scale structure arises from the additional TVD safeguard. By design, this safeguard triggers a first order upwind fallback in pockets flagged by the indicator, and the added dissipation in these zones damps emerging vortical eddies. The behaviour therefore reflects a trade-off between robustness and resolution. Relaxing the TVD or DMP thresholds, or reducing the extent of the first order fallback region (for example by expanding the CWENOZ domain through the NAD margins $\delta_m$ and $\delta_w$), recovers finer shear layer details at the cost of a higher risk of spurious oscillations. In contrast, more conservative settings enhance stability but tend to smooth delicate flow features. Compared to the CWENOZ scheme, the hybrid configuration under the default setting improves computational efficiency by 41.67\%. We also evaluated several representative configurations, including the MUSCL-CWENOZ hybrid scheme, the MUSCL-linear hybrid scheme, and the CWENOZ-linear scheme.
The results of all these test cases were consistently aligned across configurations.  Among these hybrid configurations, the MUSCL-CWENOZ configuration required the longest runtime, incurring a 118\% increase in computational cost relative to the default setting, while also capturing more small-scale flow structures. Overall, the hybrid method demonstrates excellent robustness across a range of hybrid settings.

\newpage
\begin{figure}[h!]
 \begin{centering}
 
 \captionsetup[subfigure]{width=0.22\textwidth}
  \subfloat[Default Setting]
 {\includegraphics[angle=0,width=0.25\textwidth,trim={9.8cm 9.8cm 9.8cm 9.8cm},clip]{Hybrid_2D_Riemann_Density_Fine.pdf}}
\subfloat[MUSCL CWENOZ]
 {\includegraphics[angle=0,width=0.25\textwidth,trim={9.8cm 9.8cm 9.8cm 9.8cm},clip]{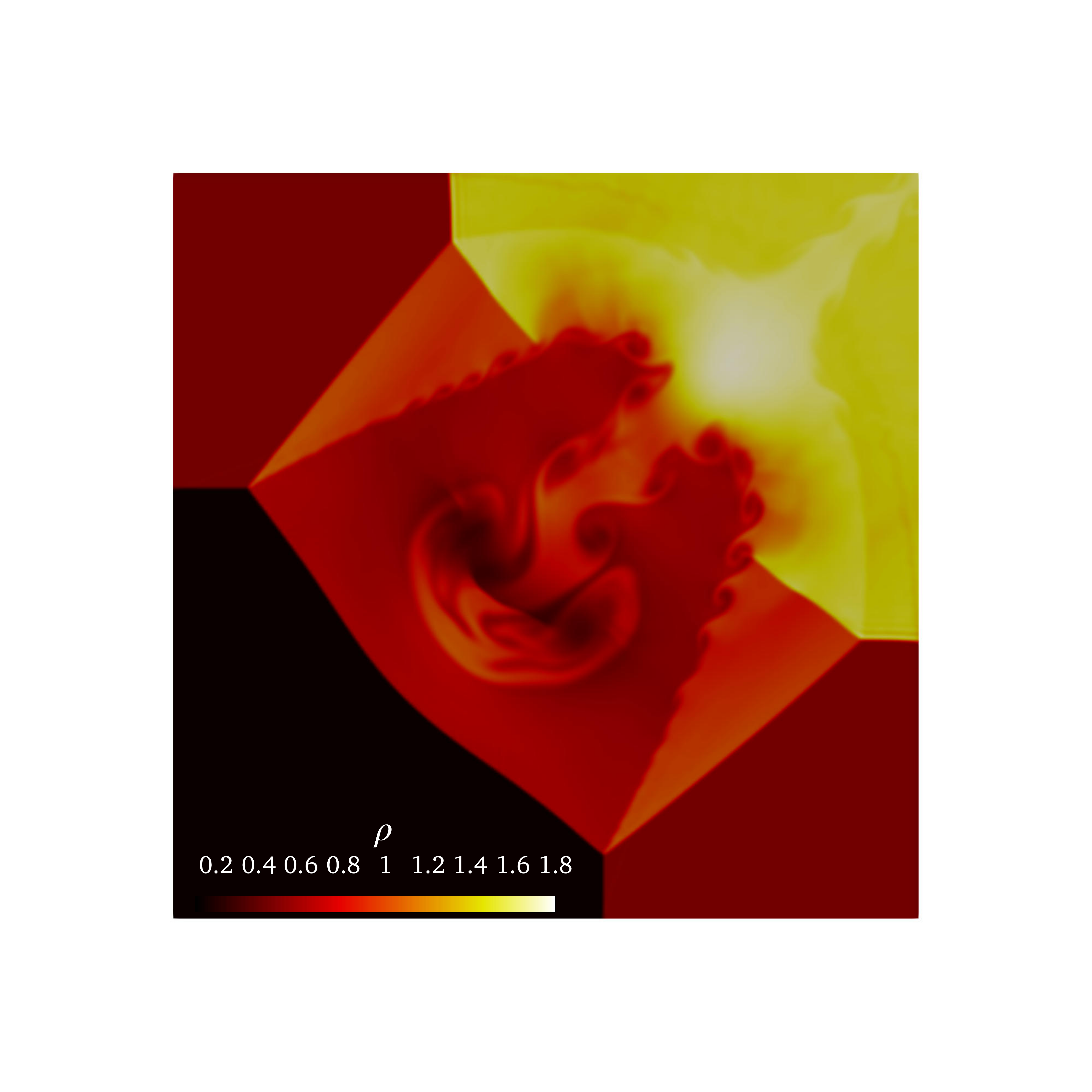}}
\subfloat[MUSCL Linear]
 {\includegraphics[angle=0,width=0.25\textwidth,trim={9.8cm 9.8cm 9.8cm 9.8cm},clip]{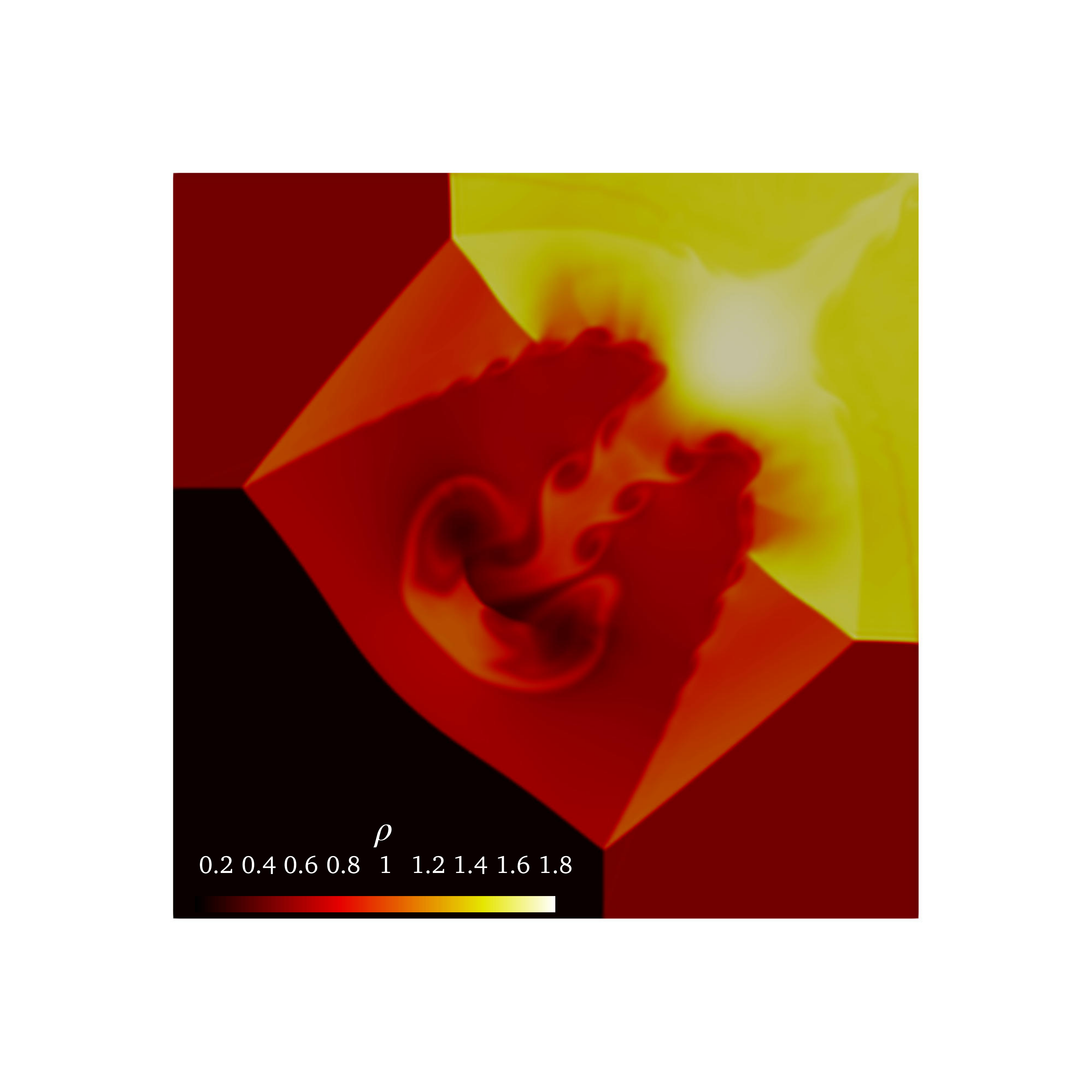}}
\subfloat[CWENOZ Linear]
 {\includegraphics[angle=0,width=0.25\textwidth,trim={9.8cm 9.8cm 9.8cm 9.8cm},clip]{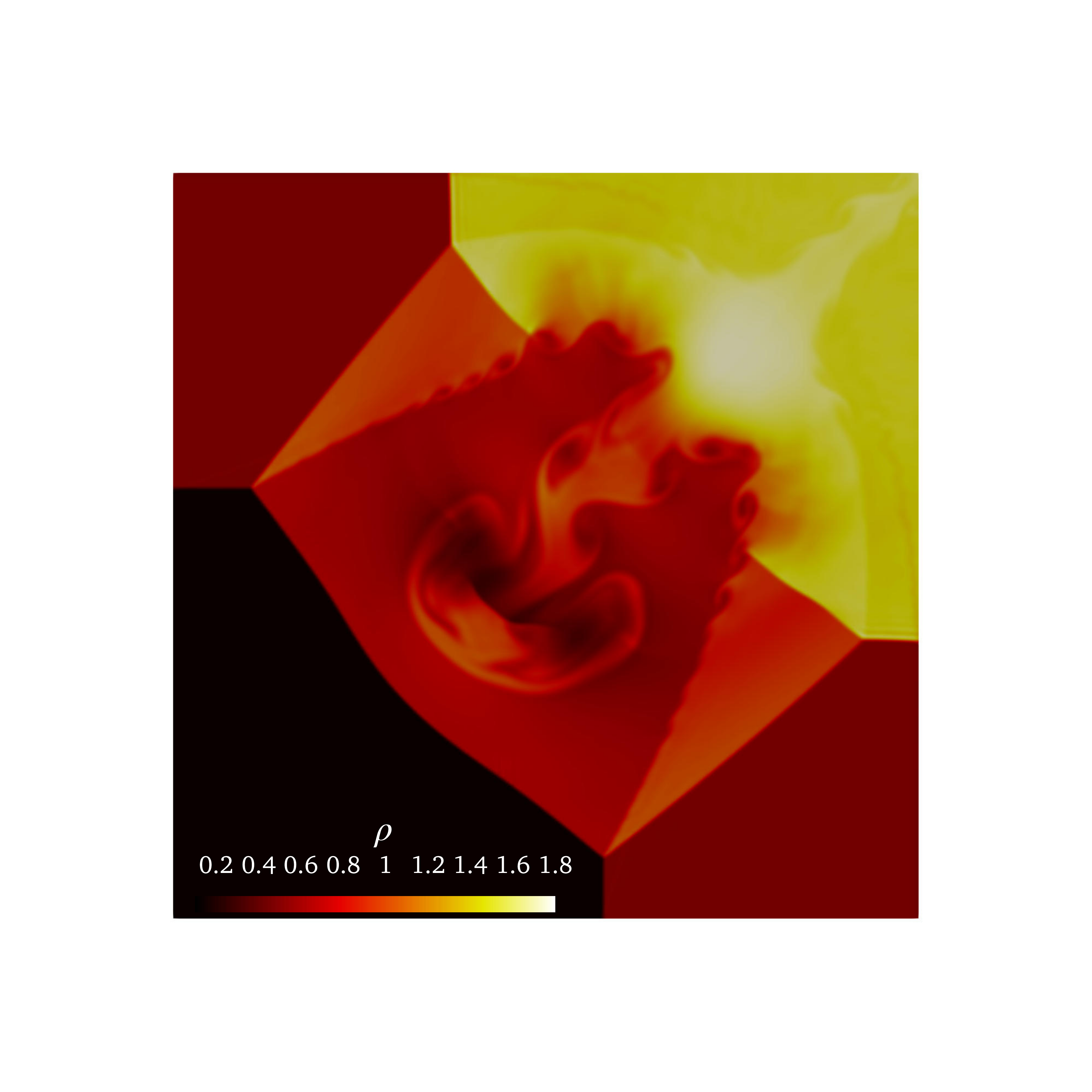}}\\
 
 \subfloat[Default Setting]
 {\includegraphics[angle=0,width=0.25\textwidth,trim={9cm 9cm 9cm 9cm},clip]{Hybrid_2D_Riemann_Schemes_Fine.pdf}} 
 \subfloat[MUSCL CWENOZ]
  {\includegraphics[angle=0,width=0.25\textwidth,trim={9cm 9cm 9cm 9cm},clip]{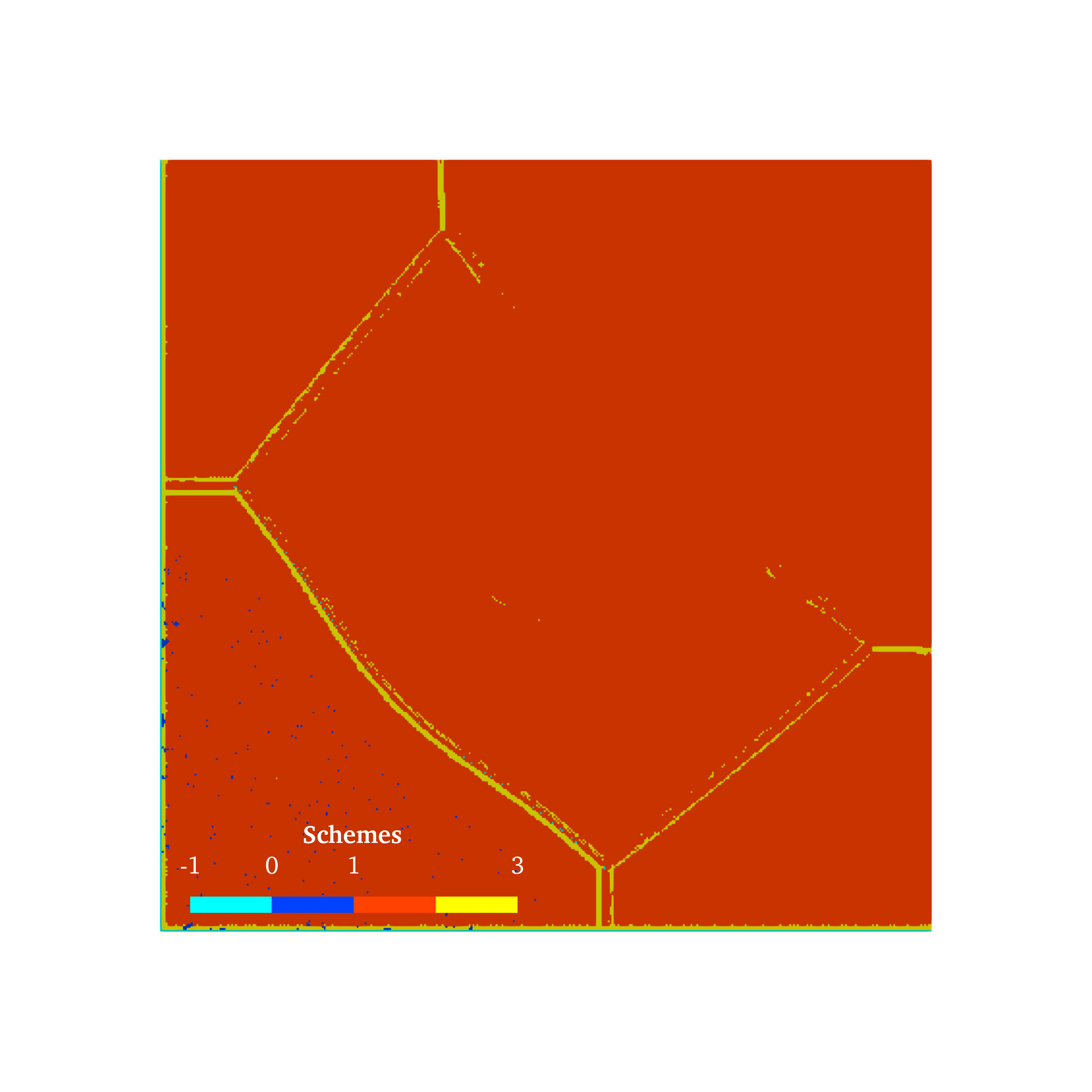}}   
 \subfloat[MUSCL Linear]
  {\includegraphics[angle=0,width=0.25\textwidth,trim={9cm 9cm 9cm 9cm},clip]{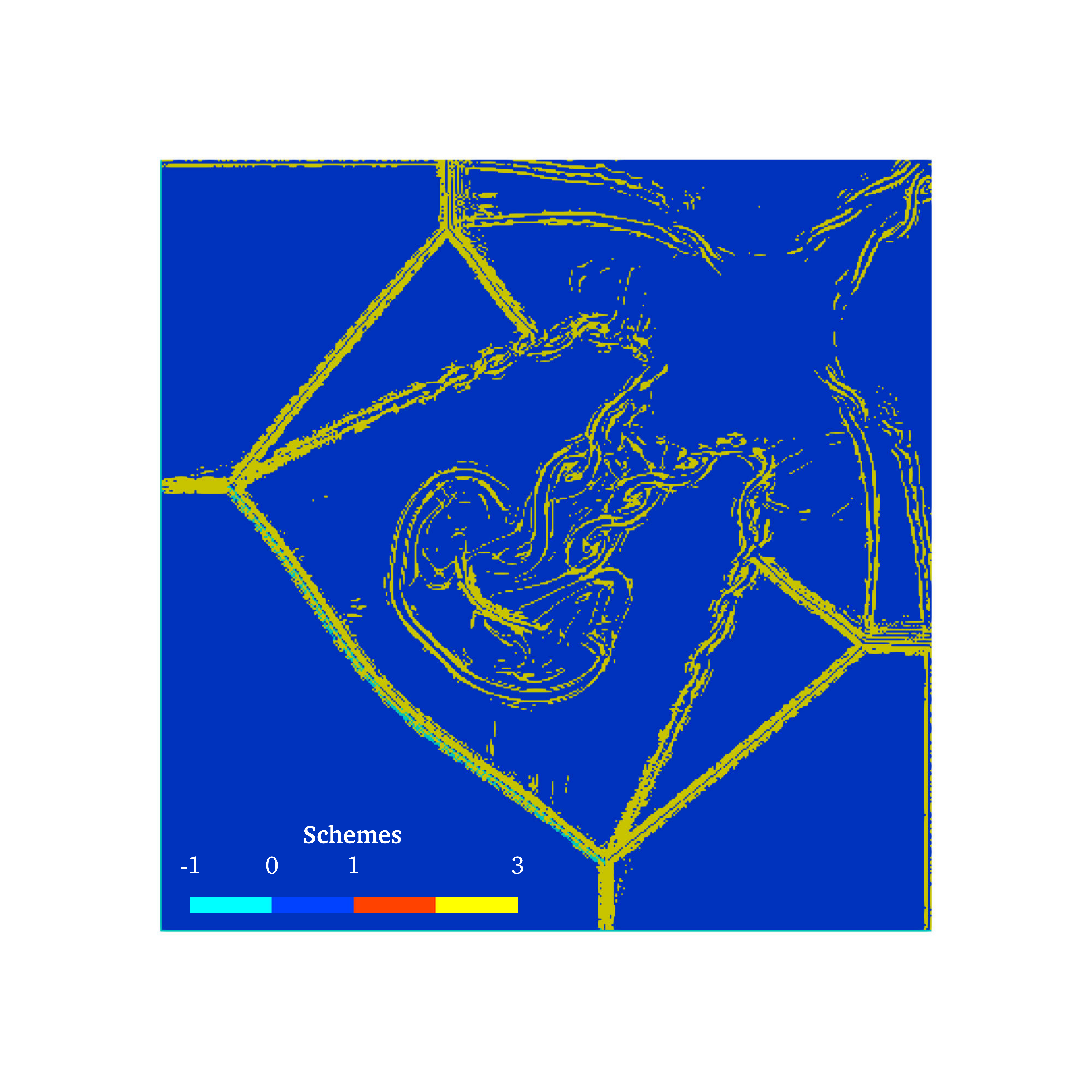}}   
 \subfloat[CWENOZ Linear]
  {\includegraphics[angle=0,width=0.25\textwidth,trim={9cm 9cm 9cm 9cm},clip]{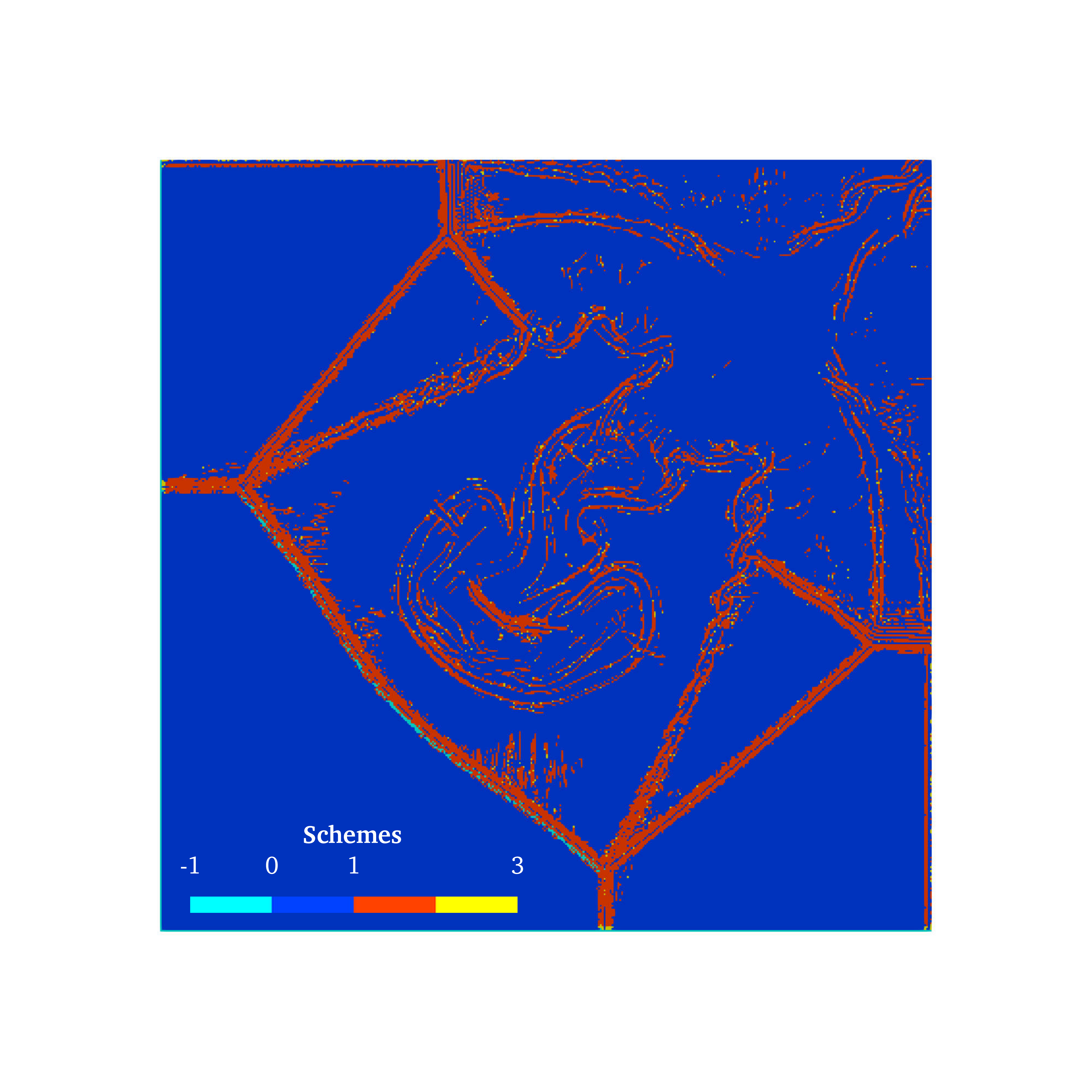}}

 \par\end{centering}\caption{Density contours for CWENOZ scheme \textcolor{black}{\textit({Top})}, and Cells using different schemes \textcolor{black}{(\textit{Bottom})}. From \textcolor{black}{{left to right}}, several configurations are presented in sequence. First, the default setting; then the hybrid result of the MUSCL and the CWENOZ schemes; then the hybrid result of the MUSCL and the \textcolor{black}{{high-order}} linear schemes; and finally, the hybrid configuration of CWENOZ and the \textcolor{black}{{high-order}} linear scheme.}\label{fig:RiemannDifferentHybrid}\end{figure}

 \clearpage
\subsection{iLES of supersonic viscous Taylor-Green vortex}

The supersonic variant of the Taylor-Green vortex flow problem is used for assessing the non-oscillatory properties and low-numerical dissipation of the developed framework.  This variant was introduced by Lusher and Sandham \citep{tgv_super}, and further investigated by Chapelier et al. \citep{Chapelier2024} for benchmarking several high-order CFD codes and is considered a well-established  benchmark for numerical methods for compressible flow problems.

The computational domain is defined as $\Omega=[0,2\pi]^3$ with periodic boundary conditions, and the present formulation of the Taylor-Green vortex problem is initialised with the following velocity, pressure and density fields:
\begin{equation}
u(x,y,z,0)=\sin(kx)\cos(ky)\cos(kz),
\label{eq:tgveq1}\end{equation}
\begin{equation}
v(x,y,z,0)=-\cos(kx)\sin(ky)\cos(kz),
\label{eq:tgveq2}
\end{equation}
\begin{equation}
w(x,y,z,0)=0,
\label{eq:tgveq3}
\end{equation}
\begin{equation}
p(x,y,z,0)=\frac{1}{\gamma M_{ref}^{2}}+\frac{1}{16}[\cos(2z)+2]\cdot[\cos(2x)+\cos(2y)].
\label{eq:tgveq5}
\end{equation}
The density is computed from the equation of state and is given by 
\begin{equation}
\rho(x,y,z,0)=\frac{\gamma M_{ref}^{2}}{T(x,y,z,0)}
\label{eq:tgveq5}
\end{equation}
with $T(x,y,z,0)=1$. 
We are using the variant with a Mach number $M_{ref}=1.25$, $Re=1600$ and the present viscous simulations were carried out on a hexahedral mesh of $128$ cells per side (2097152 total number of cells). We employ two 6th-order schemes, the CWENOZ  and the Hybrid scheme developed in this work. A CFL number of $1.4$ is used for the explicit Runge-Kutta 4th-order scheme, the HLLC Riemann solver \citep{HLLC}, the Low-Mach Correction of Simmonds et al. \citep{takis_LMC}, and a CWENOZ central stencil weight  $\lambda_1 = 10^{3}$ were employed, and the simulations were run up to $t=20$ for obtaining the required statistics. The results of Lusher and Sandham \cite{tgv_super} and the more recent comprehensive study of Chapelier et al. \cite{Chapelier2024} are used for comparison against the computed solutions.

The kinetic energy integrated over the domain is given by the following Eq. \ref{eq:tgvkinetic}

\begin{equation}
E_{K}^{\mathbf{t}}=\frac{1}{\rho_{ref}{\Omega}}\int_{\Omega}\frac{1}{2}\rho{u}_{j}{u}_{j}d{\Omega} ,
\label{eq:tgvkinetic}
\end{equation}

and the total viscous dissipation rate, which consists of the solenoidal dissipation $\varepsilon_{s}$ and dilatational dissipation $\varepsilon_{d}$ given by Eq. \ref{eq:tgv_vd} as

\begin{equation}
\textcolor{black}{\varepsilon_{t}=\varepsilon_{s}+\varepsilon_{d}=\frac{1}{\rho_{ref}{Re} {\Omega}}\int_{\Omega}\mu \left(\mathbf{\omega} \cdot \mathbf{\omega}\right)d{\Omega}+\frac{4}{3 \rho_{ref}{Re} {\Omega}}\int_{\Omega}\mu \left(\nabla \cdot \mathbf{u} \right)^{2} d{\Omega}  }.
\label{eq:tgv_vd}
\end{equation}
\textcolor{black}{Where $\mathbf{\omega}$ is the vorticity vector and $E_{K}=\left({E_{K}^t}/ {E_{K}^{0}}\right)$ is the normalised  kinetic energy}.

\begin{figure}[H]\begin{centering}
\captionsetup[subfigure]{width=0.20\textwidth}
{\includegraphics[angle=0,width=0.33\textwidth]{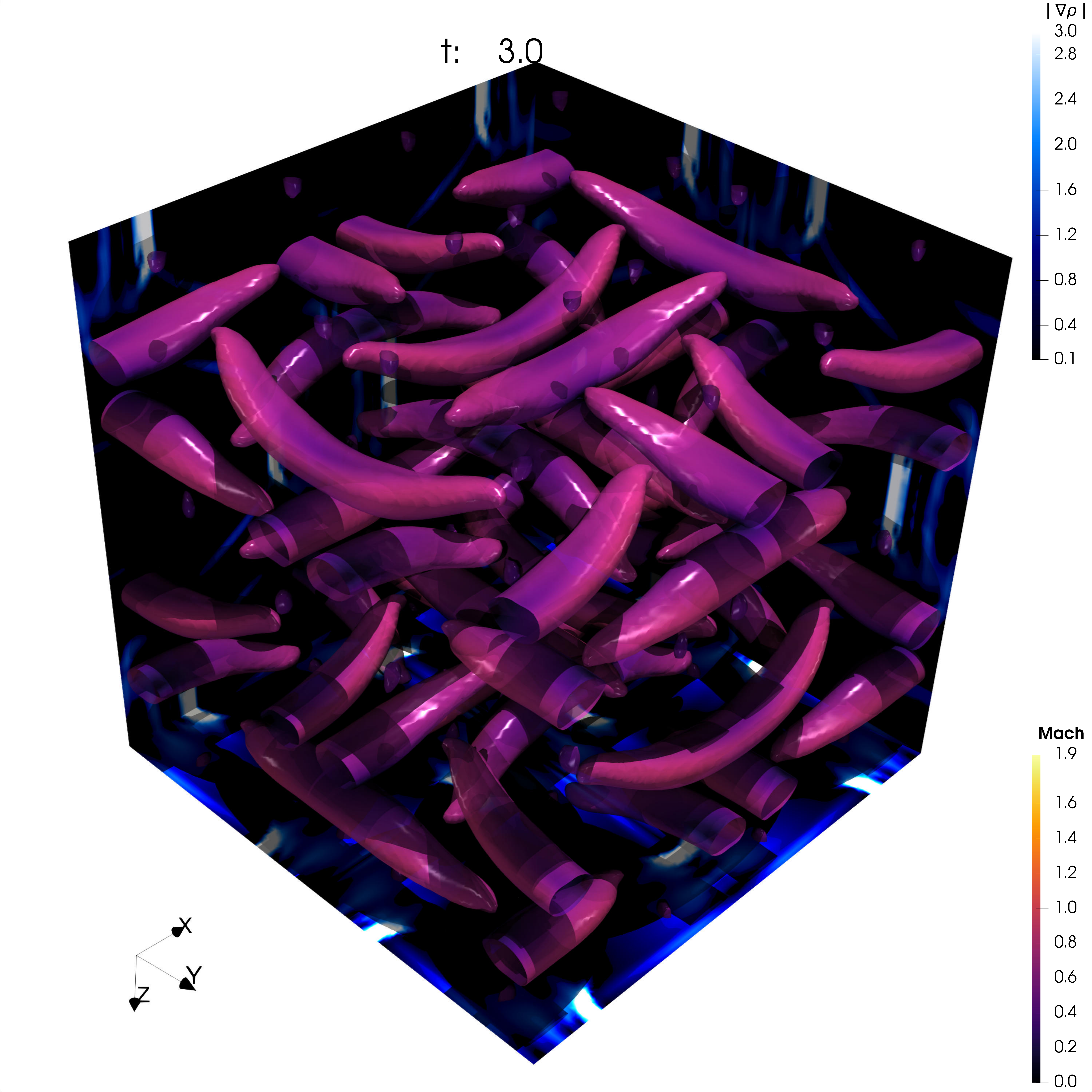}}
{\includegraphics[angle=0,width=0.33\textwidth]{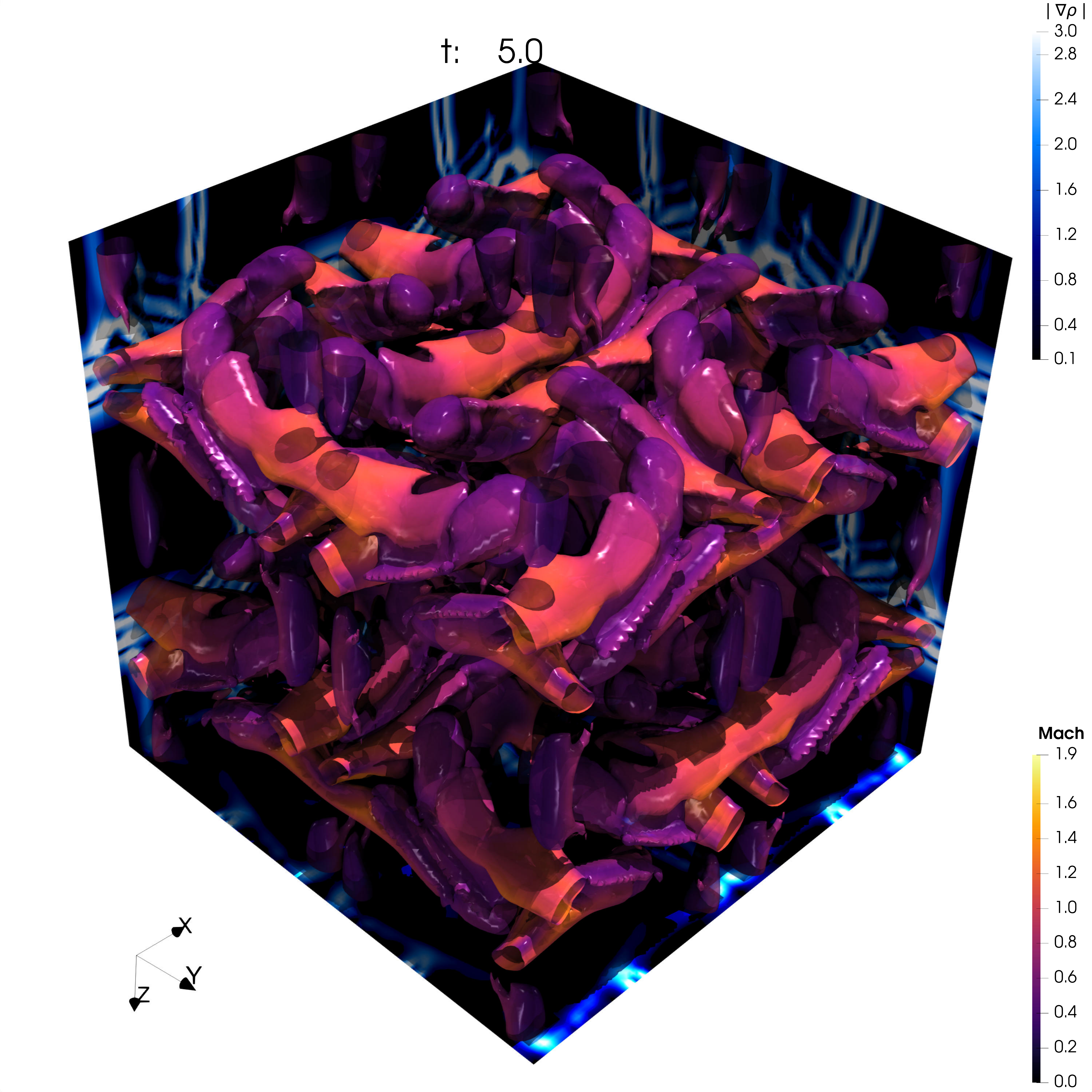}}
{\includegraphics[angle=0,width=0.33\textwidth]{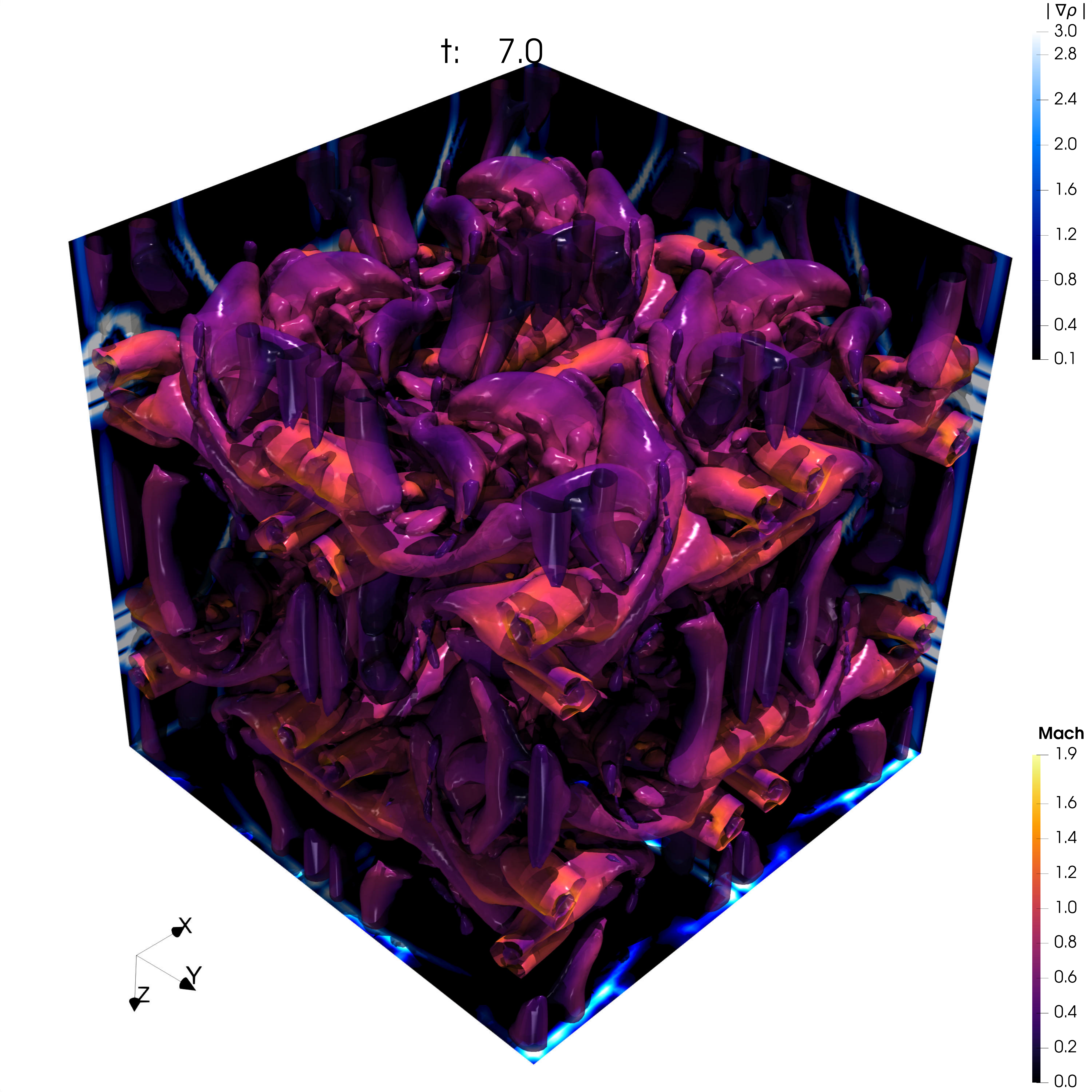}}\\
{\includegraphics[angle=0,width=0.33\textwidth]{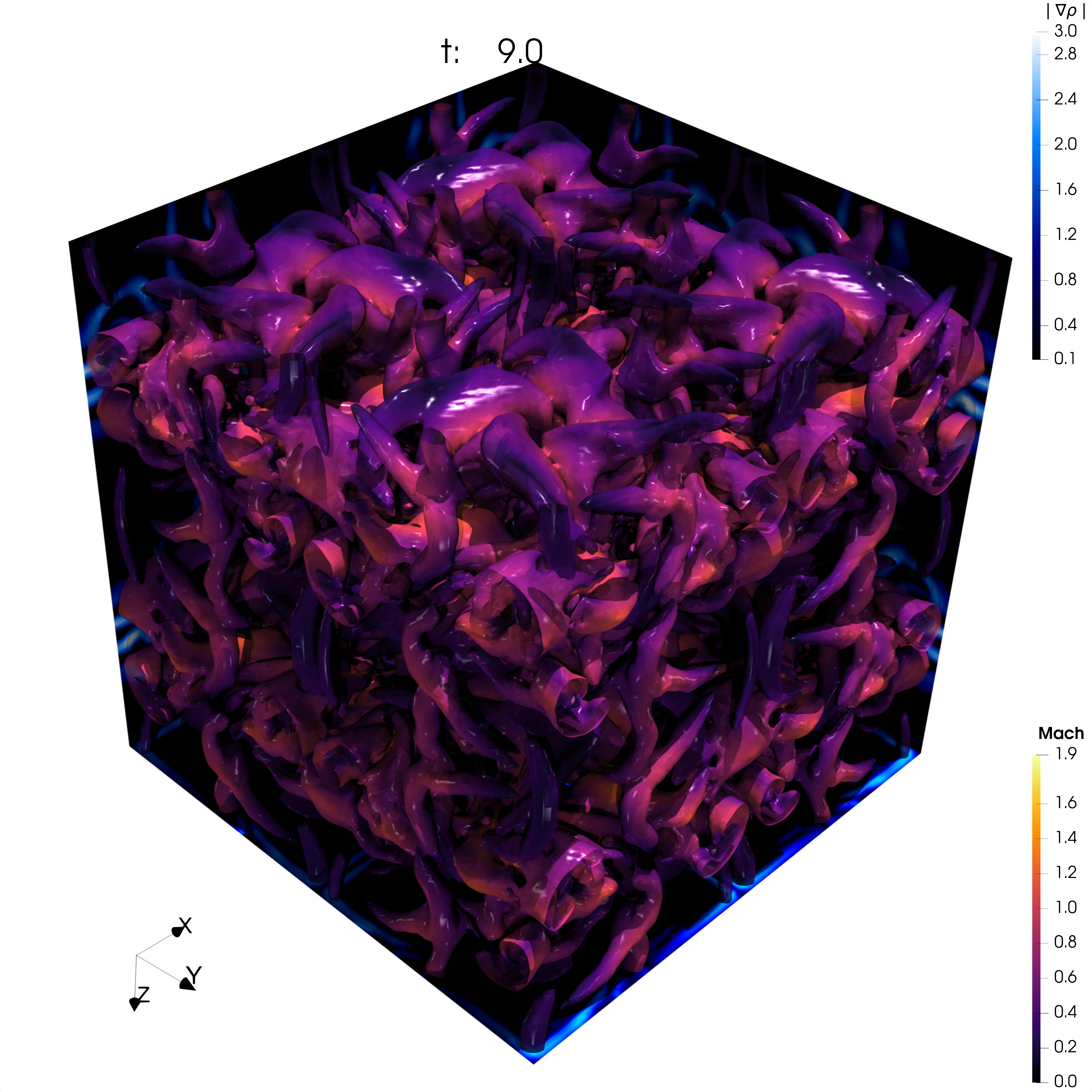}}
{\includegraphics[angle=0,width=0.33\textwidth]{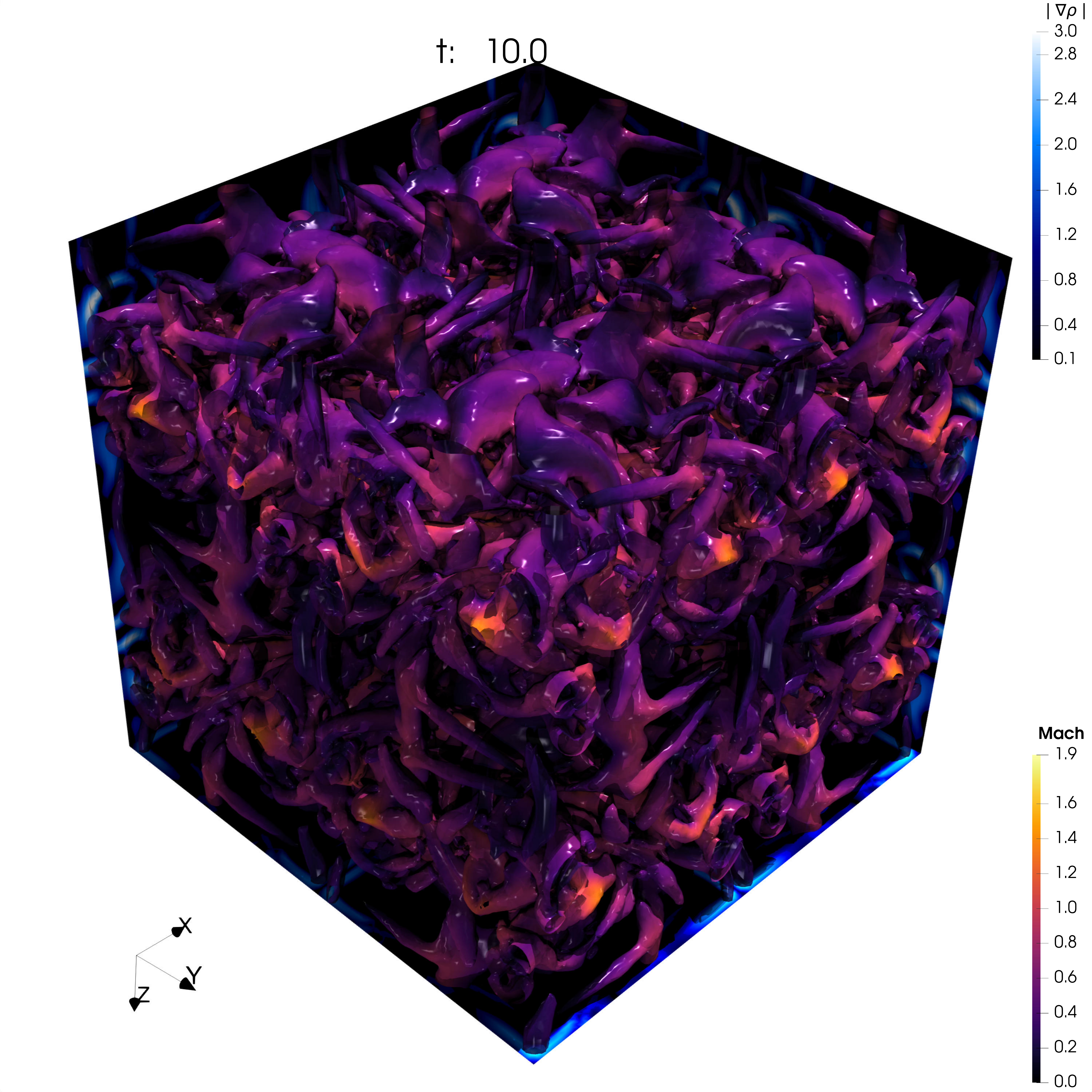}}
{\includegraphics[angle=0,width=0.33\textwidth]{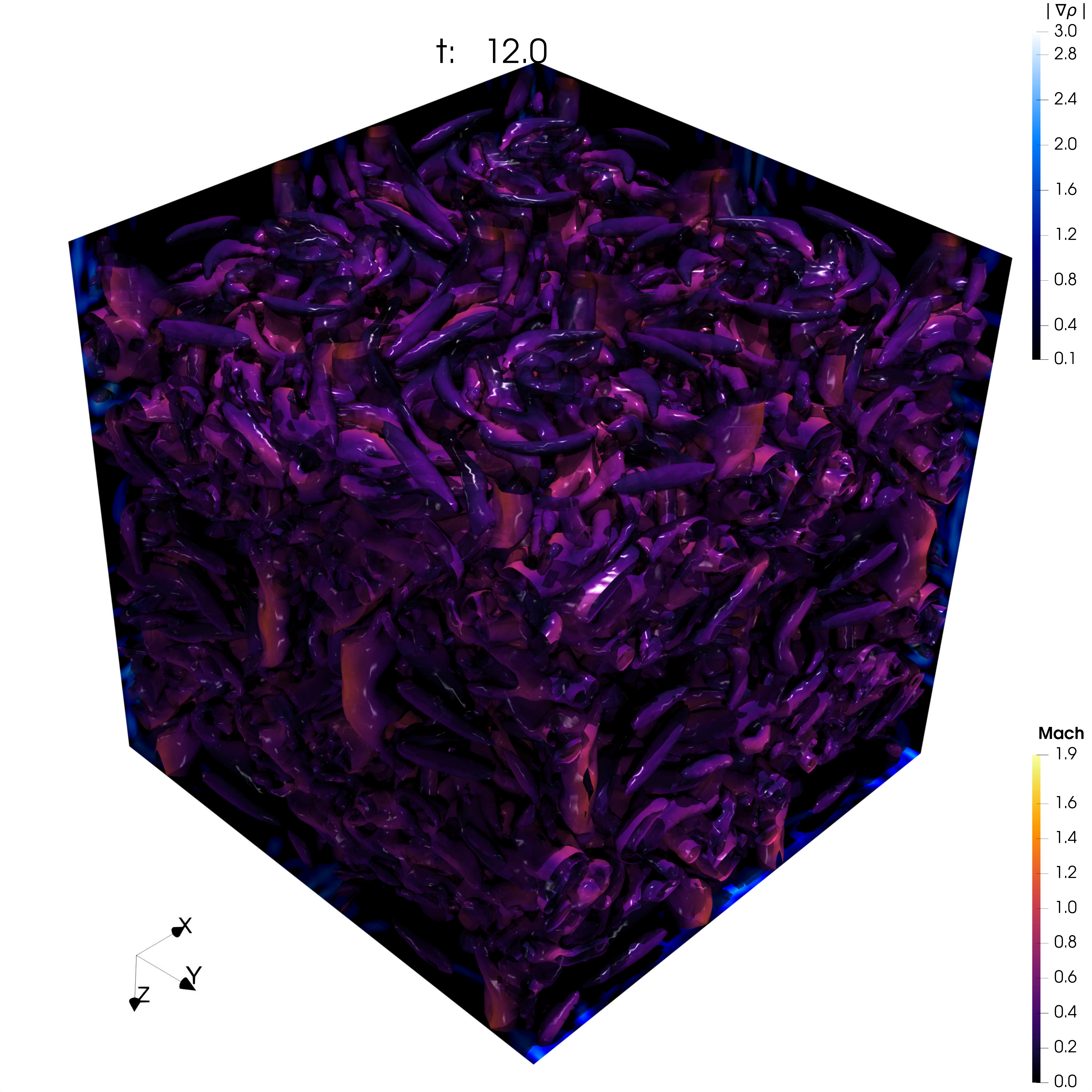}}
\par\end{centering}\caption{Contour plots at  different instants of the $|\nabla \rho|$ for the supersonic viscous Taylor-Green vortex flow at the $Y=\pi$ and $Z=\pi$ plane, and iso-surfaces of Q-criterion coloured by Mach number computed with the Hybrid 6th-order scheme on a hexahedral mesh of $128^3$.}\label{fig:TGVS2}\end{figure}

The evolution of the isosurface of Q criterion, coloured by the Mach number and slices of $|\nabla \rho|$ can be seen in Fig. \ref{fig:TGVS2}, where the representative transition to turbulence of the Taylor-Green vortex flow problem is captured.

\begin{figure}[H]\begin{centering}
\captionsetup[subfigure]{width=0.39\textwidth}
\subfloat[Kinetic energy]
{\includegraphics[angle=0,width=0.49\textwidth]{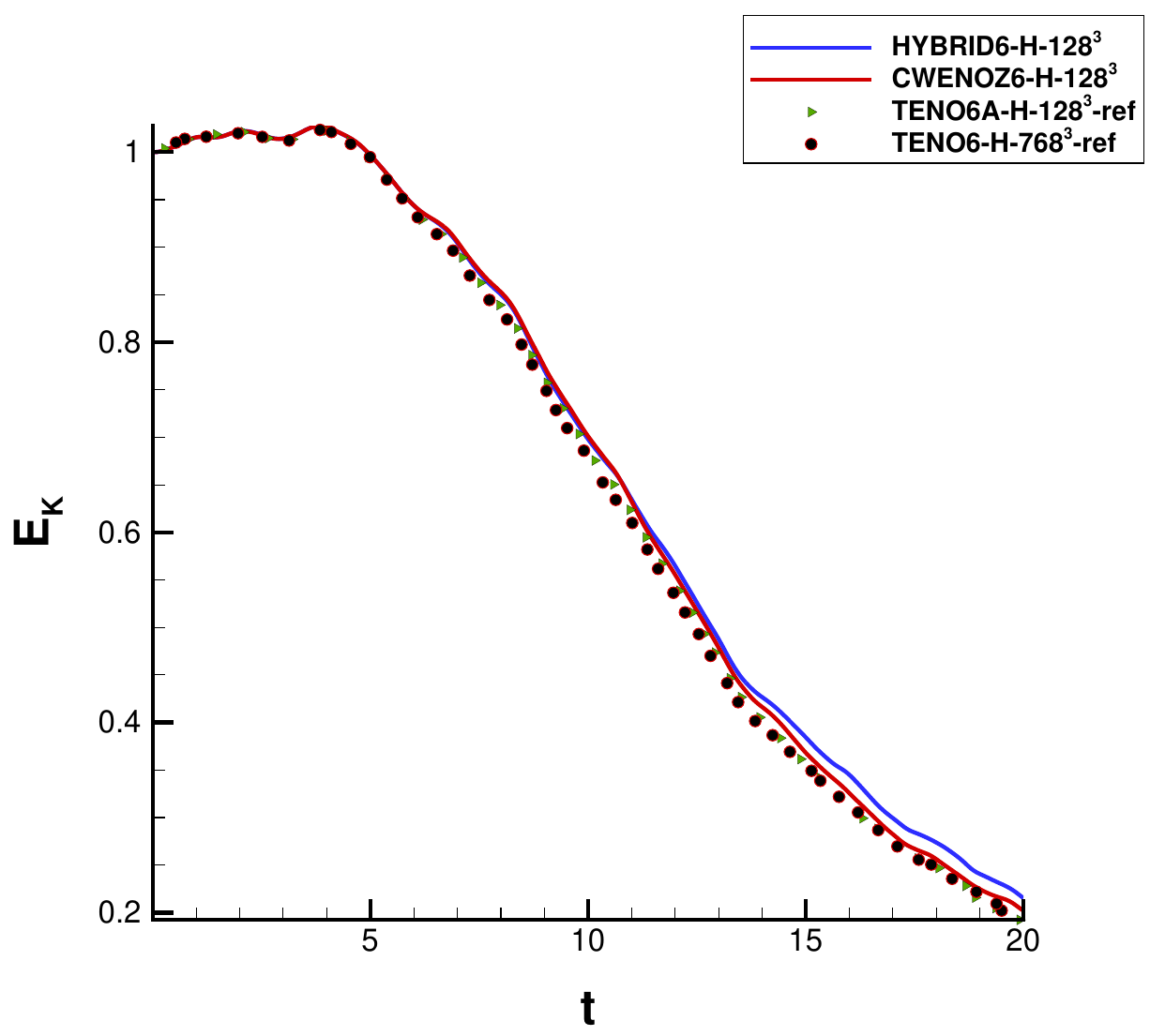}}
\subfloat[Total dissipation]
{\includegraphics[angle=0,width=0.49\textwidth]{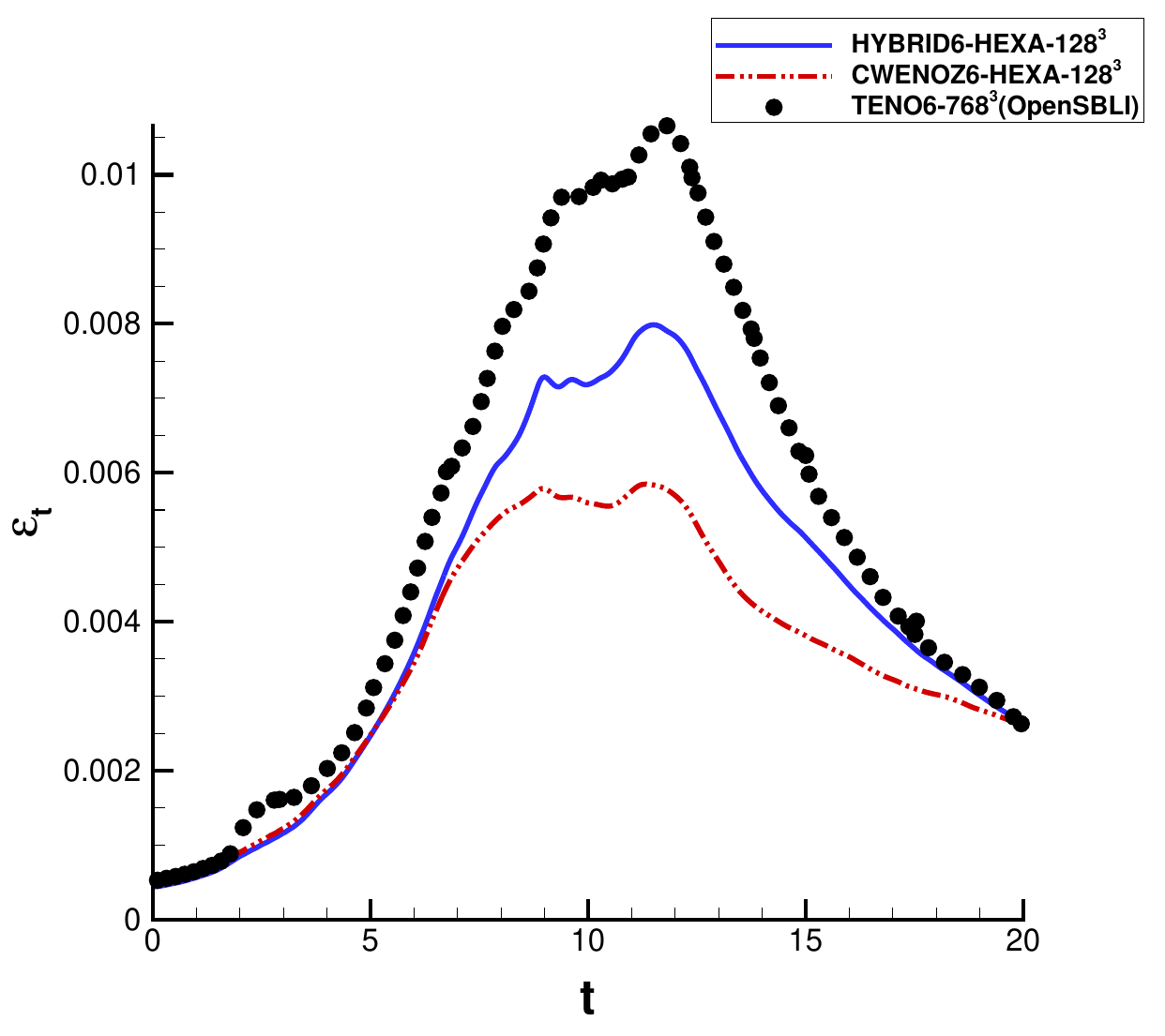}}\\
\subfloat[Solenoidal dissipation]
{\includegraphics[angle=0,width=0.49\textwidth]{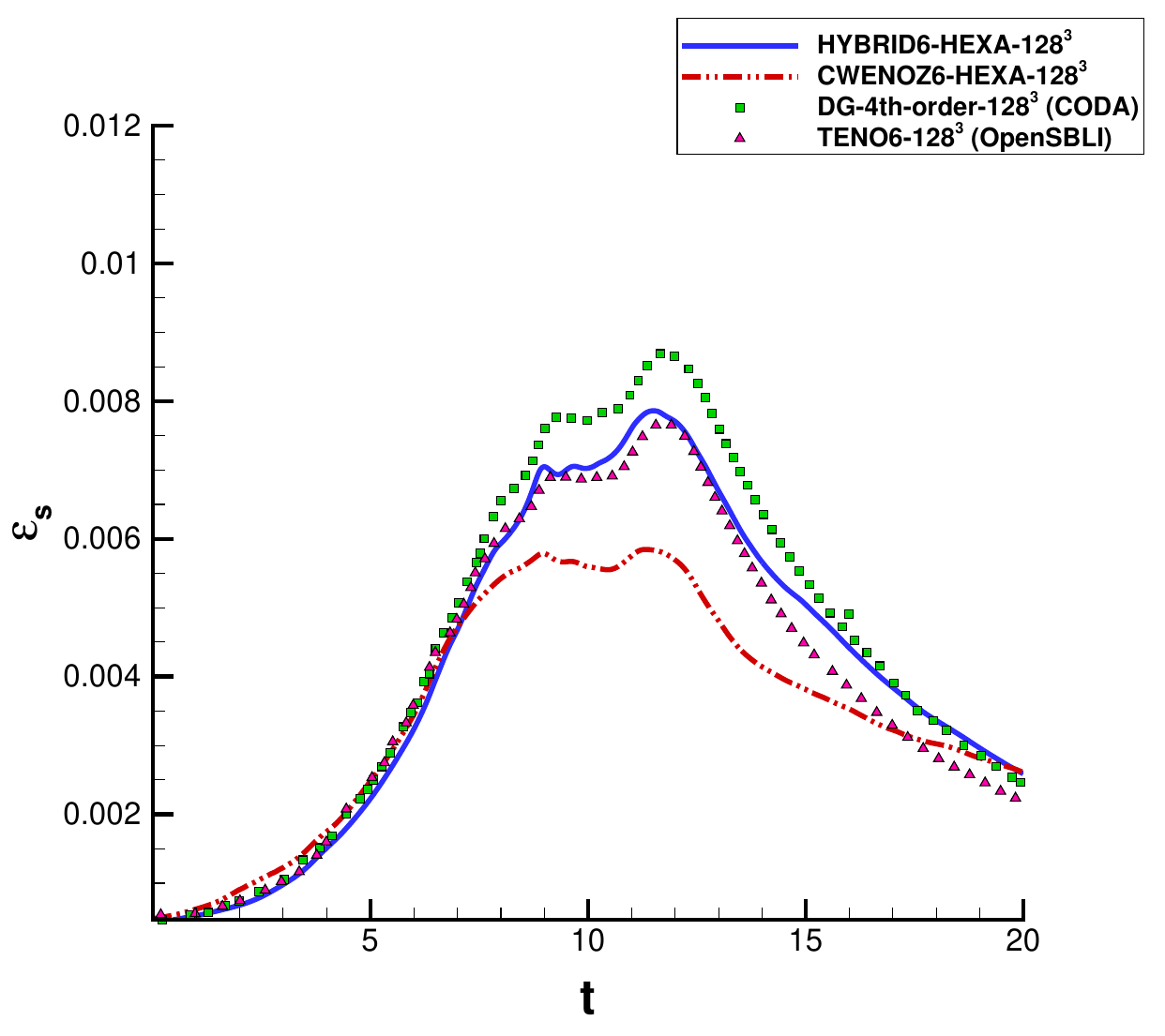}}
\subfloat[Dilatational dissipation]
{\includegraphics[angle=0,width=0.49\textwidth]{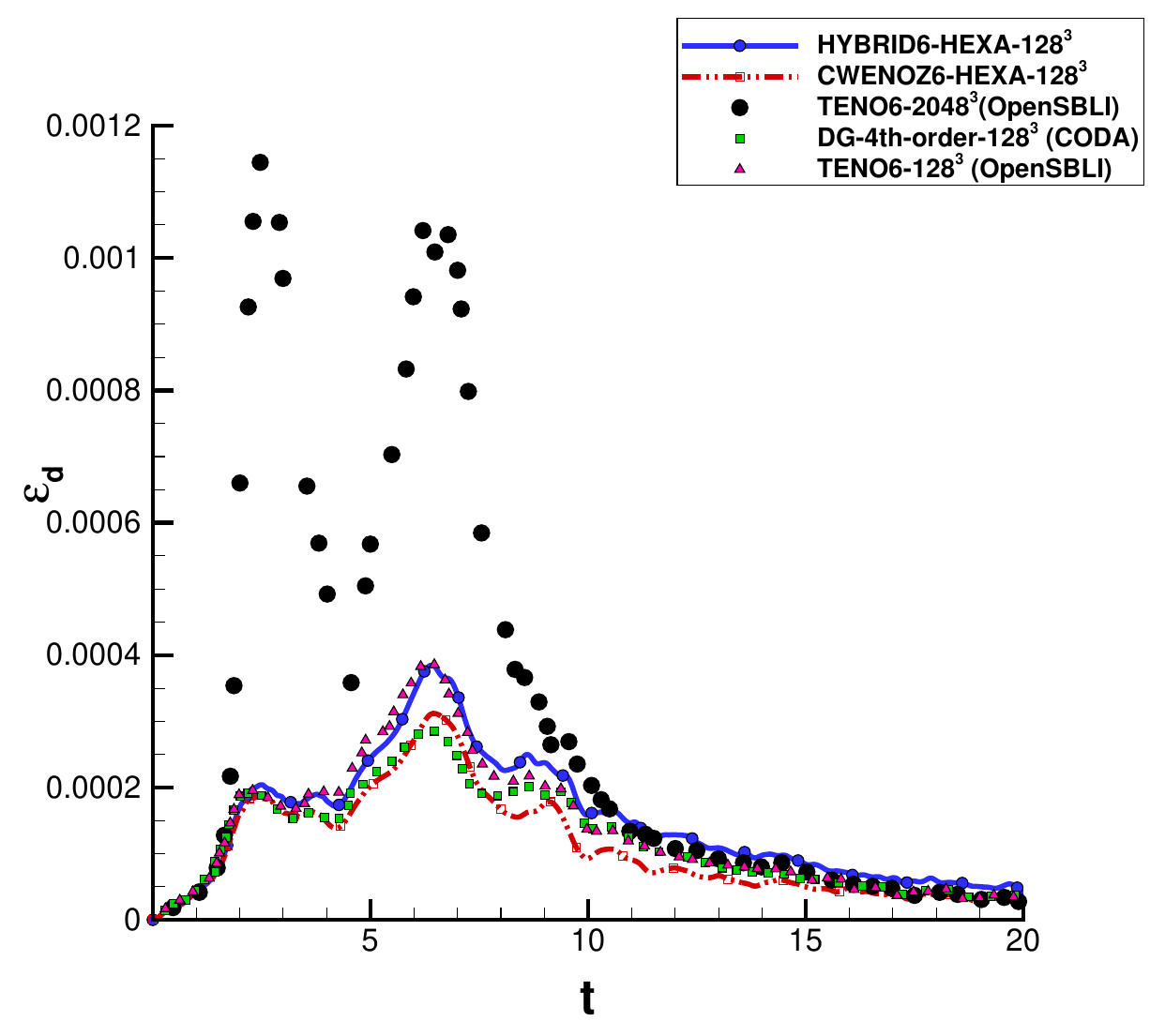}}
\par\end{centering}\caption{\textcolor{black}{Plots of the evolution of the normalised kinetic energy, total dissipation $\varepsilon_{t}$, solenoidal dissipation $\varepsilon_{s}$ and dilatational dissipation $\varepsilon_{t}$ for the solution of the supersonic viscous Taylor-Green vortex 
computed with several schemes on a hexahedral mesh  $128^3$.  It can be seen that the Hybrid6 scheme improves significantly the results compared to the CWENOZ6 scheme, and offers a good agreement with the modal DG-4th order scheme of the (CODA) CFD software and a 6th-order TENO scheme of the (OPENSBLI) CFD software as provided by Chapelier et al. \citep{Chapelier2024} at the same $(128^3)$ resolution and the results of Lusher and Sandham \cite{tgv_super}.}}\label{fig:TGVS1}\end{figure}

From the evolution of the kinetic energy as seen in Fig. \ref{fig:TGVS1}, all the schemes provide similar results and are in good agreement with the results obtained with  very-low dissipation TENO 6th-order schemes of Lusher and Sandham \cite{tgv_super,fu2019improved} at $(128^3)$ and $(768^3)$ resolutions.
From the total dissipation $\varepsilon_{t}$ evolution we notice that the Hybrid scheme is outperforming CWENOZ by  generating higher total-dissipation due to the wider range of scales being resolved.
At the $(128^3)$ resolution employed, a sudden increase and decrease of total dissipation is missing, driven by the dilatational component, occurring between (t = 2-4) seen at very fine resolutions $(768^3)$.
Finally we compare the obtained results in terms of the solenoidal dissipation $\varepsilon_{s}$ and dilatational dissipation $\varepsilon_{d}$ of our methods against the results obtained with a modal DG-4th order scheme of the (CODA) CFD software and a 6th-order TENO scheme of the (OPENSBLI) CFD software as provided by Chapelier et al. \citep{Chapelier2024} at the same $(128^3)$ resolution.
It can be seen that the Hybrid scheme has a good agreement with the (CODA) and the (OPENSBLI) results for both the solenoidal dissipation $\varepsilon_{s}$ and dilatational dissipation $\varepsilon_{d}$, and it offers a substantial improvement over the CWENOZ6 scheme.

 \begin{figure}[H]\begin{centering}
{\includegraphics[angle=0,width=0.49\textwidth]{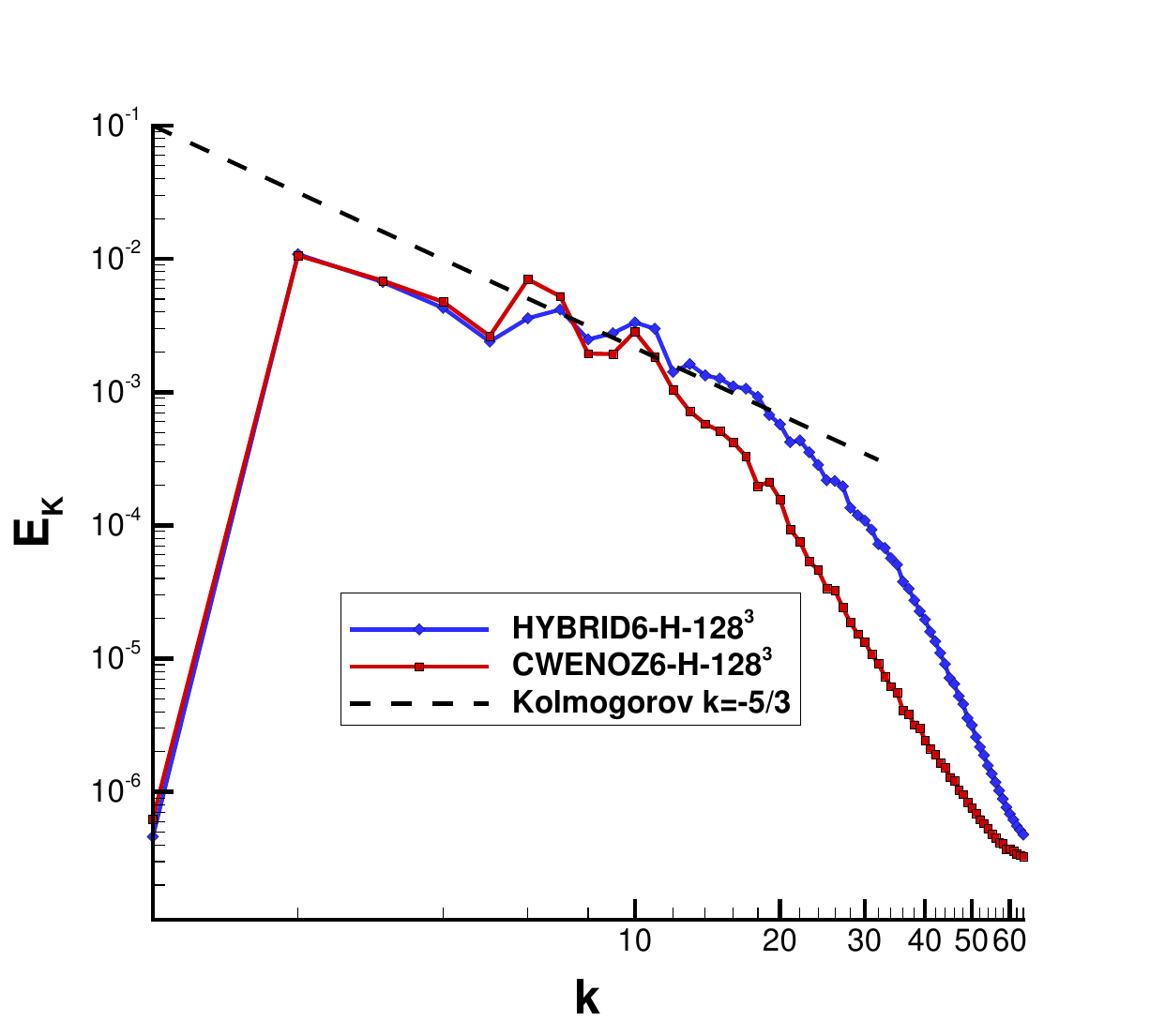}}
 {\includegraphics[angle=0,width=0.49\textwidth]{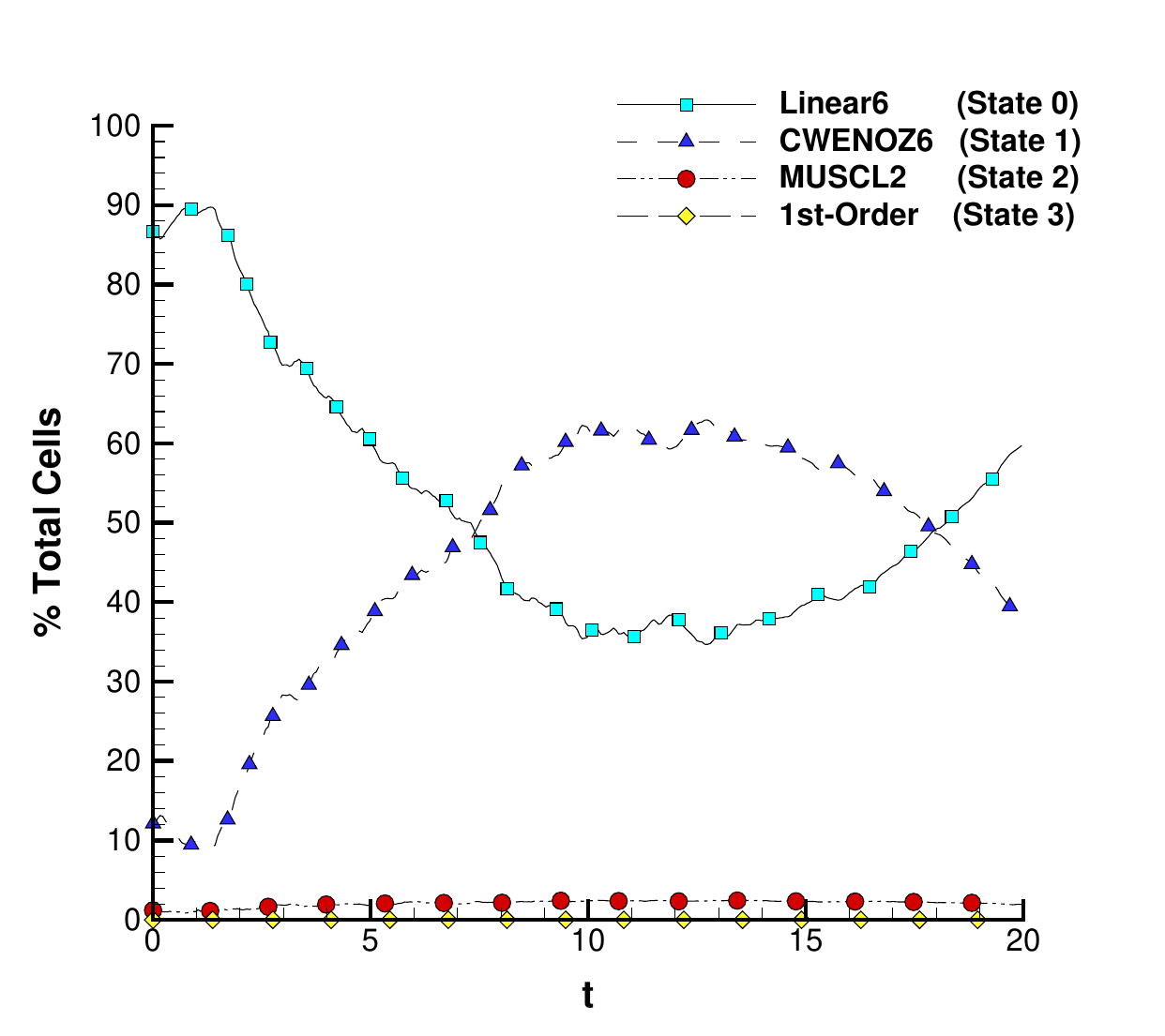}}\\
 \par\end{centering}\caption{Kinetic energy 3D spectra at $t=12$ for the solution of the supersonic viscous Taylor-Green (left) and time history of the evolution of the state of the cells in the domain (right). It  can be seen that the newly developed Hybrid6 scheme resolves a wider-range of wavenumbers compared to the CWENOZ6 scheme. The computational savings are driven by the reduced portion of the cells requiring an expensive CWENOZ reconstruction compared to the lower-cost linear 6th-order scheme, as seen in the evolution of the state of the cells.}\label{fig:TGVS_SPECTRA}\end{figure}

From the kinetic energy spectra immediately after the total dissipation peak at $t=12$ as seen in Fig. \ref{fig:TGVS_SPECTRA}, it can be noticed that the newly developed Hybrid6 scheme resolves a wider-range of wave numbers compared to CWENOZ6 and provides a better agreement with the Kolmogorov slope.  

\begin{figure}[H]\begin{centering}
\captionsetup[subfigure]{width=0.20\textwidth}
{\includegraphics[angle=0,width=0.33\textwidth]{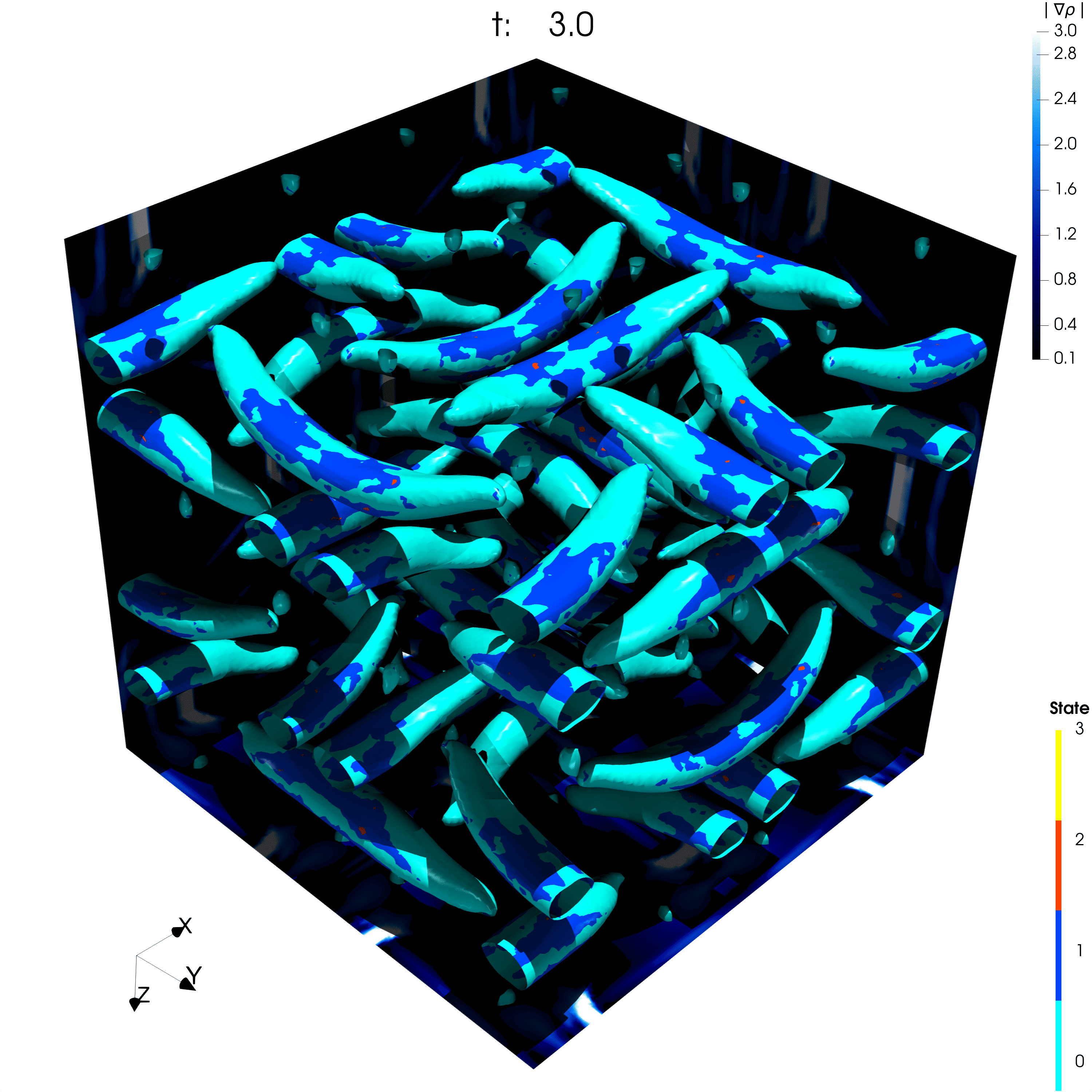}}
{\includegraphics[angle=0,width=0.33\textwidth]{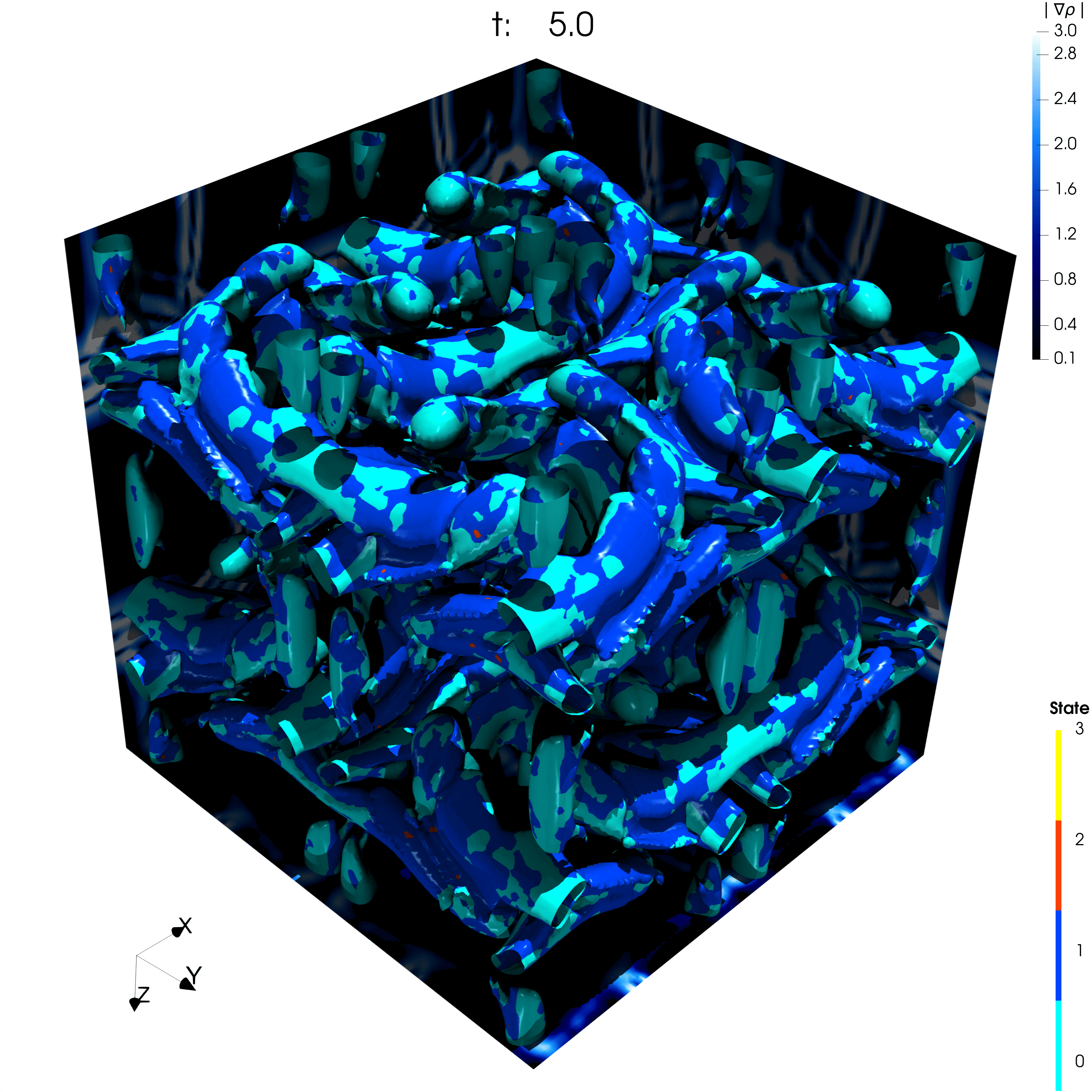}}
{\includegraphics[angle=0,width=0.33\textwidth]{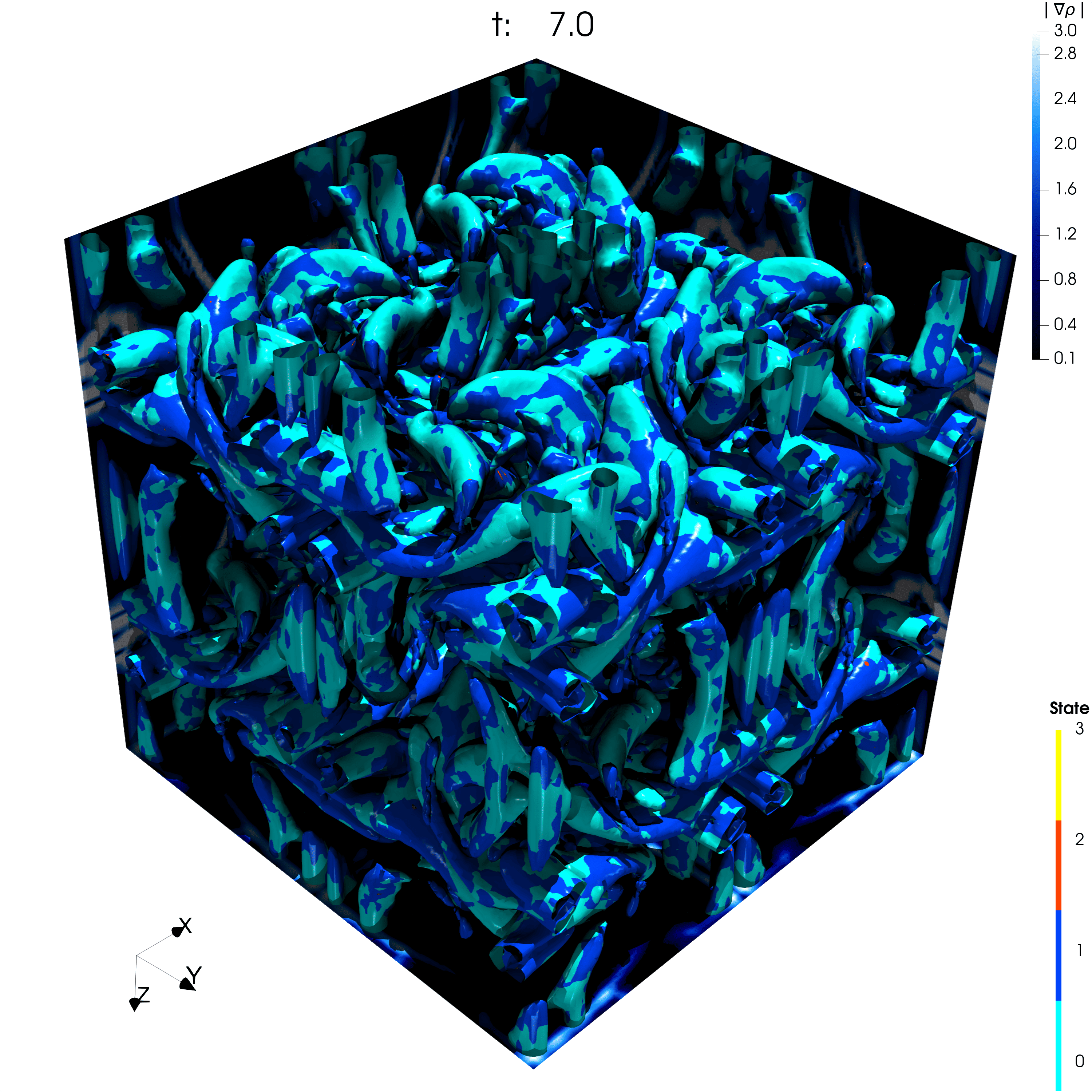}}\\
{\includegraphics[angle=0,width=0.33\textwidth]{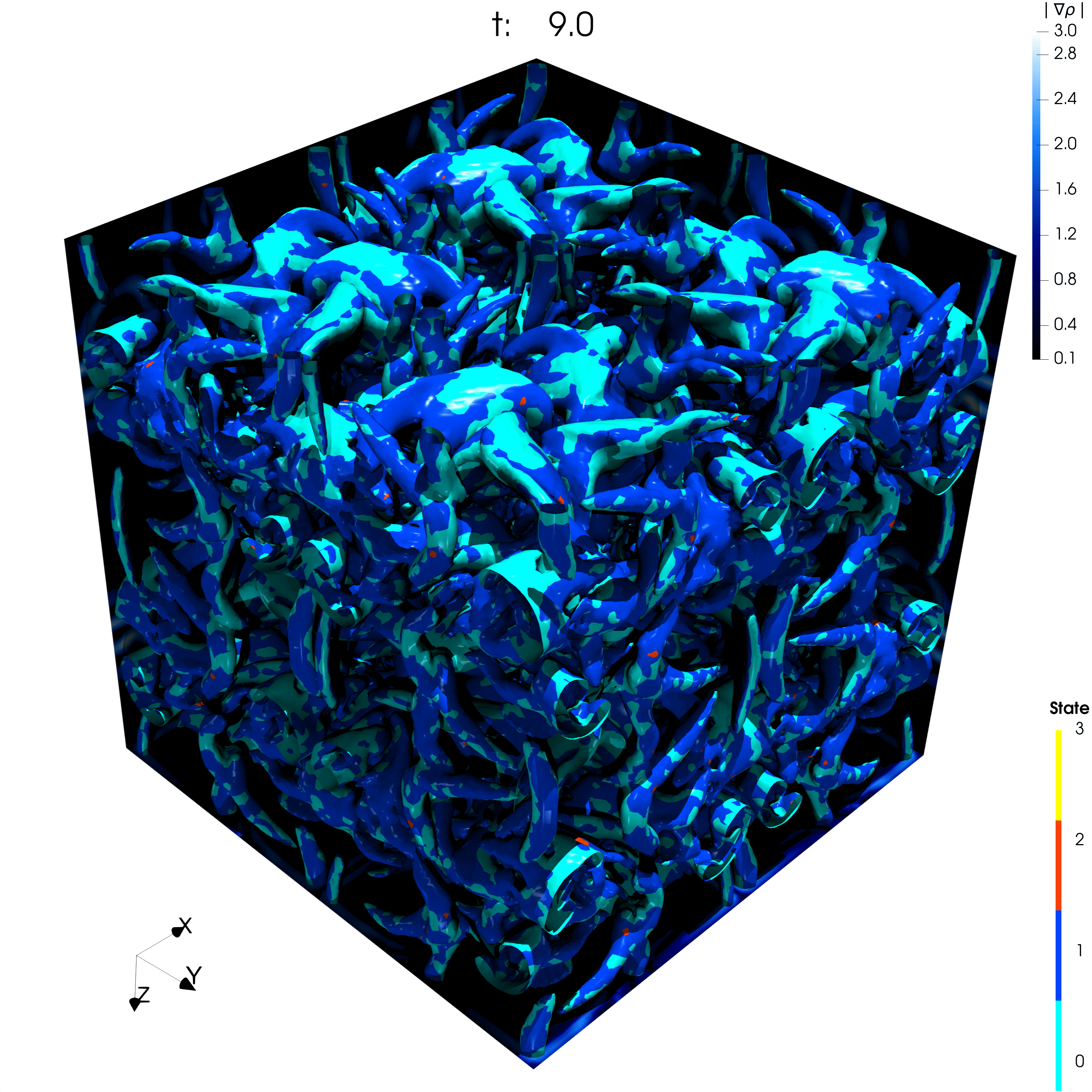}}
{\includegraphics[angle=0,width=0.33\textwidth]{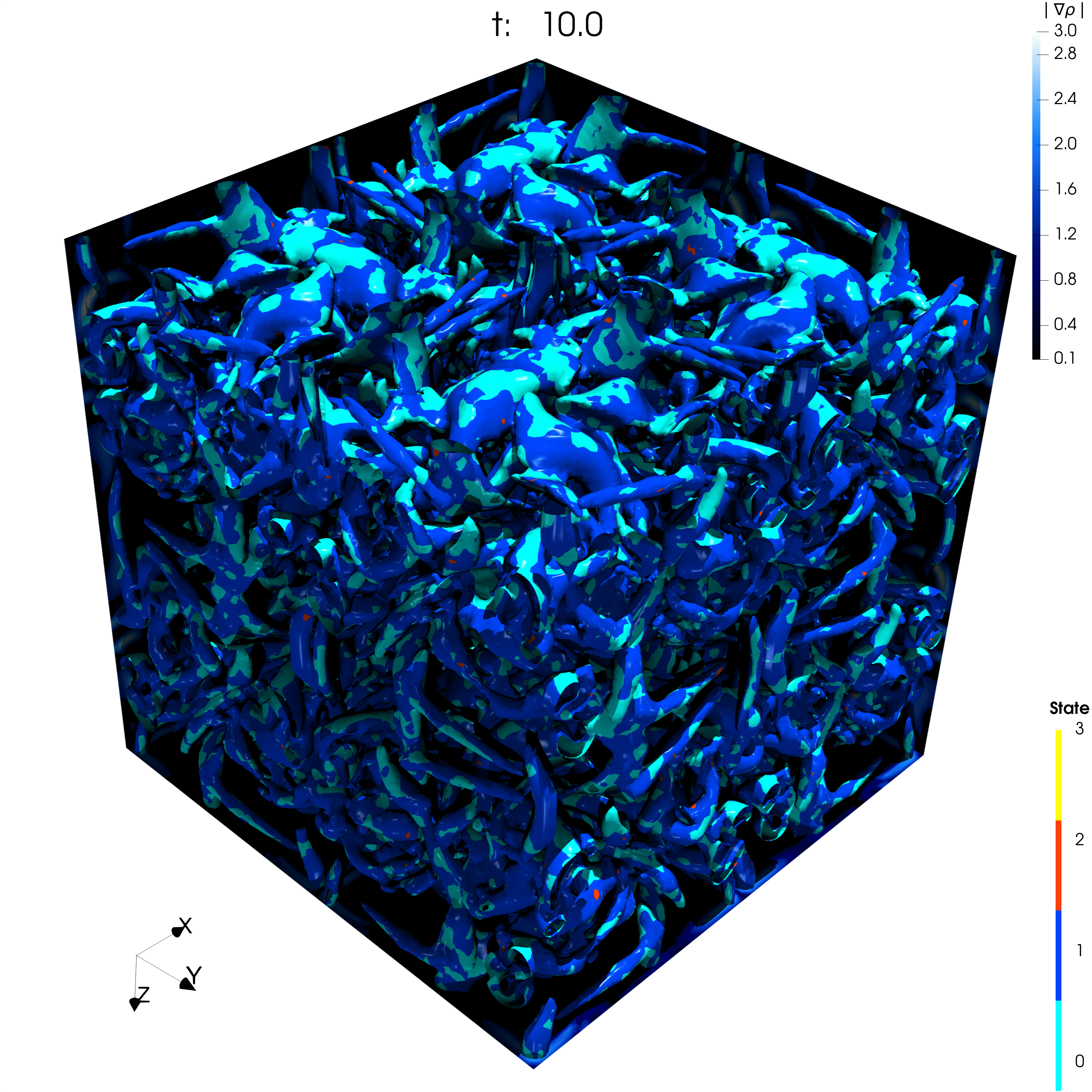}}
{\includegraphics[angle=0,width=0.33\textwidth]{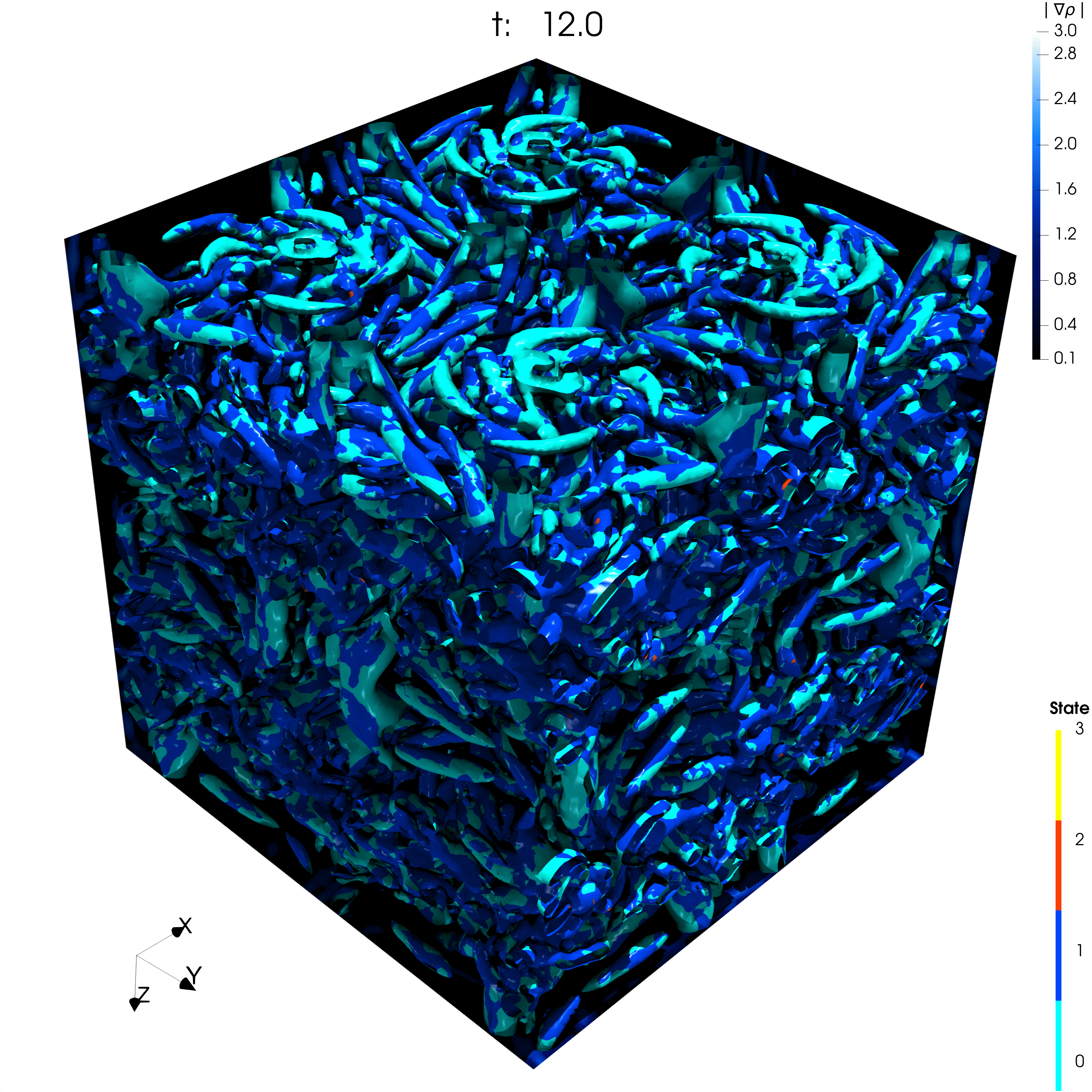}}
\par\end{centering}\caption{Contour plots at  different instants of the $|\nabla \rho|$ for the supersonic viscous Taylor-Green vortex flow at the $Y=\pi$ and $Z=\pi$ plane, and iso-surfaces of Q-criterion coloured by the state of each cell, computed with the Hybrid 6th-order scheme on a hexahedral mesh of $128^3$. The states of each method are Linear6 (State 0), CWENOZ6 (State 1), MUSCL2 (State 2), 1st-Order (State 3). It can be noticed that Linear6 dominates these vortical structures at early times (top) and as we are approaching the dissipation peak switching to CWENOZ6 method is required for ensuring robust non-oscillatory solutions.}\label{fig:TGVS3}\end{figure}

The Hybrid6 scheme is {2.5 times faster than the CWENOZ6 scheme} when using ARCHER2 and 10 nodes (128cores per node), due to the computational savings produced from the less-frequent deployment of the CWENOZ6 method. All the cells are solved with an unlimited Linear 6th-order method, and the cells that do not satisfy certain criteria previously described, switch to CWENOZ6, MUSCL2 and even first-order method. \textcolor{black}{The increase in the number of cells treated with nonlinear dissipation near the peaks of the kinetic energy dissipation rate is an expected response of the hybrid scheme. As finer-scale structures emerge and the flow becomes locally under-resolved, the switching mechanism activates the CWENOZ reconstruction in a larger fraction of the domain to control high-frequency content generated by the unlimited high-order linear discretisation and to preserve numerical stability.
}
The computational savings are therefore mainly due to the additional reconstructions needed for the directional stencils of the CWENOZ method, and the compute intensive CWENOZ weights calculations. This can also be seen in Fig. \ref{fig:TGVS3} where Linear6 dominates these vortical structures at early times and as we are approaching the dissipation peak switching to CWENOZ6 method is required for ensuring robust non-oscillatory solutions.

\subsection{Mach 3 flow over cylinder}
Finally, the inviscid supersonic flow past a circular cylinder is simulated to assess the numerical performance of the proposed hybrid discretisation approach. The computational setup follows the benchmark configuration reported in \citep{DENG2023102150}. A circular cylinder of radius $0.25$ is positioned at $(0.6,\,1)$ within a two-dimensional wind tunnel spanning the domain $[0,4] \times [0,2]$. The computational mesh consists primarily of quadrilateral elements, with a total of $1\,947\,909$ cells, of which  $2\,438$ are triangular elements. A Mach~3 supersonic inflow condition is prescribed at the left boundary, while an outflow condition is applied at the right boundary. The upper and lower boundaries are treated as inviscid slip walls. The simulation is carried out using a hybrid sixth-order spatial discretisation scheme with default settings, the HLL Riemann solver and a CFL number of $0.7$. The computation is advanced up to a final time of $t = 5.0$.

\begin{figure}[H]\begin{centering}
\captionsetup[subfigure]{width=0.20\textwidth}
{\includegraphics[angle=0,width=0.6536\textwidth]{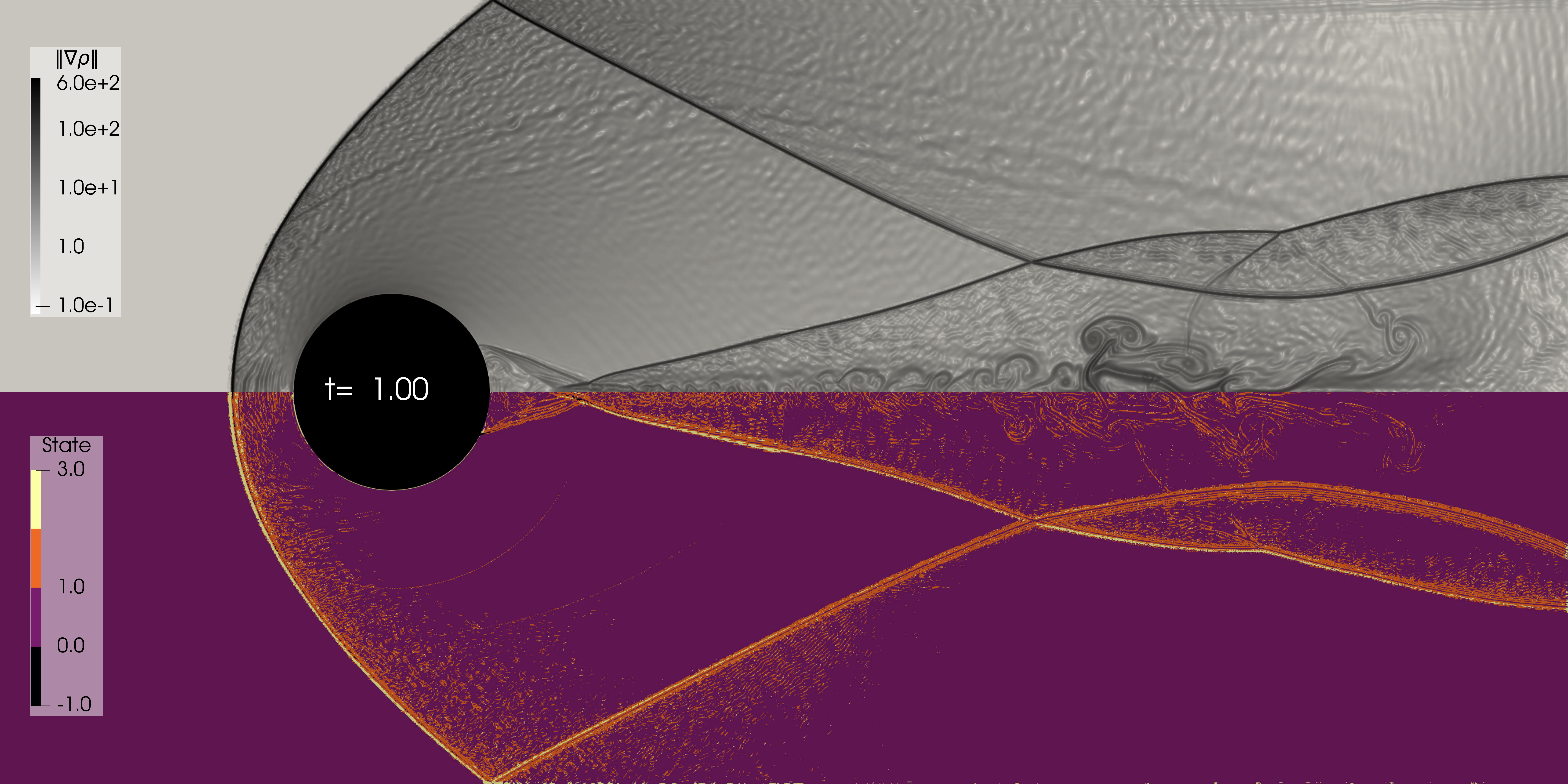}}
{\includegraphics[angle=0,width=0.3263\textwidth]{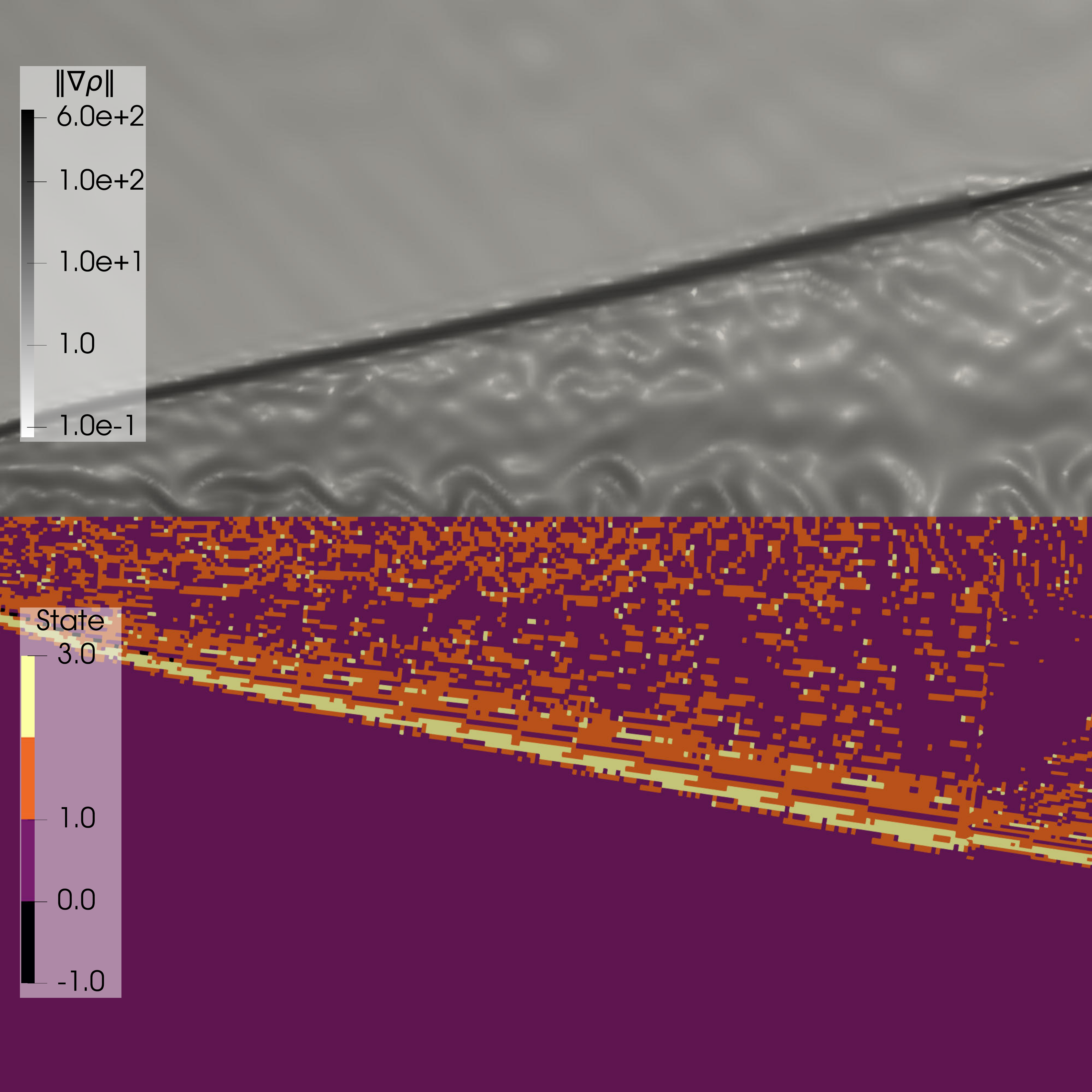}}
{\includegraphics[angle=0,width=0.6536\textwidth]{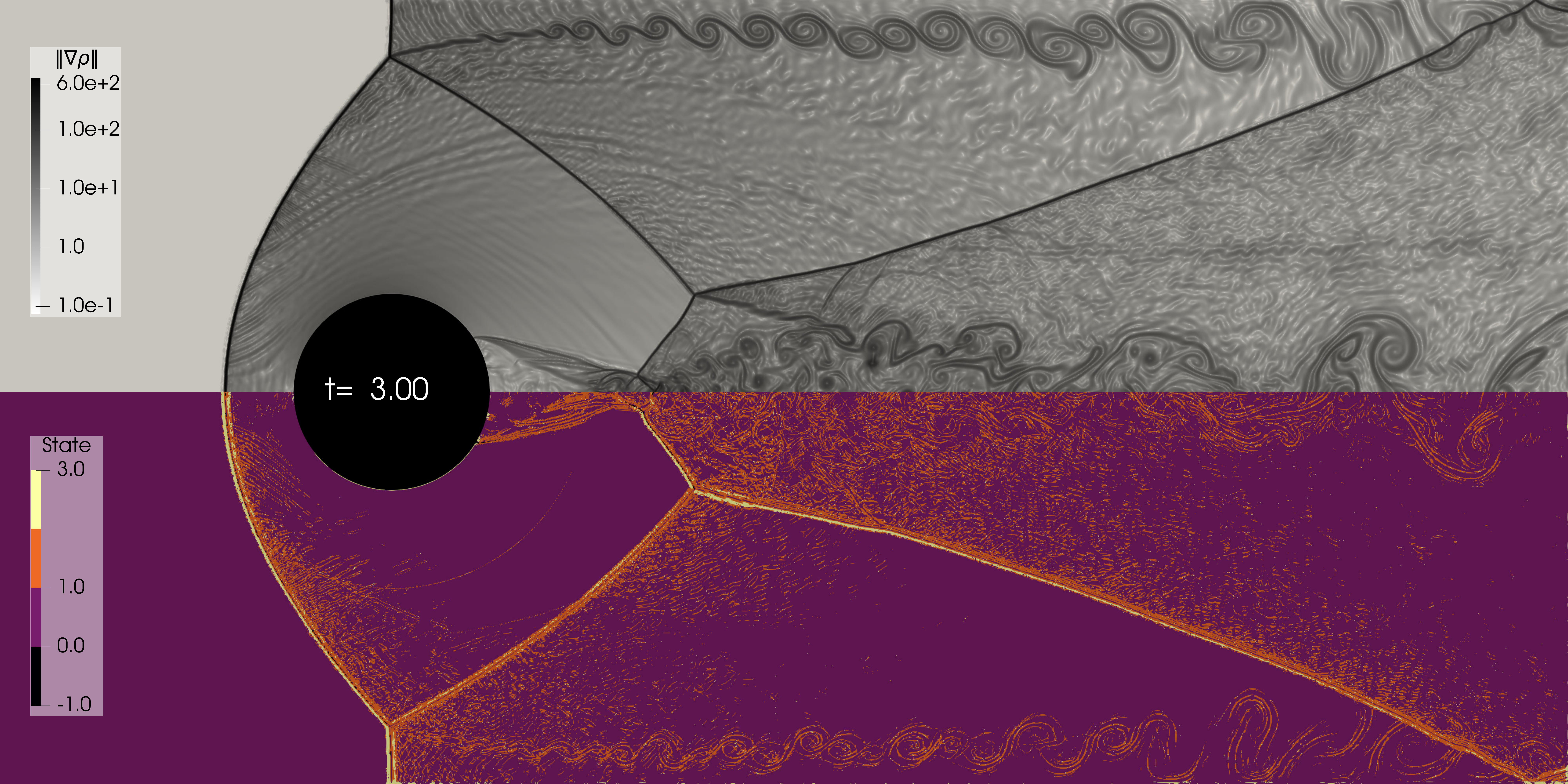}}
{\includegraphics[angle=0,width=0.3263\textwidth]{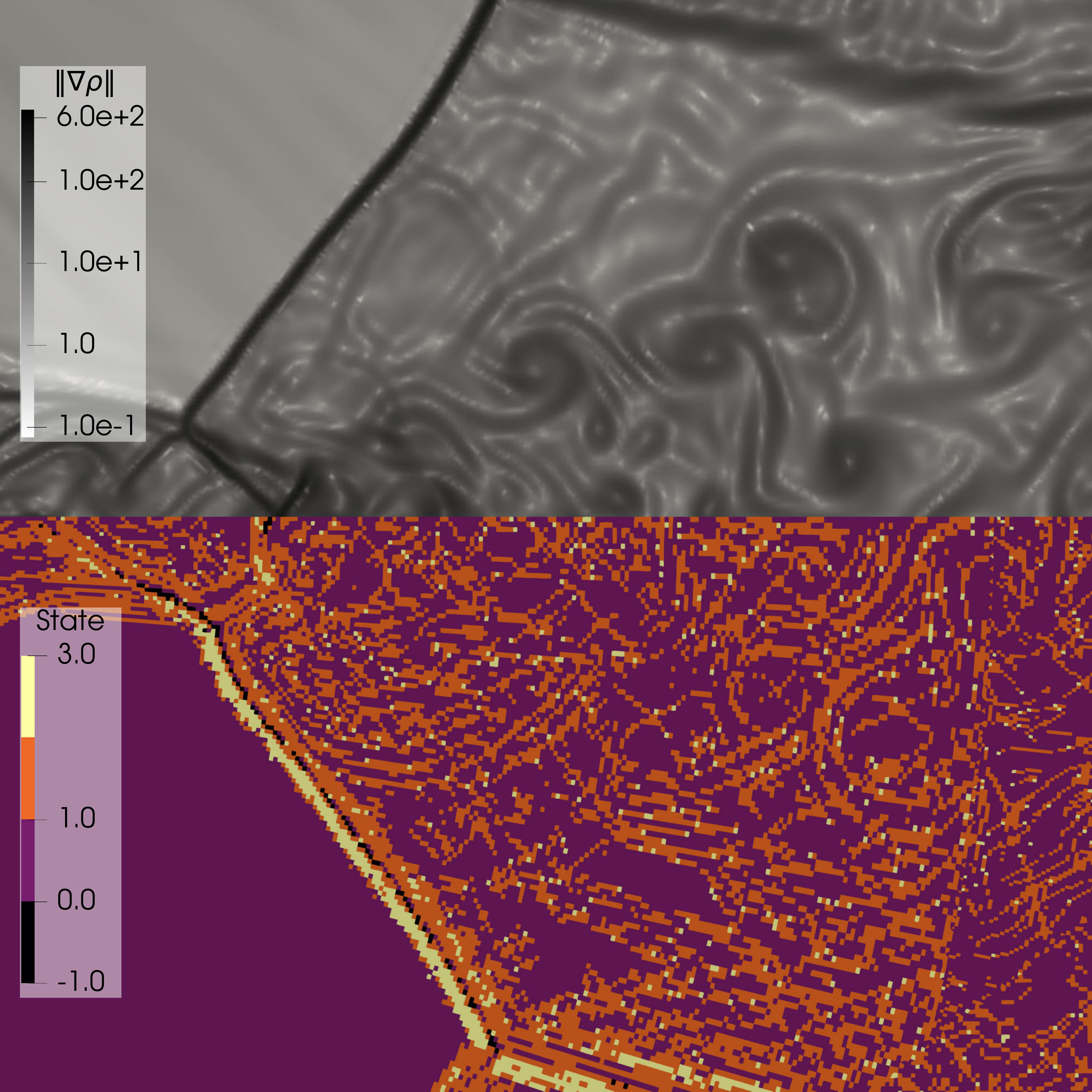}}
{\includegraphics[angle=0,width=0.6536\textwidth]{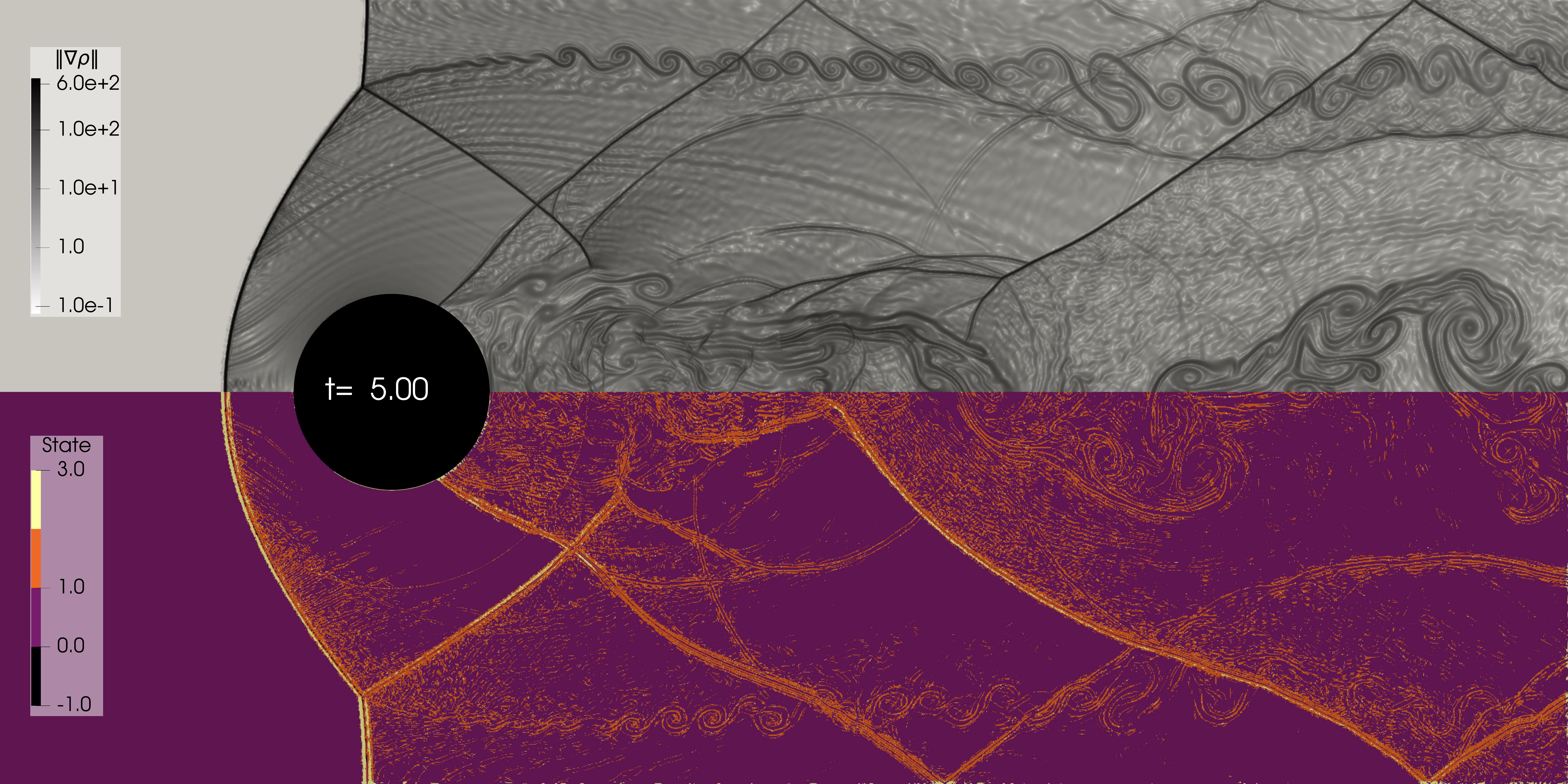}}
{\includegraphics[angle=0,width=0.3263\textwidth]{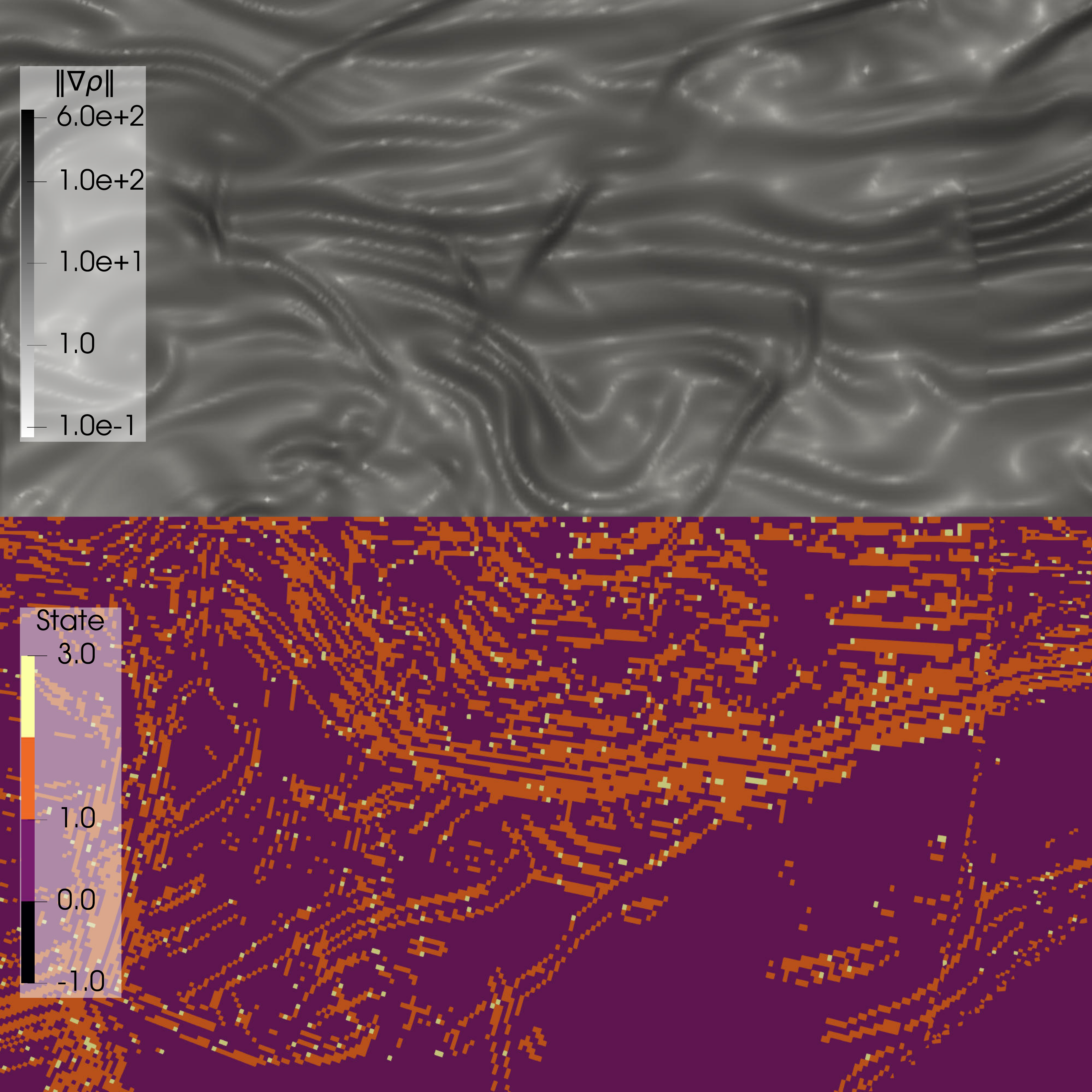}}
\par\end{centering}\caption{Contour plots at  different instants of the $|\nabla \rho|$ at the top half, and the state of each cell in the domain lower  half for the supersonic Mach 3 flow over a cylinder, computed with the Hybrid 6th-order scheme. It can be noticed that the high-order 6th-order linear scheme dominates the flow field (State 0 -purple colour), followed by the CWENOZ6 scheme (State 1 -orange colour), and the MUSCL2 (State 2 -yellow colour) is switched primarily on at the shocks, while the 1st-Order (State 3 -black colour) has been activated only across the strongest bow shock.}\label{fig:cylinder2}\end{figure}

The developed hybrid scheme captures the key flow features of the Mach~3 cylinder problem in a robust manner. In particular, the bow shock is sharply resolved without spurious oscillations, while the downstream shear layers and wake dynamics remain well defined. The method also successfully resolves the onset and growth of wake instabilities, indicating that the nonlinear switching does not introduce excessive numerical dissipation in smooth but under-resolved regions. As evidenced by the cell-state distribution in Fig. \ref{fig:cylinder2}, the sixth-order linear discretisation dominates over most of the domain, with a targeted transition to CWENOZ6 in regions of increased gradients and to MUSCL2 in the vicinity of shocks; the first-order fallback is activated only locally across the strongest portion of the bow shock. This selective blending provides a favourable balance between stability and accuracy: discontinuities are treated robustly, whereas high-order resolution is retained in smooth regions where small-scale structures develop.

From a computational-efficiency perspective, the hybrid approach delivers substantial savings relative to a uniform CWENO6 sixth-order discretisation which is approximately three times more expensive per iteration on 4 nodes on ARCHER2, using the optimum MPI/thread placement of 8 MPI processes per node and 16 OpenMP threads per MPI process.

\section{Conclusions}\label{sec:conclusions}
This work introduced a hybrid high-order finite-volume framework that integrates linear reconstruction in smooth regions, CWENOZ in weakly non-smooth zones, and MUSCL at true discontinuities, all coordinated by a redesigned single-step DMP-based detector. The detector effectively distinguishes smooth, weakly non-smooth, and discontinuous regions, minimizing over-flagging away from shocks and activating nonlinear reconstructions only when necessary. Across a suite of canonical 2D and 3D benchmarks, the method accurately resolves shocks and steep gradients while remaining free from spurious oscillations.
The hybrid strategy achieves significant computational savings relative to uniform CWENOZ, confining costly nonlinear operations to regions where they are essential. It maintains high accuracy across different mesh resolutions and orders of accuracy, with only a minor increase in dissipation at higher orders—readily mitigated through modest parameter retuning. The framework’s flexibility allows seamless switching between three-scheme and two-scheme modes, ensuring both robustness and adaptability. In particular, in the three-dimensional supersonic viscous Taylor–Green vortex, the method achieves the closest match to established high-order reference solutions, demonstrating its capability for compressible turbulence and complex flow dynamics. 
Looking ahead, this hybrid framework brings high-order accuracy closer to practical realisation in industrial-scale CFD simulations through its combination of reduced computational cost, improved robustness, and reliability. Building on this foundation, future work will extend the approach to complex three-dimensional engineering flows, where its balance of efficiency, sharp shock capturing, and resilience can have the greatest impact.

\section*{Acknowledgements}
The authors acknowledge the computing time on ARCHER2 through UK Turbulence Consortium [EP/X035484/1], and P.T. acknowledges the support provided by the EPSRC grant for ``Adaptively Tuned High-Order Unstructured Finite-Volume Methods for Turbulent Flows'' [EP/W037092/1].

\section*{Data Availability}
\added{The datasets of the test problems in this article will be available at the Cranfield Online Research Data repository}

\newpage
\bibliography{scopus}
\bibliographystyle{elsarticle-num}

\end{document}